\documentclass[a4paper]{article}
\usepackage{amssymb,amscd,amsfonts,mathrsfs,amsmath}
\usepackage[all]{xy}
\input amssym.def

\def\double{\mathbb}

\def\cc{{\double C}}
\def\nn{{\double N}}
\def\zz{{\double Z}}
\def\qq{{\double Q}}
\def\rr{{\double R}}

\newtheorem{theorem}{Theorem}[section]
\newtheorem{lemma}[theorem]{Lemma}
\newtheorem{corollary}[theorem]{Corollary}
\newtheorem{definition}[theorem]{Definition}
\newtheorem{proposition}[theorem]{Proposition}
\newtheorem{remark}[theorem]{Remark}

\newtheorem{example}[theorem]{Example}

\def\Kt{K^{\mathrm{top}}}
\def\Ka{K^{\mathrm{alg}}}

\def\limproj{\mathop{\mathrm{lim}}\limits_{\longleftarrow}}
\def\cp{\rtimes}
\def\si{\sigma}

\def\cinf{C^{\infty}}
\def\cinfc{C^{\infty}_c}

\newcommand{\be}{\begin{equation}}
\newcommand{\ee}{\end{equation}}
\newcommand{\beq}{\begin{eqnarray}}
\newcommand{\eeq}{\end{eqnarray}}
\newcommand{\om}{\omega}
\newcommand{\Om}{\Omega}
\newcommand{\al}{\alpha}
\def\nat{\natural}

\newcommand{\la}{\lambda}
\newcommand{\Ec}{{\mathscr E}}

\newcommand{\Lc}{{\mathscr L}}
\newcommand{\non}{\nonumber}
\newcommand{\eps}{\varepsilon}
\newcommand{\Sc}{{\mathscr S}}

\newcommand{\Rc}{{\mathscr R}}
\newcommand{\Mc}{{\mathscr M}}
\newcommand{\Nc}{{\mathscr N}}
\newcommand{\Ic}{{\mathscr I}}
\newcommand{\Jc}{{\mathscr J}}

\def\ch{\mathrm{ch}}
\def\ad{\mathrm{Ad\,}}
\def\tch{\mathrm{\,\slash\!\!\!\! \ch}}
\def\cs{\mathrm{cs}}

\def\chdr{\mathrm{ch}_{\mathrm{dR}}}

\newcommand{\Tr}{{\mathop{\mathrm{Tr}}}}

\newcommand{\tr}{{\mathop{\mathrm{tr}}}}
\newcommand{\Ac}{{\mathscr A}}
\newcommand{\Gc}{{\mathscr G}}
\newcommand{\te}{\theta}

\newcommand{\cqfd}{\hfill\rule{1ex}{1ex}}

\def\Id{\mathrm{Id}}

\def\d{\partial}
\def\dd{\mathrm{\bf d}}

\def\Bc{{\mathscr B}}
\def\Cc{{\mathscr C}}
\def\Jc{{\mathscr J}}
\def\Kc{{\mathscr K}}
\def\Fc{{\mathscr F}}

\def\im{\mathop{\mathrm{Im}}}
\def\ker{\mathop{\mathrm{Ker}}}

\def\bb{\overline{b}}

\def\hom{{\mathop{\mathrm{Hom}}}}

\def\hotimes{\hat{\otimes}}
\def\hotimesp{\hat{\otimes}_{\pi}}

\def\phit{\widetilde{\phi}}

\def\pt{\tilde{p}}

\def\gt{\widetilde{g}}
\def\ut{\widetilde{u}}
\def\omt{\widetilde{\om}}
\def\Gh{\widehat{G}}
\def\Uh{\widehat{U}}

\def\et{\widetilde{e}}

\def\ft{\widetilde{f}}
\def\Omh{\widehat{\Omega}}
\def\Th{\widehat{T}}
\def\chih{\widehat{\chi}}
\def\etah{\widehat{\eta}}

\def\Rch{\widehat{\mathscr R}}
\def\Mch{\widehat{\mathscr M}}

\def\Omb{\underline{\Omega}}
\def\zzb{\underline{\zz}}

\def\Xh{\widehat{X}}
\def\Jh{\widehat{J}}
\def\Tc{{\mathscr T}}
\def\Dc{{\mathscr D}}

\def\Kc{{\mathscr K}}

\def\mod{\ \mathrm{mod}\ }

\def\Hdr{H_{\mathrm{dR}}}
\def\Zdr{Z_{\mathrm{dR}}}
\def\Bdr{B_{\mathrm{dR}}}
\def\Hch{\check{H}}
\def\supp{\mathrm{supp}\,}

\def\eh{\hat{e}}
\def\gh{\hat{g}}
\def\uh{\hat{u}}
\def\fh{\hat{f}}
\def\vh{\hat{v}}
\def\omh{\hat{\om}}

\begin{document}

\begin{center}

{\bf   SECONDARY INVARIANTS FOR FRECHET ALGEBRAS\\[1mm] 
AND QUASIHOMOMORPHISMS}
\vskip 1cm
{\bf Denis PERROT}
\vskip 0.5cm
Universit\'e de Lyon, Universit\'e Lyon 1,\\
CNRS, UMR 5208 Institut Camille Jordan,\\
43, bd du 11 novembre 1918, 69622 Villeurbanne Cedex, France \\[2mm]
{\tt perrot@math.univ-lyon1.fr}\\[2mm]
\today
\end{center}
\vskip 0.5cm
\begin{abstract}
A Fr\'echet algebra endowed with a multiplicatively convex topology has two types of invariants: homotopy invariants (topological $K$-theory and periodic cyclic homology) and secondary invariants (multiplicative $K$-theory and the non-periodic versions of cyclic homology). The aim of this paper is to establish a Riemann-Roch-Grothendieck theorem relating direct images for homotopy and secondary invariants of Fr\'echet $m$-algebras under finitely summable quasihomomorphisms. 
\end{abstract}

\vskip 0.5cm


\section{Introduction}

For a noncommutative space described by an associative Fr\'echet algebra $\Ac$ over $\cc$, we distinguish two types of invariants. The first type are (smooth) homotopy invariants, for example topological $K$-theory \cite{Ph} and periodic cyclic homology \cite{C1}. The other type are secondary invariants; they are no longer stable under homotopy and carry a finer information about the ``geometry'' of the space $\Ac$. Typical examples of secondary invariants are algebraic $K$-theory \cite{R} (which will not be used here), multiplicative $K$-theory \cite{K2} and the unstable versions of cyclic homology \cite{L}. The aim of this paper is to define push-forward maps for homotopy and secondary invariants  between two Fr\'echet algebras $\Ac$ and $\Bc$, induced by a smooth finitely summable quasihomomorphism \cite{Cu}. The compatibility between the different types of invariants is expressed through a noncommutative Riemann-Roch-Grothendieck theorem (Theorem \ref{trr}). The present paper is the first part of a series on secondary characteristic classes. In the second part we will show how to obtain \emph{local formulas} for push-forward maps, following a general principle inspired by renormalization which establishes the link with chiral anomalies in quantum field theory \cite{P4}; in order to keep a reasonable size to the present paper, these methods will be published in a separate survey with further examples \cite{P5}.\\ 

We deal with Fr\'echet algebras endowed with a multiplicatively convex topology, or Fr\'echet $m$-algebras for short. These algebras can be presented as inverse limits of sequences of Banach algebras, and as a consequence many constructions valid for Banach algebras carry over Fr\'echet $m$-algebras. In particular Phillips \cite{Ph} defines topological $K$-theory groups $\Kt_n(\Ac)$ for any such algebra $\Ac$ and $n\in\zz$. The fundamental properties of interest for us are (smooth) homotopy invariance and Bott periodicity, i.e. $\Kt_{n+2}(\Ac) \cong \Kt_n(\Ac)$. Hence there are essentially two topological $K$-theory groups for any Fr\'echet $m$-algebra, $\Kt_0(\Ac)$ whose elements are roughly represented by idempotents in the stabilization of $\Ac$ by the algebra $\Kc$ of "smooth compact operators", and $\Kt_1(\Ac)$ whose elements are represented by invertibles. Fr\'echet $m$-algebras naturally arise in many situations related to differential geometry, commutative or not, and the formulation of index problems. In the latter situation one usually encounters an algebra $\Ic$ of "finitely summable operators", for us a Fr\'echet $m$-algebra provided with a continuous trace on its $p$-th power for some $p\geq 1$. A typical example is the Schatten class $\Ic=\Lc^p(H)$ of $p$-summable operators on an infinite-dimensional separable Hilbert space $H$. $\Ac$ can be stabilized by the completed projective tensor product $\Ic\hotimes\Ac$ and its topological $K$-theory $\Kt_n(\Ic\hotimes\Ac)$ is the natural receptacle for indices. Other important topological invariants of $\Ac$ (as a locally convex algebra) are provided by the periodic cyclic homology groups $HP_n(\Ac)$, which is the correct version sharing the properties of smooth homotopy invariance and periodicity mod 2 with topological $K$-theory \cite{C1}. For any finitely summable algebra $\Ic$ the Chern character $\Kt_n(\Ic\hotimes\Ac)\to HP_n(\Ac)$ allows to obtain cohomological formulations of index theorems. \\
If one wants to go beyond differential topology and detect \emph{secondary} invariants as well, which are no longer stable under homotopy, one has to deal with algebraic $K$-theory \cite{R} and the unstable versions of cyclic homology \cite{L}. In principle the algebraic $K$-theory groups $\Ka_n(\Ac)$ defined for any $n\in\zz$ provide interesting secondary invariants for any ring $\Ac$, but are very hard to calculate. It is also unclear if algebraic $K$-theory can be linked to index theory in a way consistent with topological $K$-theory, and in particular if it is possible to construct direct images of algebraic $K$-theory classes in a reasonable context. Instead, we will generalize slightly an idea of Karoubi \cite{K1, K2} and define for any Fr\'echet $m$-algebra $\Ac$ the multiplicative $K$-theory groups $MK^{\Ic}_n(\Ac)$, $n\in \zz$, indexed by a given finitely summable Fr\'echet $m$-algebra $\Ic$. Depending on the parity of the degree $n$, multiplicative $K$-theory classes are represented by idempotents or invertibles in certain extensions of $\Ic\hotimes\Ac$, together with a transgression of their Chern character in certain quotient complexes. Multiplicative $K$-theory is by definition a mixture of the topological $K$-theory $\Kt_n(\Ic\hotimes\Ac)$ and the non-periodic cyclic homology $HC_n(\Ac)$. It provides a "good" approximation of algebraic $K$-theory but is much more tractable. In addition, the Jones-Goodwillie Chern character in \emph{negative cyclic homology} $\Ka_n(\Ac)\to HN_n(\Ac)$ factors through multiplicative $K$-theory. The precise relations between topological, multiplicative $K$-theory and the various versions of cyclic homology are encoded in a commutative diagram whose rows are long exact sequences of abelian groups
\be
\vcenter{\xymatrix{
\Kt_{n+1}(\Ic\hotimes\Ac) \ar[r] \ar[d] & HC_{n-1}(\Ac) \ar[r]^{\delta} \ar@{=}[d] & MK^{\Ic}_n(\Ac)  \ar[r] \ar[d] & \Kt_n(\Ic\hotimes\Ac)  \ar[d]  \\
HP_{n+1}(\Ac) \ar[r]^S & HC_{n-1}(\Ac) \ar[r]^{\widetilde{B}} & HN_n(\Ac) \ar[r]^I & HP_n(\Ac) }} \label{invar}
\ee
The particular case $\Ic=\cc$ was already considered by Karoubi \cite{K1, K2} after the construction by Connes and Karoubi of regulator maps on algebraic $K$-theory \cite{CK}. The incorporation of a finitely summable algebra $\Ic$ is rather straightforward. This diagram describes the primary and secondary invariants associated to the noncommutative "manifold" $\Ac$. We mention that the restriction to Fr\'echet $m$-algebras is mainly for convenience. In principle these constructions could be extended to all locally convex algebras over $\cc$, however the subsequent results, in particular the proof of the Riemann-Roch-Grothendieck theorem would become much more involved. \\

If now $\Ac$ and $\Bc$ are two Fr\'echet $m$-algebras, it is natural to consider the adequate "morphisms" mapping the primary and secondary invariants from $\Ac$ to $\Bc$. Let $\Ic$ be a $p$-summable Fr\'echet $m$-algebra. By analogy with Cuntz' description of bivariant $K$-theory for $C^*$-algebras \cite{Cu}, if $\Ec\triangleright\Ic\hotimes\Bc$ denotes a Fr\'echet $m$-algebra containing $\Ic\hotimes\Bc$ as a (not necessarily closed) two-sided ideal, we define a \emph{$p$-summable quasihomomorphism} from $\Ac$ to $\Bc$ as a continuous homomorphism
$$
\rho: \Ac \to \Ec^s\triangleright \Ic^s\hotimes \Bc\ ,
$$
where $\Ec^s$ and $\Ic^s$ are certain $\zz_2$-graded algebras obtained from $\Ec$ and $\Ic$ by a standard procedure. Quasihomomorphisms come equipped with a parity (even or odd) depending on the construction of $\Ec^s$ and $\Ic^s$. In general, we may suppose that the parity is $p\mod 2$. We say that $\Ic$ is multiplicative if it is provided with a homomorphism $\boxtimes: \Ic\hotimes\Ic\to\Ic$, possibly defined up to adjoint action of multipliers on $\Ic$, and compatible with the trace. A basic example of multiplicative $p$-summable algebra is, once again, the Schatten class $\Lc^p(H)$. Then it is easy to show that such a quasihomomorphism induces a pushforward map in topological $K$-theory $\rho_!:\Kt_n(\Ic\hotimes\Ac)\to \Kt_{n-p}(\Ic\hotimes\Bc)$, whose degree coincides with the parity of the quasihomomorphism. This is what one expects from bivariant $K$-theory and is not really new. Our goal is to extend this map to the entire diagram (\ref{invar}). Direct images for the unstable versions of cyclic homology are necessarily induced by a bivariant \emph{non-periodic} cyclic cohomology class $\ch^p(\rho)\in HC^p(\Ac,\Bc)$. This bivariant Chern character exists only under certain admissibility properties about the algebra $\Ec$ (note that it is sufficient for $\Ic$ to be $(p+1)$-summable instead of $p$-summable). In particular, the bivariant Chern caracter constructed by Cuntz for any quasihomomorphism in \cite{Cu1, Cu2} cannot be used here because it provides a bivariant periodic cyclic cohomology class, which does not detect the secondary invariants of $\Ac$ and $\Bc$. We give the precise definition of an admissible quasihomomorphism and construct the bivariant Chern character $\ch^p(\rho)$ in section \ref{sbiv}, on the basis of previous works \cite{P2}. An analogous construction was obtained by Nistor \cite{Ni1, Ni2} or by Cuntz and Quillen \cite{CQ1}. However the bivariant Chern character of \cite{P2} is related to other constructions involving the heat operator and can be used concretely for establishing local index theorems, see for example \cite{P3}. The pushforward map in topological $K$-theory combined with the bivariant Chern character leads to a pushforward map in multiplicative $K$-theory $\rho_!:MK^{\Ic}_n(\Ac)\to MK^{\Ic}_{n-p}(\Bc)$. Our first main result is the following non-commutative version of the Riemann-Roch-Grothendieck 
theorem (see Theorem \ref{trr} for a precise statement):

\begin{theorem}
Let $\rho:\Ac\to \Ec^s\triangleright \Ic^s\hotimes\Bc$ be an admissible quasihomomorphism of parity $p$ mod 2. Suppose that $\Ic$ is $(p+1)$-summable in the even case and $p$-summable in the odd case. Then one has a graded-commutative diagram
$$
\vcenter{\xymatrix{
\Kt_{n+1}(\Ic\hotimes\Ac) \ar[r] \ar[d]^{\rho_!} & HC_{n-1}(\Ac) \ar[r] \ar[d]^{\ch^p(\rho)} & MK^{\Ic}_n(\Ac)  \ar[r] \ar[d]^{\rho_!}  & \Kt_n(\Ic\hotimes\Ac)  \ar[d]^{\rho_!}  \\
\Kt_{n+1-p}(\Ic\hotimes\Bc) \ar[r]  & HC_{n-1-p}(\Bc) \ar[r]  & MK^{\Ic}_{n-p}(\Bc)  \ar[r]  & \Kt_{n-p}(\Ic\hotimes\Bc) }} 
$$
compatible with the cyclic homology $SBI$ exact sequences after taking the Chern characters $\Kt_*(\Ic\hotimes\cdot )\to HP_*$ and $MK^{\Ic}_*\to HN_*$.
\end{theorem}

\noindent At this point it is interesting to note that the pushforward maps $\rho_!$ and the bivariant Chern character $\ch^p(\rho)$ enjoy some invariance properties with respect to equivalence relations among quasihomomorphisms. Two types of equivalence relations are defined: smooth homotopy and conjugation by invertibles. The second relation corresponds to "compact perturbation" in Kasparov $KK$-theory for $C^*$-algebras \cite{Bl}. In the latter situation, the $M_2$-stable version of conjugation essentially coincide with homotopy, at least for separable $\Ac$ and $\si$-unital $\Bc$. For Fr\'echet algebras however, $M_2$-stable conjugation is \emph{strictly stronger} than homotopy as an equivalence relation. This is indeed in the context of Fr\'echet algebras that secondary invariants appear. The pushforward maps in topological $K$-theory and periodic cyclic homology are invariant under homotopy of quasihomomorphisms. The maps in multiplicative $K$-theory and the non-periodic versions of cyclic homology $HC_*$ and $HN_*$ are only invariant under conjugation and not homotopy. Also note that in contrast with the $C^*$-algebra situation, the Kasparov product of two quasihomomorphisms $\rho:\Ac\to \Ec^s\triangleright \Ic^s\hotimes \Bc$ and $\rho':\Bc\to \Fc^s\triangleright \Ic^s\hotimes \Cc$ is not defined as a quasihomomorphism from $\Ac$ to $\Cc$. The various bivariant $K$-theories constructed for $m$-algebras \cite{Cu1, Cu2} or even for general bornological algebras \cite{CMR} cannot be used here, again because they are homotopy invariant by construction. We leave the construction of a bivariant $K$-theory compatible with secondary invariants as an open problem.\\

In the last part of the paper we illustrate the Riemann-Roch-Grothendieck theorem by constructing assembly maps for certain crossed product algebras. If $\Gamma$ is a discrete group acting on a Fr\'echet $m$-algebra $\Ac$, under certain conditions the crossed product $\Ac\cp\Gamma$ is again a Fr\'echet $m$-algebra and one would like to obtain multiplicative $K$-theory classes out of a geometric model inspired by the Baum-Connes construction \cite{BC}. Thus let $P\stackrel{\Gamma}{\to} M$ be a principal $\Gamma$-bundle over a compact manifold $M$, and denote by $\Ac_P$ the algebra of smooth sections of the associated $\Ac$-bundle. If $D$ is a $K$-cycle for $M$ represented by a pseudodifferential operator, we obtain a quasihomomorphism from $\Ac_P$ to $\Ac\cp\Gamma$ and hence a map
$$
MK^{\Ic}_n(\Ac_P)\to MK^{\Ic}_{n-p}(\Ac\cp\Gamma)
$$
for suitable $p$ and Schatten ideal $\Ic$. In general this map cannot exhaust the entire multiplicative $K$-theory of the crossed product but nevertheless interesting secondary invariants arise in this way. In the case where $\Ac$ is the algebra of smooth functions on a compact manifold, $\Ac_P$ is commutative and its secondary invaiants are closely related to (smooth) Deligne cohomology. From this point of view the pushforward map in multiplicative $K$-theory should be considered as a non-commutative version of ``integrating Deligne classes along the fibers'' of a submersion. We perform the computations for the simple example provided by the noncommutative torus. \\

The paper is organized as follows. In section \ref{scy} we review the Cuntz-Quillen formulation of (bivariant) cyclic cohomology \cite{CQ1} in terms of quasi-free extensions for $m$-algebras. Nothing is new but we take the opportunity to fix the notations and recall a proof of generalized Goodwillie theorem. In section \ref{sbiv} we define quasihomomorphisms and construct the bivariant Chern character. The formulas are identical to those found in \cite{P2} but in addition we carefully establish their adic properties and conjugation invariance. In section \ref{st} we recall Phillips' topological $K$-theory for Fr\'echet $m$-algebras, and introduce the periodic Chern character $\Kt_n(\Ic\hotimes\Ac)\to HP_n(\Ac)$ when $\Ic$ is a finitely summable algebra. The essential point here is to give explicit and simple formulas for subsequent use. Section \ref{sm} is devoted to the definition of the multiplicative $K$-theory groups $MK^{\Ic}_n(\Ac)$ and the proof of the long exact sequence relating them with topological $K$-theory and cyclic homology. We also construct the negative Chern character $MK^{\Ic}_n(\Ac)\to HN_n(\Ac)$ and show the compatibility with the $SBI$ exact sequence. Direct images of topological and multiplicative $K$-theory under quasihomomorphisms are constructed in section \ref{srr} and the Riemann-Roch-Grothendieck theorem is proved. The example of assembly maps and crossed products is treated in section \ref{sassem}.

\section{Cyclic homology}\label{scy}

Cyclic homology can be defined for various classes of associative algebras over $\cc$, in particular complete locally convex algebras. For us, a locally convex algebra $\Ac$ has a topology induced by a family of continuous seminorms $p:\Ac\to \rr_+$, for which the multiplication $\Ac\times\Ac\to\Ac$ is jointly continuous. Hence for any seminorm $p$ there exists a seminorm $q$ such that $p(a_1a_2)\leq q(a_1)q(a_2)$, $\forall a_i\in \Ac$. For technical reasons however we shall restrict ourselves to \emph{multiplicatively convex} algebras \cite{C1}, whose topology is generated by a family of submultiplicative seminorms 
$$
p(a_1a_2)\leq p(a_1)p(a_2)\quad \forall a_i \in \Ac\ .
$$
A complete multiplicatively convex algebra is called $m$-algebra, and may equivalently be described as a projective limit of Banach algebras. The unitalization $\Ac^+=\cc\oplus\Ac$ of an $m$-algebra $\Ac$ is again an $m$-algebra, for the seminorms $\tilde{p}(\la 1 + a)=|\la| + p(a)$, $\forall \la\in\cc, a\in\Ac$. In the same way, if $\Bc$ is another $m$-algebra, the direct sum $\Ac\oplus\Bc$ is an $m$-algebra for the topology generated by the seminorms $(p\oplus q)(a,b)=p(a)+q(b)$, where $p$ is a seminorm on $\Ac$ and $q$ a seminorm on $\Bc$. Also, the algebraic tensor product $\Ac\otimes\Bc$ may be endowed with the projective topology induced by the seminorms 
\be
(p\otimes q)(c)=\inf \Big\{\sum_{i=1}^np(a_i)q(b_i)\ \mbox{such that}\ c=\sum_{i=1}^na_i\otimes b_i\ \in\Ac\otimes\Bc \Big\}\ .\label{proj}
\ee 
The \emph{completion} $\Ac\hotimes\Bc=\Ac\hotimesp\Bc$ of the algebraic tensor product under this family of seminorms is the projective tensor product of Grothendieck \cite{Gr}, and is again an $m$-algebra. \\
Cyclic homology, cohomology and bivariant cyclic cohomology for $m$-algebras can be defined either within the cyclic bicomplex formalism of Connes \cite{C1}, or the $X$-complex of Cuntz and Quillen \cite{CQ1}. We will make an extensive use of both formalisms throughout this paper. In general, we suppose that all linear maps or homomorphims between $m$-algebras are continuous, tensor products are completed projective tensor products, and extensions of $m$-algebras $0\to \Ic \to \Rc \to \Ac \to 0$ always admit a continuous linear splitting $\si: \Ac \to \Rc$. 

\subsection{Cyclic bicomplex}

\noindent {\bf Non-commutative differential forms.} Let $\Ac$ be an $m$-algebra. The space of non-commutative differential forms over $\Ac$ is the algebraic direct sum $\Om\Ac=\bigoplus_{n\ge 0}\Om^n\Ac$ of the $n$-forms subspaces $\Om^n\Ac=\Ac^+\hotimes\Ac^{\hotimes n}$ for $n\ge 1$ and $\Om^0\Ac=\Ac$, where $\Ac^+$ is the unitalization of $\Ac$. Each of the subspaces $\Om^n\Ac$ is complete but we do not complete the direct sum. It is customary to use the differential notation $a_0da_1\ldots da_n$ (resp. $da_1\ldots da_n$) for the string $a_0\otimes a_1\ldots\otimes a_n$ (resp. $1\otimes a_1\ldots\otimes a_n)$. A continuous differential $d: \Om^n\Ac\to\Om^{n+1}\Ac$ is uniquely specified by $d(a_0da_1\ldots da_n)=da_0da_1\ldots da_n$ and $d^2=0$. A continuous and associative product $\Om^n\Ac\times \Om^m\Ac\to \Om^{n+m}\Ac$ is defined as usual and fulfills the Leibniz rule $d(\om_1\om_2)=d\om_1\om_2 +(-)^{|\om_1|}\om_1d\om_2$, where $|\om_1|$ is the degree of $\om_1$. This turns $\Om\Ac$ into a differential graded (DG) algebra. \\
On $\Om\Ac$ are defined various operators. First of all, the Hochschild boundary map $b: \Om^{n+1}\Ac\to \Om^n\Ac$ reads $b(\om da)=(-)^n[\om,a]$ for $\om\in\Om^n\Ac$, and $b=0$ on $\Om^0\Ac=\Ac$. One easily shows that $b$ is continuous and $b^2=0$, hence $\Om\Ac$ is a complex graded over $\nn$. The Hochschild homology of $\Ac$ (with coefficients in the bimodule $\Ac$) is the homology of this complex:
\be
HH_n(\Ac)=H_n(\Om\Ac,b)\ ,\qquad \forall n\in\nn\ .
\ee
Then the Karoubi operator $\kappa:\Om^n\Ac\to \Om^n\Ac$ is defined by $1-\kappa=db+bd$. Therefore $\kappa$ is continuous and commutes with $b$ and $d$. One has $\kappa(\om\, da)=(-)^n da\,\om$ for any $\om\in \Om^n\Ac$ and $a\in\Ac$. The last operator is Connes' $B:\Om^n\Ac\to \Om^{n+1}\Ac$, equal to $(1+\kappa+\ldots+\kappa^n)d$ on $\Om^n\Ac$. It is also continuous and verifies $B^2=0=Bb+bB$ and $B\kappa=\kappa B=B$. Thus $\Om\Ac$ endowed with the two anticommuting differentials $(b,B)$ becomes a bicomplex. It splits as a direct sum $\Om\Ac=\Om\Ac^+\oplus\Om\Ac^-$ of even and odd degree differential forms, hence is a $\zz_2$-graded complex for the total boundary map $b+B$. However its homology is trivial \cite{L}. The various versions of cyclic homology are defined using the natural filtrations on $\Om\Ac$. Following Cuntz and Quillen \cite{CQ1}, we define the \emph{Hodge filtration} on $\Om\Ac$ as the decreasing family of $\zz_2$-graded subcomplexes for the total boundary $b+B$
$$
F^n\Om\Ac= b\Om^{n+1}\Ac\oplus \bigoplus_{k>n}\Om^k\Ac\ ,\qquad \forall n\in\zz\ ,
$$
with the convention that $F^n\Om\Ac=\Om\Ac$ for $n<0$. The completion of $\Om\Ac$ is defined as the projective limit of $\zz_2$-graded complexes 
\be
\Omh\Ac = \varprojlim_n \Om\Ac/F^n\Om\Ac = \prod_{n\geq 0}\Om^n\Ac\ .
\ee
Hence $\Omh\Ac=\Omh^+\Ac\oplus\Omh^-\Ac$ is a $\zz_2$-graded complex endowed with the total boundary map $b+B$. It is itself filtered by the decreasing family of $\zz_2$-graded subcomplexes $F^n\Omh\Ac= \ker (\Omh\Ac\to \Om\Ac/F^n\Om\Ac)$, which may be written 
\be
F^n\Omh\Ac= b\Om^{n+1}\Ac\oplus \prod_{k>n}\Om^k\Ac\ ,\qquad \forall n\in\zz\ .
\ee
In particular the quotient $ \Omh\Ac/F^n\Omh\Ac$ is a $\zz_2$-graded complex isomorphic to $\Om\Ac/F^n\Om\Ac$, explicitly
\be
\Omh\Ac/F^n\Omh\Ac=\bigoplus_{k=0}^{n-1}\Om^k\Ac\oplus\Om^n\Ac/b(\Om^{n+1}\Ac)\ ,
\ee
and it vanishes for $n<0$. As a topological vector space, $\Omh\Ac/F^n\Omh\Ac$ may fail to be separated because the image $b(\Om^{n+1}\Ac)$ is not closed in general.  
\begin{definition}
In any degree $n\in\zz$, the \emph{periodic, non-periodic} and \emph{negative} cyclic homologies are respectively the $(b+B)$-homologies 
\beq
HP_n(\Ac) &=& H_{n+2\zz}(\Omh\Ac)\ ,\non\\
HC_n(\Ac) &=& H_{n+2\zz}(\Omh\Ac/F^n\Omh\Ac)\ ,\\
HN_n(\Ac) &=& H_{n+2\zz}(F^{n-1}\Omh\Ac)\ .\non
\eeq
\end{definition} 
Hence $HP_n(\Ac)\cong HP_{n+2}(\Ac)$ is 2-periodic, $HC_n(\Ac)=0$ for $n<0$ and $HN_n(\Ac)\cong HP_n(\Ac)$ for $n\leq 0$. By construction these cyclic homology groups fit into a long exact sequence
\be
\ldots \longrightarrow  HP_{n+1}(\Ac) \stackrel{S}{\longrightarrow}   HC_{n-1}(\Ac) \stackrel{B}{\longrightarrow}   HN_n(\Ac) \stackrel{I}{\longrightarrow}  HP_n(\Ac) \longrightarrow  \ldots  \label{SBI}
\ee
where $S$ is induced by projection, $I$ by inclusion, and the connecting map corresponds to the operator $B$. The link between cyclic and Hochshild homology may be obtained through \emph{non-commutative de Rham homology} \cite{K1}, defined as
\be
HD_n(\Ac) := H_{n+2\zz}(\Omh\Ac/F^{n+1}\Omh\Ac)\ ,\quad \forall n\in\zz\ .
\ee
This yields for any $n\in \zz$ a short exact sequence of $\zz_2$-graded complexes 
$$
0 \longrightarrow G_n(\Ac) \longrightarrow \Omh\Ac/F^n\Omh\Ac \longrightarrow \Omh\Ac/F^{n-1}\Omh\Ac \longrightarrow 0 \ ,
$$
where $G_n$ is $\Om^n\Ac/b\Om^{n+1}\Ac$ in degree $n$ mod 2, and $b\Om^n\Ac$ in degree $n-1$ mod 2. One has $H_{n+2\zz}(G_n)=HH_n(\Ac)$ and $H_{n-1+2\zz}(G_n)=0$, so that the associated six-term cyclic exact sequence in homology reduces to
$$
0 \to HD_{n-1}(\Ac) \to HC_{n-1}(\Ac) \to HH_n(\Ac) \to HC_n(\Ac) \to HD_{n-2}(\Ac) \to 0 \ ,
$$
and Connes's $SBI$ exact sequence \cite{C0} for cyclic homology is actually obtained by splicing together the above sequences for all $n\in\zz$:
\be
\ldots \longrightarrow HC_{n+1}(\Ac) \stackrel{S}{\longrightarrow}  HC_{n-1}(\Ac) \stackrel{B}{\longrightarrow}   HH_n(\Ac) \stackrel{I}{\longrightarrow}  HC_n(\Ac) \longrightarrow \ldots   \label{exa}
\ee
Hence the non-commutative de Rham homology group $HD_n(\Ac)$ may be identified with the image of the periodicity shift $S:HC_{n+2}(\Ac)\to HC_n(\Ac)$. Clearly the exact sequence (\ref{SBI}) can be transformed to (\ref{exa}) by taking the natural maps $HP_n(\Ac)\to HC_n(\Ac)$ and $HN_n(\Ac)\to HH_n(\Ac)$.\\

Passing to the dual theory, let $\hom(\Omh\Ac,\cc)$ be the $\zz_2$-graded complex of linear maps $\Omh\Ac\to\cc$ which are continuous for the adic topology on $\Omh\Ac$ induced by the Hodge filtration. It is concretely described as the direct sum
$$
\hom(\Omh\Ac,\cc)=\bigoplus_{n\geq 0}\hom(\Om^n\Ac,\cc)\ ,
$$
where $\hom(\Om^n\Ac,\cc)$ is the space of continuous linear maps $\Om^n\Ac\to\cc$. The space $\hom(\Omh\Ac,\cc)$ is endowed with the transposed of the boundary operator $b+B$ on $\Omh\Ac$. Then the periodic cyclic cohomology of $\Ac$ is the cohomology of this complex:
\be
HP^n(\Ac)= H^{n+2\zz}(\hom(\Omh\Ac,\cc))\ .
\ee
One defines analogously the non-periodic and negative cyclic cohomologies which fit into an $IBS$ long exact sequence.

\subsection{$X$-complex and quasi-free algebras}

We now turn to the description of the $X$-complex. It first appeared in the coalgebra context in Quillen's work \cite{Q2}, and subsequently was used by Cuntz and Quillen in their formulation of cyclic homology \cite{CQ1}. Here we recall the $X$-complex construction for $m$-algebras. \\
Let $\Rc$ be an $m$-algebra. The space of non-commutative one-forms $\Om^1\Rc$ is a $\Rc$-bimodule, hence we can take its quotient $\Om^1\Rc_{\nat}$ by the subspace of commutators $[\Rc,\Om^1\Rc]=b\Om^2\Rc$. $\Om^1\Rc_{\nat}$ may fail to be separated in general. However, it is automatically separated when $\Rc$ is \emph{quasi-free}, see below. In order to avoid confusions in the subsequent notations, we always write a one-form $x_0\dd x_1\in\Om^1\Rc$ with a bold $\dd$ when dealing with the $X$-complex of $\Rc$. The latter is the $\zz_2$-graded complex \cite{CQ1}
\be
X(\Rc):\quad \Rc\ \mathop{\rightleftarrows}^{\nat \dd}_{\bb}\ \Om^1\Rc_{\nat}\ , 
\ee
where $\Rc=X_+(\Rc)$ is located in even degree and $\Om^1\Rc_{\nat}=X_-(\Rc)$ in odd degree. The class of the generic element $(x_0\dd x_1 \mod [,])\in \Om^1\Rc_{\nat}$ is usually denoted by $\nat x_0\dd x_1$. The map $\nat \dd:\Rc\to \Om^1\Rc_{\nat}$ thus sends $x\in\Rc$ to $\nat \dd x$. Also, the Hochschild boundary $b:\Om^1\Rc\to \Rc$ vanishes on the commutator subspace $[\Rc,\Om^1\Rc]$, hence passes to a well-defined map $\bb:\Om^1\Rc_{\nat}\to\Rc$. Explicitly the image of $\nat x_0\dd x_1$ by $\bb$ is the commutator $[x_0,x_1]$. These maps are continuous and satisfy $\nat \dd\circ \bb=0$ and $\bb\circ\nat \dd=0$, so that $(X(\Rc),\nat\dd\oplus\bb)$ indeed defines a $\zz_2$-graded complex. We mention that everything can be formulated when $\Rc$ itself is a $\zz_2$-graded algebra: we just have to replace everywhere the ordinary commutators by graded commutators, and the differentials anticommute with elements of odd degree. In particular one gets $\bb\nat x\dd y=(-)^{|x|}[x,y]$, where $|x|$ is the degree of $x$ and $[x,y]$ is the graded commutator. The $X$-complex is obviously a functor from $m$-algebras to $\zz_2$-graded complexes: if $\rho:\Rc\to\Sc$ is a continuous homomorphism, it induces a chain map of even degree $X(\rho):X(\Rc)\to X(\Sc)$, by setting $X(\rho)(x)=\rho(x)$ and $X(\rho)(\nat x_0\dd x_1)=\nat \, \rho(x_0)\dd\rho(x_1)$.\\
In fact the $X$-complex may be identified with the quotient of the $(b+B)$-complex $\Omh\Rc$ by the subcomplex $F^1\Omh\Rc= b\Om^{2}\Rc\oplus \prod_{k\geq 2}\Om^k\Rc$ of the Hodge filtration, i.e. there is an exact sequence
$$
0\to F^1\Omh\Rc \to \Omh\Rc \to X(\Rc) \to 0 \ .
$$
It turns out that the $X$-complex is especially designed to compute the cyclic homology of algebras for which the subcomplex $F^1\Omh\Rc$ is contractible. This led Cuntz and Quillen to the following definition:
\begin{definition}[\cite{CQ1}]
An $m$-algebra $\Rc$ is called \emph{quasi-free} if there exists a continuous linear map $\phi: \Rc\to \Om^2\Rc$ with property
\be
\phi(xy)= \phi(x) y + x\phi(y) +\dd x \dd y\ ,\quad \forall x,y\in\Rc\ .
\ee
\end{definition}
We refer to \cite{CQ1,Me} for many other equivalent definitions of quasi-free algebras. Let us just observe that a quasi-free algebra has dimension $\leq 1$ with respect to Hochschild cohomology. Indeed, the map $\phi$ allows to contract the Hochschild complex of $\Rc$ in dimensions $> 1$, and this contraction carries over to the cyclic bicomplex. First, the linear map 
$$
\si:\Om^1\Rc_{\nat}\to \Om^1\Rc\ ,\qquad \nat x\dd y \mapsto x\dd y + b(x\phi(y))
$$
is well-defined because it vanishes on the commutator subspace $[\Rc,\Om^1\Rc]=b\Om^2\Rc$ by the algebraic property of $\phi$. Hence $\si$ is a continuous linear splitting of the exact sequence $0\to b\Om^2\Rc \to \Om^1\Rc \to \Om^1\Rc_{\nat} \to 0$. By the way, this implies that for a quasi-free algebra $\Rc$, the topological vector space $\Om^1\Rc$ splits into the direct sum of two closed subspaces $b\Om^2\Rc$ and $\Om^1\Rc_{\nat}$. Then, we extend $\phi$ to a continuous linear map $\phi : \Om^n\Rc\to \Om^{n+2}\Rc$ in all degrees $n\geq 1$ by the formula
$$
\phi(x_0\dd x_1\ldots \dd x_n)=\sum_{i=0}^n (-)^{ni}\phi(x_i)\dd x_{i+1}\ldots \dd x_n\dd x_0\ldots \dd x_{i-1}\ .
$$
The following proposition gives a chain map $\gamma: X(\Rc)\to \Omh\Rc$ which is inverse to the natural projection $\pi: \Omh\Rc\to X(\Rc)$ up to homotopy. Remark that the infinite sum $(1-\phi)^{-1}:=\sum_{n=0}^{\infty}\phi^n$ makes sense as a linear map $\Rc\to \Omh^+\Rc$ or $\Om^1\Rc\to \Omh^-\Rc$.

\begin{proposition}\label{pgamma}
Let $\Rc$ be a quasi-free $m$-algebra. Then\\
i) The map $\gamma: X(\Rc)\to \Omh\Rc$ defined for $x,y\in\Rc$ by
\beq
\gamma(x) &=& (1-\phi)^{-1}(x)\\
\gamma(\nat x\dd y) &=& (1-\phi)^{-1}(x\dd y + b(x\phi(y)))\non
\eeq
is a chain map of even degree from the $X$-complex to the $(b+B)$-complex.\\
ii) Let $\pi:\Omh\Rc\to X(\Rc)$ be the natural projection. There is a contracting homotopy of odd degree $h:\Omh\Rc\to\Omh\Rc$ such that
\beq
\pi\circ\gamma &=& \Id\quad \mbox{on}\quad X(\Rc)\ , \non\\
\gamma\circ\pi &=& \Id +[b+B,h] \quad \mbox{on}\quad \Omh\Rc\ .\non
\eeq
Hence $X(\Rc)$ and $\Omh\Rc$ are homotopy equivalent.
\end{proposition}
{\it Proof:} See the proof of \cite{P1}, Proposition 4.2. There the result was stated in the particular case of a tensor algebra $\Rc=T\Ac$, but the general case of a quasi-free algebra is strictly identical (with the tensor algebra the image of $\gamma$ actually lands in the subcomplex $\Om T\Ac\subset \Omh T\Ac$ for a judicious choice of $\phi$, but for generic quasi-free algebras it is necessary to take the completion $\Omh\Rc$ of $\Om\Rc$). \cqfd\\

\noindent {\bf Extensions.} Let $\Ac$ be an $m$-algebra. By an extension of $\Ac$ we mean an exact sequence of $m$-algebras $0\to \Ic \to \Rc \to \Ac \to 0$ provided with a continuous linear splitting $\Ac\to \Rc$, and the topology of the ideal $\Ic$ is induced by its inclusion in $\Rc$. Hence as a topological vector space $\Rc$ is the direct sum of the closed subspaces $\Ic$ and $\Ac$. By convention, the powers $\Ic^n$ of the ideal $\Ic$ will always denote the image in $\Rc$ of the $n$-th tensor power $\Ic\hotimes\ldots\hotimes\Ic$ by the multiplication map. For $n\leq 0$, we define $\Ic^0$ as the algebra $\Rc$. Now let us suppose that all the powers $\Ic^n$ are closed and direct summands in $\Rc$ (this is automatically satisfied if $\Rc$ is quasi-free). Then the quotients $\Rc/\Ic^n$ are $m$-algebras and give rise to an inverse system with surjective homomorphisms
$$
0\leftarrow \Ac= \Rc/\Ic \leftarrow \Rc/\Ic^2 \leftarrow \ldots \leftarrow \Rc/\Ic^n \leftarrow \ldots 
$$
We denote by $\Rch=\varprojlim_n \Rc/\Ic^n$ the projective limit and view it as a \emph{pro-algebra} indexed by the directed set $\zz$ (see \cite{Me}). Since the bicomplex of non-commutative differential forms $\Omh\Rc$ and the $X$-complex $X(\Rc)$ are functorial in $\Rc$, we can define $\Omh\Rch$ and $X(\Rch)$ as the $\zz_2$-graded pro-complexes
\beq
\Omh\Rch &=& \varprojlim_n \Omh(\Rc/\Ic^n)= \varprojlim_{m,n} \Om(\Rc/\Ic^n)/F^m\Om(\Rc/\Ic^n)\ , \non\\
X(\Rch) &=& \varprojlim_n X(\Rc/\Ic^n)\ ,  \non
\eeq
endowed respectively with the total boundary maps $b+B$ and $\partial=\nat\dd\oplus\bb$. When $\Rc$ is quasi-free, a refinement of Proposition \ref{pgamma} yields a chain map $\gamma: X(\Rch)\to \Omh\Rch$ inverse to the projection $\pi:\Omh \Rch\to X(\Rch)$ up to homotopy, which we call a \emph{ generalized Goodwillie equivalence}:

\begin{proposition}\label{pgood}
Let $0\to \Ic \to \Rc \to \Ac \to 0$ be an extension of $m$-algebras, with $\Rc$ quasi-free. Then the chain map $\gamma: X(\Rc)\to \Omh\Rc$ extends to a homotopy equivalence of pro-complexes $X(\Rch)\to \Omh\Rch$.
\end{proposition}
{\it Proof:} We recall the proof because it will be useful for establishing Proposition \ref{padic}. Let us introduce the following decreasing filtration of the space $\Om^m\Rc$ by the subspaces $H^k\Om^m\Rc$, $k\in \zz$:
$$
H^k\Om^m\Rc=\sum_{k_0+\ldots+ k_m \geq k} \Ic^{k_0}\dd\Ic^{k_1}\ldots \dd\Ic^{k_m}\ .
$$
Clearly $H^{k+1}\Om^m\Rc\subset H^k\Om^m\Rc$, and for $k\leq 0$ $H^k\Om^m\Rc=\Om^m\Rc$. Morally, $H^k\Om^m\Rc$ contains at least $k$ powers of the ideal $\Ic$. The direct sum $\bigoplus_{m} H^k\Om^m\Rc$ is stable by the operators $d$, $b$, $\kappa$, $B$ for any $k$. We have to establish how $k$ changes when the linear map $\phi: \Om^m\Rc\to \Om^{m+2}\Rc$ is applied. First consider $\phi:\Rc\to \Om^2\Rc$. If $x_1,\ldots, x_k$ denote $k$ elements in $\Rc$, one has by the algebraic property of $\phi$ (see \cite{CQ1})
\beq
\lefteqn{\phi(x_1\ldots x_k)=\sum_{i=1}^k x_1\ldots x_{i-1}\phi(x_i)x_{i+1}\ldots x_k } \non\\ 
&& \qquad \qquad + \sum_{1\leq i < j\leq k} x_1\ldots x_{i-1}\dd x_i x_{i+1}\ldots x_{j-1} \dd x_j x_{j+1}\ldots x_k\ .\non
\eeq
Taking the elements $x_i$ in the ideal $\Ic$ and using that $\phi(\Ic)\subset \Om^2\Rc$ yields
$$
\phi(\Ic^k)\subset \sum_{i=1}^k \Ic^{i-1} \dd\Rc\dd\Rc \Ic^{k-i} + \sum_{1\leq i < j\leq k} \Ic^{i-1}\dd\Ic\Ic^{j-i-1}\dd\Ic\Ic^{k-j}\ .
$$
Therefore $\phi(\Ic^k)\subset H^{k-1}\Om^2\Rc$ for any $k$. Now from the definition of $\phi$ on $\Om^m\Rc$, one has 
$$
\phi(\Ic^{k_0}\dd\Ic^{k_1}\ldots \dd\Ic^{k_m})\subset \sum_{l=0}^m\phi(\Ic^{k_l})\dd\Ic^{k_{l+1}}\ldots \dd\Ic^{k_{l-1}} \subset H^{k-1}\Om^{m+2}\Rc
$$
whenever $k=k_0+\ldots+k_m$, hence $\phi(H^k\Om^m\Rc)\subset H^{k-1}\Om^{m+2}\Rc$. Now let us evaluate the even part of the chain map $\gamma:\Rc\to \Omh^+\Rc$. The part of $\gamma$ landing in $\Om^{2m}\Rc$ is the $m$-th power $\phi^m$. One gets $\phi^m(\Ic^k)\subset H^{k-m}\Om^{2m}\Rc$, hence $\phi^m$ sends the quotient algebra $\Rc/\Ic^k$ to $\Om^{2m}(\Rc/\Ic^n)$ provided $(1+2m) n \leq k-m$ (indeed $1+2m$ is the maximal number of factors $\Rc$ in the tensor product $\Om^{2m}\Rc$). Passing to the projective limits, $\phi^m$ induces a well-defined map $\Rch\to \Om^{2m}\Rch$, and summing over all degrees $2m$ yields $\gamma:\Rch\to \Omh^+\Rch$. \\
Let us turn to the odd part of the chain map $\gamma:\Om^1\Rc_{\nat}\to \Omh^-\Rc$. By construction, it is the composition of the linear map $\si:\Om^1\Rc_{\nat}\to \Om^1\Rc$ with all the powers $\phi^m:\Om^1\Rc\to \Om^{2m+1}\Rc$. One has $\nat(\Ic^k\dd\Rc+ \Rc\dd(\Ic^k)) \subset \nat(\Ic^k\dd\Rc+\Ic^{k-1}\dd\Ic)$, and by the definition of $\si$,
\beq
\si\, \nat(\Ic^k\dd\Rc+\Ic^{k-1}\dd\Ic) &\subset& \Ic^k\dd\Rc+\Ic^{k-1}\dd\Ic + b(\Ic^k\phi(\Rc)+\Ic^{k-1}\phi(\Ic))\non\\
&\subset& H^k\Om^1\Rc + bH^{k-1}\Om^2\Rc\ \subset\ H^{k-1}\Om^1\Rc\ .\non
\eeq
Therefore $(\phi^m\circ\si) \nat(\Ic^k\dd\Rc+\Ic^{k-1}\dd\Ic)\subset H^{k-m-1}\Om^{2m+1}\Rc$. Since $\Rc$ is the direct sum of $\Rc/\Ic^k$ and $\Ic^k$ as a topological vector space, the quotient $\Om^1(\Rc/\Ic^k)_{\nat}$ coincides with $\Om^1\Rc/(\Ic^k\dd\Rc+ \Rc\dd(\Ic^k) + [\Rc,\Om^1\Rc])$, and the map $\phi^m\circ\si: \Om^1(\Rc/\Ic^k)_{\nat} \to \Om^{2m+1}(\Rc/\Ic^n)$ is well-defined provided $(2m+2)n \leq k-m-1$. Thus passing to the projective limits induces a map $\Om^1\Rch_{\nat}\to \Om^{2m+1}\Rch$, and summing over $m$ yields $\gamma: \Om^1\Rch_{\nat}\to \Omh^-\Rch$.\\ 
Finally, the contracting homotopy $h$ of Proposition \ref{pgamma} is also constructed from $\phi$ (see \cite{P1} Proposition 4.2), hence extends to a contracting homotopy $h:\Omh\Rch\to \Omh\Rch$. The relations $\pi\circ\gamma = \Id$ on $X(\Rch)$ and $\gamma\circ\pi = \Id +[b+B,h]$ on $\Omh\Rch$ follow immediately. \cqfd \\

\noindent {\bf Adic filtration.} Suppose that $\Ic$ is a (not necessarily closed) two-sided ideal in $\Rc$, provided with its own topology of $m$-algebra for which the inclusion $\Ic\to\Rc$ is continuous and the multiplication map $\Rc^+\times \Ic\times\Rc^+\to\Ic$ is jointly continuous. As usual we define the powers $\Ic^n$ as the two-sided ideals corresponding to the image in $\Rc$ of the $n$-fold tensor products $\Ic\hotimes\ldots\hotimes\Ic$ under multiplication. Following \cite{CQ1}, we introduce the adic filtration of $X(\Rc)$ by the subcomplexes
\beq 
F_{\Ic}^{2n}X(\Rc) &:& \Ic^{n+1}+[\Ic^n,\Rc] \ \rightleftarrows \ \nat \Ic^n\dd\, \Rc \label{filtration}\\
F_{\Ic}^{2n+1}X(\Rc) &:& \Ic^{n+1}\ \rightleftarrows\ \nat(\Ic^{n+1}\dd\, \Rc + \Ic^n\dd\, \Ic)\ ,\non
\eeq
where the commutator $[\Ic^n,\Rc]$ is by definition the image of $\Ic^n\dd\Rc$ under the Hochschild operator $b$, and $\Ic^n$ is defined as the unitalized algebra $\Rc^+$ for $n\leq 0$. This is a decreasing filtration because $F^{n+1}_{\Ic}X(\Rc)\subset F^{n}_{\Ic}X(\Rc)$, and for $n<0$ one has $F^{n}_{\Ic}X(\Rc)=X(\Rc)$. Denote by $X_n(\Rc,\Ic)=X(\Rc)/F^{n}_{\Ic}X(\Rc)$ the quotient complex. It is generally not separated. One gets in this way an inverse system of $\zz_2$-graded complexes $\{X_n(\Rc,\Ic)\}_{n\in\zz}$ with projective limit $\Xh(\Rc,\Ic)$. \\

Now suppose that we start from an extension of $m$-algebras $0\to \Ic \to \Rc \to \Ac \to 0$ with continuous linear splitting, and assume that any power $\Ic^n$ is a direct summand in $\Rc$. Then, the sequence $X_n(\Rc,\Ic)$ is related to the $X$-complexes of the quotient $m$-algebras $\Rc/\Ic^n$:
\beq
&&  0 \leftarrow X_0(\Rc,\Ic)=\Ac/[\Ac,\Ac] \leftarrow  X_{1}(\Rc,\Ic) = X(\Ac)\leftarrow \ldots \non\\
&& \ldots\leftarrow  X(\Rc/\Ic^{n-1}) \leftarrow X_{2n-1}(\Rc,\Ic) \leftarrow  X_{2n}(\Rc,\Ic) \leftarrow X(\Rc/\Ic^n) \leftarrow \ldots  \non 
\eeq
Hence the projective limit of the system $\{X_n(\Rc,\Ic)\}_{n\in\zz}$ is isomorphic to the $X$-complex of the pro-algebra $\Rch$: 
\be
\Xh(\Rc,\Ic)= \varprojlim_n X_n(\Rc,\Ic)= \varprojlim_n X(\Rc/\Ic^n)=X(\Rch)\ .
\ee
The pro-complex $\Xh(\Rc,\Ic)$ is naturally filtered by the family of subcomplexes $F^n\Xh(\Rc,\Ic)=\ker(\Xh\to X_n)$. If $0\to \Jc \to \Sc \to \Bc \to 0$ is another extension of $m$-algebras with continuous linear splitting, then the space of linear maps between the two pro-complexes $\Xh(\Rc,\Ic)$ and $\Xh(\Sc,\Jc)$, or between $\Xh$ and $\Xh'$ for short, is given by
\be
\hom(\Xh,\Xh')=\varprojlim_m\Big(\varinjlim_n \hom(X_n,X'_m)\Big)\ ,\label{hom}
\ee
where $\hom(X_n,X'_m)$ is the space of continuous linear maps between the $\zz_2$-graded complexes $X_n(\Rc,\Ic)$ and $X_m(\Sc,\Jc)$. Thus $\hom(\Xh,\Xh')$ is a $\zz_2$-graded complex. It corresponds to the space of linear maps $\{f:\Xh\to\Xh'\ |\ \forall k,\,\exists n:\ f(F^n\Xh)\subset F^k\Xh'\}$; the boundary of an element $f$ of parity $|f|$ is given by the graded commutator $\partial \circ f-(-)^{|f|}f\circ\partial$ with the bounary maps $\d=\nat\dd\oplus\bb$ on $\Xh$ and $\Xh'$. $\hom(\Xh,\Xh')$ itself is filtered by the subcomplexes of linear maps of order $\le n$ for any $n\in\nn$:
\be
\hom^n(\Xh,\Xh')=\{ f:\Xh\to\Xh'\ |\ \forall k,\ f(F^{k+n}\Xh)\subset F^k\Xh'\}\ .\label{homn}
\ee
These $\hom$-complexes will be used in the various definitions of bivariant cyclic cohomology, once the relation between the adic filtration over the $X$-complex of a quasi-free algebra $\Rc$ and the Hodge filtration of the cyclic bicomplex over the quotient algebra $\Ac=\Rc/\Ic$ is established.

\subsection{The tensor algebra}

Taking $\Rc$ as the tensor algebra of an $m$-algebra $\Ac$ provides the link with cyclic homology \cite{Cu1,CQ1}. The (non-unital) tensor algebra $T\Ac$ is the completion of the algebraic direct sum $\bigoplus_{n\ge 1}\Ac^{\hotimes n}$ with respect to the family of seminorms
$$
\widehat{p}=\bigoplus_{n\geq 1}p^{\otimes n}= p \oplus (p\otimes p) \oplus (p\otimes p\otimes p) \oplus \ldots\ ,
$$
where $p$ runs through all the submultiplicative seminorms on $\Ac$. Of course $p^{\otimes n}$ is the projective seminorm on $\Ac^{\otimes n}$ defined by a generalization of (\ref{proj}). These seminorms are submultiplicative with respect to the tensor product $\Ac^{\hotimes n}\times \Ac^{\hotimes m}\to \Ac^{\hotimes n+m}$ and therefore the completion $T\Ac$ is an $m$-algebra. It is free, hence quasi-free: a linear map $\phi: T\Ac\to \Om^2 T\Ac$ with the property $\phi(xy)=\phi(x)y+x\phi(y)+\dd x\dd y$ may be canonically constructed by setting $\phi(a)=0$ on the generators $a\in\Ac\subset T\Ac$, and then recursively $\phi(a_1\otimes a_2)=\dd a_1 \dd a_2$, $\phi(a_1\otimes a_2\otimes a_3)=(\dd a_1\dd a_2)a_3 + \dd(a_1 \otimes a_2)\dd a_3$, and so on...\\
The multiplication map $T\Ac\to\Ac$, sending $a_1\otimes\ldots\otimes a_n$ to the product $a_1\ldots a_n$, is continuous and we denote by $J\Ac$ its kernel. Since the inclusion $\si_{\Ac}:\Ac\to T\Ac$ is a continuous linear splitting of the multiplication map, the two-sided ideal $J\Ac$ is a direct summand in $T\Ac$. This implies a linearly split quasi-free extension $0\to J\Ac\to T\Ac \to \Ac \to 0$. It is the universal free extension of $\Ac$ in the following sense: let $0\to \Ic\to \Rc \to \Ac \to 0$ be any other extension ($\Rc$ is not necessarily quasi-free), provided with a continuous linear splitting $\si:\Ac\to \Rc$. Then one gets a commutative diagram
$$
\xymatrix{
0 \ar[r]  & J\Ac \ar[r] \ar[d]_{\rho_*} & T\Ac  \ar[r] \ar[d]_{\rho_*} & \Ac \ar[r]  \ar@{=}[d] & 0  \\
0 \ar[r] & \Ic \ar[r] & \Rc \ar[r] & \Ac \ar[r] \ar@/_/@{.>}[l]_{\si} & 0 }
$$
where $\rho_*:T\Ac \to \Rc$ is the continuous algebra homomorphism obtained by setting $\rho_*(a)=\si(a)$ on the generators $a\in\Ac\subset T\Ac$. Moreover, the homomorphism $\rho_*$ is independent of the linear splitting $\si$ up to homotopy (two splittings can always be connected by a linear homotopy). \\
As remarked by Cuntz and Quillen \cite{CQ1}, the tensor algebra is closely related to a deformation of the algebra of even-degree noncommutative differential forms $\Om^+\Ac$. Endow the space $\Om^+\Ac$ with the {\it Fedosov product}
\be
\om_1\odot\om_2 :=\om_1\om_2 -d\om_1d\om_2\ ,\quad \om_i\in\Om^+\Ac\ .
\ee
Then $(\Om^+\Ac,\odot)$ is a dense subalgebra of $T\Ac$, with the explicit correspondence 
$$
\Om^+\Ac\ni a_0da_1\ldots da_{2n} \longleftrightarrow a_0\otimes\om(a_1,a_2)\otimes\ldots \otimes\om(a_{2n-1},a_{2n}) \in T\Ac\ . 
$$
It turns out that the Fedosov product $\odot$ extends to the projective limit $\Omh^+\Ac$ and the latter is isomorphic to the pro-algebra
\be
\Th\Ac= \mathop{\limproj}\limits_n T\Ac/(J\Ac)^n\ .
\ee
Moreover, $\Omh\Ac$ and $X(\Th\Ac)=\Xh(T\Ac,J\Ac)$ are isomorphic as $\zz_2$-graded pro-vector spaces \cite{CQ1}, and this isomorphism identifies the Hodge filtration $F^n\Omh\Ac$ with the adic filtration $F^n\Xh(T\Ac,J\Ac)$. By a fundamental result of Cuntz and Quillen, all these identifications are \emph{homotopy equivalences} of pro-complexes, i.e. the boundary $b+B$ on $\Omh\Ac$ corresponds to the boundary $\nat\dd\oplus \bb$ on $\Xh(T\Ac,J\Ac)$ up to homotopy and rescaling (see \cite{CQ1}). Hence the periodic and negative cyclic homologies of $\Ac$ may be computed respectively by $\Xh(T\Ac,J\Ac)$ and $F^n\Xh(T\Ac,J\Ac)$. Also, the non-periodic cyclic homology of $\Ac$ may be computed by the quotient complex $X_n(T\Ac,J\Ac)$ which is homotopy equivalent to the complex $\Omh\Ac/F^n\Omh\Ac$. More generally the same result holds if tensor algebra $T\Ac$ is replaced by any quasi-free extension of $\Ac$. Indeed if $0\to \Ic\to \Rc\to \Ac \to 0$ is a quasi-free extension with continuous linear splitting, the classifying homomorphism $\rho_*:T\Ac\to \Rc$ induces a map of pro-complexes $X(\rho_*): \Xh(T\Ac,J\Ac)\to \Xh(\Rc,\Ic)$ compatible with the adic filtrations induced by the ideals $J\Ac$ and $\Ic$. It turns out to be a homotopy equivalence, irrespective to the choice of $\Rc$:

\begin{theorem}[Cuntz-Quillen \cite{CQ1}]
For any linearly split extension of $m$-algebras $0\to \Ic\to \Rc\to \Ac \to 0$ with $\Rc$ quasi-free, one has isomorphisms 
\beq
HP_n(\Ac) &=& H_{n+2\zz}(\Xh(\Rc,\Ic))\ ,\non\\
HC_n(\Ac) &=& H_{n+2\zz}(X_n(\Rc,\Ic))\ ,\\
HD_n(\Ac) &=& H_{n+2\zz}(X_{n+1}(\Rc,\Ic))\ ,\non\\
HN_n(\Ac) &=& H_{n+2\zz}(F^{n-1}\Xh(\Rc,\Ic))\ .\non
\eeq
\end{theorem}

These filtrations also allow to define various versions of bivariant cyclic cohomology, which may be formulated either within the $X$-complex framework or by means of the $(b+B)$-complex of differential forms.

\begin{definition}[\cite{CQ1}]
Let $\Ac$ and $\Bc$ be $m$-algebras, and choose arbitrary (linearly split) quasi-free extensions $0\to \Ic\to \Rc\to \Ac \to 0$ and $0\to \Jc\to \Sc\to \Bc \to 0$. The bivariant periodic cyclic cohomology of $\Ac$ and $\Bc$ is the homology of the $\zz_2$-graded complex (\ref{hom}) of linear maps between the pro-complexes $\Xh(\Rc,\Ic)$ and $\Xh(\Sc,\Jc)$:
\be
HP^n(\Ac,\Bc)=H_{n+2\zz}\big(\hom(\Xh(\Rc,\Ic),\Xh(\Sc,\Jc))\big)\ ,\quad \forall n\in\zz\ .
\ee
For any $n\in\zz$, the non-periodic cyclic cohomology group $HC^n(\Ac,\Bc)$ of degree $n$ is the homology, in degree $n$ mod $2$, of the $\zz_2$-graded subcomplex (\ref{homn}) of linear maps of order $\le n$:
\be
HC^n(\Ac,\Bc)=H_{n+2\zz}\big(\hom^n(\Xh(\Rc,\Ic),\Xh(\Sc,\Jc))\big)\ .
\ee
The embedding $\hom^n\hookrightarrow \hom^{n+2}$ induces, for any $n$, the $S$-operation in bivariant cyclic cohomology $S:HC^{n}(\Ac,\Bc)\to HC^{n+2}(\Ac,\Bc)$, and $\hom^n\hookrightarrow \hom$ yields a natural map $HC^n(\Ac,\Bc)\to HP^n(\Ac,\Bc)$.
\end{definition}
Of course the bivariant periodic theory has period two: $HP^{n+2}=HP^n$. Let us look at particular cases. The algebra $\cc$ is quasi-free hence $\Xh(T\cc,J\cc)$ is homotopically equivalent to $X(\cc): \cc\rightleftarrows 0$, and the periodic cyclic homology of $\cc$ is simply $HP_0(\cc)=\cc$ and $HP_1(\cc)=0$. This implies that for any $m$-algebra $\Ac$, we get the usual isomorphisms $HP^n(\cc,\Ac)\cong HP_{-n}(\Ac)$ and $HP^n(\Ac,\cc)\cong HP^n(\Ac)$ in any degree $n$. For the non-periodic theory, one has the isomorphism $HC^n(\cc,\Ac)\cong HN_{-n}(\Ac)$ with negative cyclic homology, and $HC^n(\Ac,\cc)\cong HC^n(\Ac)$ is the non-periodic cyclic cohomology of Connes \cite{C0}. \\ 
Finally, since any class $\varphi\in HC^p(\Ac,\Bc)$ is represented by a chain map sending the subcomplex $F^n\Xh(T\Ac,J\Ac)$ to $F^{n-p}\Xh(T\Bc,J\Bc)$ for any $n\in\zz$, it is not difficult to check that $\varphi$ induces a transformation of degree $-p$ between the $SBI$ exact sequences for $\Ac$ and $\Bc$, i.e. a graded-commutative diagram
$$
\vcenter{\xymatrix{
HP_{n+1}(\Ac) \ar[r]^S \ar[d]^{\varphi} & HC_{n-1}(\Ac) \ar[r]^B \ar[d]^{\varphi} & HN_n(\Ac)  \ar[r]^I \ar[d]^{\varphi} & HP_n(\Ac)  \ar[d]^{\varphi} \\
HP_{n-p+1}(\Bc) \ar[r]^S & HC_{n-p-1}(\Bc) \ar[r]^B & HN_{n-p}(\Bc) \ar[r]^I & HP_{n-p}(\Bc) }}
$$
The graded-commutativity comes from the fact that the middle square is actually \emph{anti}commutative when $\varphi$ is of odd degree, for in this case the connecting morphism $B$ anticommutes with the chain map representing $\varphi$.

\section{Quasihomomorphisms and Chern character}\label{sbiv}

In this section we define quasihomomorphisms for \emph{metrizable} (or Fr\'echet) $m$-algebras and construct a bivariant Chern character. The topology of a Fr\'echet $m$-algebra is defined by a countable family of submultiplicative seminorms, and can alternatively be considered as the projective limit of a sequence of Banach algebras \cite{Ph}. In particular, the projective tensor product of two Fr\'echet $m$-algebras is again a Fr\'echet $m$-algebra.\\
We say that a Fr\'echet $m$-algebra $\Ic$ is $p$-summable (with $p\geq 1$ an integer), if there is a continuous trace $\Tr: \Ic^p\to \cc$ on the $p$th power of $\Ic$. Recall that by definition, $\Ic^p$ is the image in $\Ic$ of the $p$-th completed tensor product $\Ic\hotimes\ldots\hotimes \Ic$ under the multiplication map. Hence the trace is understood as a continuous linear map $\Ic\hotimes\ldots\hotimes\Ic\to\cc$, and the tracial property means that it vanishes on the image of $1-\la$, where the operator $\la$ is the backward cyclic permutation $\la(i_1\otimes\ldots\otimes i_p)=i_p\otimes i_1\ldots \otimes i_{p-1}$. In the low degree $p=1$ we interpret the trace just as a linear map $\Ic\to\cc$ vanishing on the subspace of commutators $[\Ic,\Ic]:= b\Om^1\Ic$.\\
Now consider any Fr\'echet $m$-algebra $\Bc$ and form the completed tensor product $\Ic\hotimes\Bc$. Suppose that $\Ec$ is a Fr\'echet $m$-algebra containing $\Ic\hotimes\Bc$ as a (not necessarily closed) two-sided ideal, in the sense that the inclusion $\Ic\hotimes\Bc\to \Ec$ is continuous. The left and right multiplication maps $\Ec\times \Ic\hotimes\Bc \times \Ec\to \Ic\hotimes\Bc$ are then automatically jointly continuous (see \cite{CT}). As in \cite{CuT}, we define the semi-direct sum $\Ec\ltimes \Ic\hotimes\Bc$ as the algebra modeled on the vector space $\Ec\oplus \Ic\hotimes\Bc$, where the product is such that as many elements as possible are put in the summand $\Ic\hotimes\Bc$. The semi-direct sum is a Fr\'echet algebra but it may fail to be multiplicatively convex in general. \emph{The situation when $\Ec\ltimes \Ic\hotimes\Bc$ is a Fr\'echet $m$-algebra} will be depicted as
\be
\Ec\triangleright \Ic\hotimes\Bc
\ee 
to stress the analogy with \cite{Cu}. The definition of quasihomomorphisms involves a $\zz_2$-graded version of $\Ec\triangleright \Ic\hotimes\Bc$, depending only on a choice of parity (even or odd). It is constructed as follows:\\

\noindent{\bf 1) Even case:} Define $\Ec^s_+$ as the Fr\'echet $m$-algebra $\Ec\ltimes \Ic\hotimes\Bc$. It is endowed with a linear action of the group $\zz_2$ by automorphisms: the image of an element $(a,b)\in \Ec\oplus \Ic\hotimes\Bc$ under the generator $F$ of the group is $(a+b,-b)$. We define the $\zz_2$-graded algebra $\Ec^s$ as the crossed product $\Ec^s_+\cp \zz_2$. Hence $\Ec^s$ splits as the direct sum $\Ec^s_+\oplus\Ec^s_-$ where $\Ec^s_+$ is the subalgebra of even degree elements and $\Ec^s_-=F\Ec^s_+$ is the odd subspace. \\
This definition is rather abstract but there is a concrete description of $\Ec^s$ in terms of $2\times 2$ matrices. Consider $M_2(\Ec)=M_2(\cc)\hotimes\Ec$ as a $\zz_2$-graded algebra with grading operator $\bigl(
\begin{smallmatrix}
1 & 0 \\
0 & -1 \end{smallmatrix} \bigr)$. Thus diagonal elements are of even degree and off-diagonal elements are odd. $\Ec^s$ can be identified with a (non-closed) subalgebra of $M_2(\cc)\hotimes\Ec$ in the following way. Any element $x+Fy\in\Ec^s$ may be decomposed to its even and odd parts $x,y\in \Ec\oplus \Ic\hotimes\Bc$, with $x=(a,b)$ and $y=(c,d)$. Then $x+Fy$ is represented by the matrix 
$$
x+Fy=\left(\begin{matrix}
a+b & c \\
c+d & a \end{matrix} \right)\quad \mbox{with} \quad a,c\in \Ec\ ,\ b,d \in \Ic\hotimes\Bc\ .
$$
The action of $\zz_2$ on $\Ec^s_+$ is implemented by the adjoint action of the following odd-degree multiplier of $M_2(\Ec)$:
\be
F=\left(\begin{matrix}
0 & 1 \\
1 & 0 \end{matrix} \right)\in M_2(\cc)\ ,\quad F^2=1\ .
\ee
Denote by $\Ic^s=\Ic^s_+\oplus \Ic^s_-$ the $\zz_2$-graded algebra $M_2(\cc)\hotimes\Ic$, with $\Ic^s_+$ the subalgebra of diagonal elements and $\Ic^s_-$ the off-diagonal subspace. We thus have an inclusion of $\Ic^s\hotimes\Bc$ as a (non-closed) two-sided ideal in $\Ec^s$, with $\Ec^s \triangleright \Ic^s\hotimes\Bc$. The commutator $[F,\Ec^s_+]$ is contained in $\Ic^s_-\hotimes\Bc$. Finally, we denote by $\tr_s$ the supertrace of even degree on $M_2(\cc)$:
$$
\tr_s:M_2(\cc)\to\cc\ ,\qquad \tr_s\left(
\begin{matrix}
a' & c \\
c' & a \end{matrix} \right)=a'-a\ .
$$

\noindent {\bf 2) Odd case:} Now regard $M_2(\Ec)$ as a trivially graded algebra and define $\Ec^s_+$ as the (non-closed) subalgebra
\be
\Ec^s_+=\left(\begin{matrix}
\Ec & \Ic\hotimes\Bc \\
\Ic\hotimes\Bc & \Ec \end{matrix} \right)
\ee
provided with its own topology of complete $m$-algebra. Let $C_1=\cc\oplus \eps\cc$ be the complex Clifford algebra of the one-dimensional euclidian space. $C_1$ is the $\zz_2$-graded algebra generated by the unit $1\in\cc$ in degree zero and $\eps$ in degree one with $\eps^2=1$. We define the $\zz_2$-graded algebra $\Ec^s$ as the tensor product $C_1\hotimes\Ec^s_+$. Hence $\Ec^s=\Ec^s_+\oplus \Ec^s_-$ where $\Ec^s_+$ is the subalgebra of even degree and $\Ec^s_-=\eps\Ec^s_+$ is the odd subspace. Similarly, define $\Ic^s=M_2(C_1)\hotimes \Ic=\Ic^s_+\oplus \Ic^s_-$. Then $\Ic^s\hotimes\Bc$ is a (non-closed) two-sided ideal of $\Ec^s$ and one has $\Ec^s \triangleright \Ic^s\hotimes\Bc$. The matrix
\be
F=\eps\left(\begin{matrix}
1 & 0 \\
0 & -1 \end{matrix} \right) \in M_2(C_1)\ ,\quad F^2=1
\ee 
is an odd multiplier of $\Ec^s$ and the commutator $[F,\Ec^s_+]$ is contained in $\Ic^s_-\hotimes\Bc$. Finally, we define a supertrace $\tr_s$ of odd degree on $C_1$ by sending the generators $1$ to $0$ and $\eps$ to $\pm\sqrt{2i}$. The normalization $\pm\sqrt{2i}$ is chosen for compatibility with Bott periodicity, see \cite{P1}. We will choose conventionally the ``sign'' as $-\sqrt{2i}$ in order to simplify the subsequent formulas. One thus has
$$
\tr_s:M_2(C_1)\to\cc\ ,\quad \tr_s\left(
\begin{matrix}
a+\eps a' & b+\eps b' \\
c+\eps c' & d+\eps d' \end{matrix} \right)=-\sqrt{2i}\, (a'+d')\ .
$$

The objects $F$, $\Ec^s$ and $\Ic^s$ are defined in such a way that we can handle the even and odd case simultaneously. This allows to give the following synthetic definition of quasihomomorphisms.

\begin{definition}
Let $\Ac$, $\Bc$, $\Ic$, $\Ec$ be Fr\'echet $m$-algebras. Assume that $\Ic$ is $p$-summable and $\Ec\triangleright \Ic\hotimes\Bc$. A {\bf quasihomomorphism} from $\Ac$ to $\Bc$ is a continuous homomorphism
\be
\rho: \Ac\to \Ec^s \triangleright \Ic^s\hotimes\Bc 
\ee
sending $\Ac$ to the even degree subalgebra $\Ec^s_+$. The quasihomomorphism comes equipped with a degree (even or odd) depending on the degree chosen for the above construction of $\Ec^s$. In particular, the linear map $a\in\Ac \mapsto [F,\rho(a)]\in \Ic^s_-\hotimes\Bc$ is continuous.
\end{definition}

In other words, a quasihomomorphism of even degree $\rho=\bigl( \begin{smallmatrix} \rho_+ & 0 \\ 0 & \rho_- \end{smallmatrix} \bigr)$ is a pair of homomorphisms $(\rho_+,\rho_-):\Ac\rightrightarrows \Ec$ such that the difference $\rho_+(a)-\rho_-(a)$ lies in the ideal $\Ic\hotimes\Bc$ for any $a\in\Ac$. A quasihomomorphism of odd degree is a homomorphism $\rho:\Ac\to M_2(\Ec)$ such that the off-diagonal elements land in $\Ic\hotimes\Bc$. \\

The Cuntz-Quillen formalism for bivariant cyclic cohomology $HC^n(\Ac,\Bc)$ requires to work with quasi-free extensions of $\Ac$ and $\Bc$. Hence let us suppose that we choose such extensions of $m$-algebras
$$
0\to \Gc \to \Fc \to \Ac \to 0\ ,\qquad 0\to \Jc \to \Rc \to \Bc \to 0\ ,
$$
with $\Fc$ and $\Rc$ quasi-free. We always take $\Fc=T\Ac$ as the universal free extension of $\Ac$, but we leave the possibility to take any quasi-free extension $\Rc$ for $\Bc$ since the tensor algebra $T\Bc$ will not be an optimal choice in general. The first step toward the bivariant Chern character is to lift a given quasihomomorphism $\rho: \Ac\to \Ec^s\triangleright\Ic^s\hotimes\Bc$ to a quasihomomorphism from $\Fc$ to $\Rc$, compatible with the filtrations by the ideals $\Gc\subset \Fc$, $\Jc\subset \Rc$. This requires to fix some admissibility conditions on the intermediate algebra $\Ec$:
\begin{definition}\label{dadm}
Let $0\to \Jc \to \Rc \to \Bc \to 0$ be a quasi-free extension of $\Bc$, and let $\Ic$ be $p$-summable with trace $\Tr:\Ic^p\to \cc$. We say that $\Ec\triangleright \Ic\hotimes\Bc$ is provided with an \emph{$\Rc$-admissible extension} if there are algebras $\Mc\triangleright \Ic\hotimes\Rc$ and $\Nc\triangleright\Ic\hotimes\Jc$ and a commutative diagram of extensions
\be
\vcenter{\xymatrix{
0\ar[r] & \Nc \ar[r]  & \Mc \ar[r] & \Ec \ar[r] & 0 \\
0 \ar[r] & \Ic\hotimes \Jc \ar[r] \ar[u] & \Ic\hotimes \Rc \ar[r] \ar[u] & \Ic\hotimes\Bc \ar[r] \ar[u] & 0}}
\label{adm}
\ee
with the following properties:\\

\noindent i) Any power $\Nc^n$ is a direct summand in $\Mc$ (as a topological vector space), and $\Nc^n\cap \Ic\hotimes\Rc = \Ic\hotimes\Jc^n$;\\

\noindent ii) For any degree $n\geq \max(1,p-1)$, the linear map $(\Ic\hotimes\Rc)^n\dd(\Ic\hotimes\Rc)\to \Om^1\Rc_{\nat}$ induced by the trace $\Ic^{n+1}\to\cc$ factors through the quotient 
$$
\nat (\Ic\hotimes\Rc)^n\dd(\Ic\hotimes\Rc)= (\Ic\hotimes\Rc)^n\dd(\Ic\hotimes\Rc) \mod [\Mc,\Om^1\Mc]\ ,
$$
and the chain map $\Tr: F^{2n+1}_{\Ic\hotimes\Rc}X(\Mc) \to X(\Rc)$ thus obtained is of order zero with respect to the adic filtration induced by the ideals $\Nc\subset \Mc$ and $\Jc\subset\Rc$.
\end{definition}

\noindent In the following we will say that $\Ec$ is $\Rc$-admissible, keeping in mind that the extension $\Mc$ is given. The chain map $\Tr$ of condition {\it ii)} is constructed as follows. For $n\geq 1$ one has the inclusion $\nat (\Ic\hotimes\Rc)^{n+1} \dd \Mc \subset \nat (\Ic\hotimes\Rc)^n\dd(\Ic\hotimes\Rc)$, so that the subcomplex of the $\Ic\hotimes\Rc$-adic filtration reads
$$
F^{2n+1}_{\Ic\hotimes\Rc}X(\Mc)\ :\ (\Ic\hotimes\Rc)^{n+1} \rightleftarrows \nat (\Ic\hotimes\Rc)^n\dd(\Ic\hotimes\Rc)\ .
$$
Then, the trace $\Ic^{n+1}\to\cc$ induces a partial trace $(\Ic\hotimes\Rc)^{n+1}\to \Rc$. The latter combined with $\nat (\Ic\hotimes\Rc)^n\dd(\Ic\hotimes\Rc)\to \Om^1\Rc_{\nat}$ yields a linear map in any degree $n\geq \max(1,p-1)$
\be
\Tr: F^{2n+1}_{\Ic\hotimes\Rc}X(\Mc) \to X(\Rc)\ , \label{tr}
\ee
compatible with the inclusions $F^{2n+3}_{\Ic\hotimes\Rc}X(\Mc)\subset F^{2n+1}_{\Ic\hotimes\Rc}X(\Mc)$. The trace over $\Ic^{n+1}$ ensures that (\ref{tr}) is a chain map. Order zero with respect to the $\Nc$-adic and $\Jc$-adic filtrations means that the intersection $F^{2n+1}_{\Ic\hotimes\Rc}X(\Mc)\cap F^k_{\Nc}X(\Mc)$ is mapped to $F^k_{\Jc}X(\Rc)$ for any $k\in\zz$.\\ 
Remark that the case $p=1$, $n=0$ is pathological, since there is no canonical way to map the space $\nat ((\Ic\hotimes\Rc)\dd \Mc + \Mc^+\dd(\Ic\hotimes\Rc))$ to $\Om^1\Rc_{\nat}$ using only the trace over $\Ic$. In this situation, it seems preferable to impose directly the existence of a chain map $\Tr: F^{1}_{\Ic\hotimes\Rc}X(\Mc) \to X(\Rc)$ in the definition of admissibility.

\begin{example}\label{ek}\textup{When $\Ac$ is arbitrary and $\Bc=\cc$, a $p$-summable quasihomomorphism represents a $K$-homology class of $\Ac$ in the sense of \cite{C0, C1}. Here we take $\Ic=\Lc^p(H)$ as the Schatten ideal of $p$-summable operators on a separable infinite-dimensional Hilbert space $H$. Recall that $\Ic$ is a two-sided ideal in the algebra of all bounded operators $\Lc=\Lc(H)$. $\Ic$ is a Banach algebra for the norm $\| x\|_p=(\Tr(|x|^p))^{1/p}$, $\Lc$ is provided with the operator norm, and the products $\Ic\times \Lc\times\Ic\to\Ic$ are jointly continuous. Since $\Lc$ and $\Ic$ are Banach algebras, the semi-direct sum $\Lc\ltimes \Ic$ is automatically a Banach algebra and we can write $\Lc\triangleright\Ic$. A $p$-summable $K$-homology class of even degree is represented by a pair of continuous homomorphisms $(\rho_+,\rho_-):\Ac\rightrightarrows \Lc$ such that the difference $\rho_+-\rho_-$ lands to $\Ic$. We get in this way an even degree quasihomomorphism $\rho:\Ac\to \Lc^s\triangleright\Ic^s$. Here it is important to note that by a slight modification of the intermediate algebra $\Lc$, it is always possible to consider $\Ic$ as a closed ideal \cite{CuT}. Indeed if we define $\Ec$ as the Banach algebra
$$
\Ec=\Lc\ltimes\Ic
$$
then one clearly has $\Ec\triangleright\Ic$ and $\Ic$ is closed in $\Ec$ by construction. The pair of homomorphisms $(\rho_+,\rho_-):\Ac\rightrightarrows \Lc$ may be replaced with a new pair $(\rho'_+,\rho'_-):\Ac\rightrightarrows \Ec$ by setting $\rho'_+(a)=(\rho_-(a),\rho_+(a)-\rho_-(a))$ and $\rho'_-(a)=(\rho_-(a),0)$ in $\Lc\oplus\Ic$. The above $K$-homology class is then represented by the new quasihomomorphism $\rho':\Ac\to\Ec^s\triangleright\Ic^s$.\\
In the odd case, a $p$-summable $K$-homology class is represented by a continuous homomorphism $\rho:\Ac\to \bigl(\begin{smallmatrix}
\Lc & \Ic \\
\Ic & \Lc \end{smallmatrix} \bigr)$, which can be equivalently described as a homomorphism $\rho':\Ac\to M_2(\Ec)$ with off-diagonal elements in $\Ic$.\\
Concerning cyclic homology, the algebra $\cc$ is quasi-free, hence the quasi-free extension $\Rc=\cc$ and $\Jc=0$ computes the cyclic homology of $\cc$. Therefore by choosing $\Mc=\Ec$ and $\Nc=0$, the algebra $\Ec\triangleright \Ic$ is $\cc$-admissible  (condition {\it ii)} is trivial since $\Om^1\cc_{\nat}=0$).}
\end{example}

\begin{example}\label{ebiv}\textup{More generally, if $\Ic$ is a $p$-summable Fr\'echet $m$-algebra contained as a (not necessarily closed) two-sided ideal in a unital Fr\'echet $m$-algebra $\Lc$, with $\Lc\triangleright\Ic$, a $p$-summable quasihomomorphism from $\Ac$ to $\Bc$ could be constructed from the generic intermediate algebra $\Ec=\Lc\hotimes\Bc$, provided that the map $\Ic\hotimes\Bc\to \Ec$ is injective. If $0\to \Jc \to \Rc \to \Bc \to 0$ is a quasi-free extension of $\Bc$, the choice $\Mc=\Lc\hotimes\Rc$ and $\Nc=\Lc\hotimes\Jc$ shows that $\Ec$ is $\Rc$-admissible provided that the maps $\Ic\hotimes\Rc\to \Mc$ and $\Ic\hotimes\Jc\to \Nc$ are injective. In fact it is easy to get rid of these injectivity conditions by redefining the algebra
$$
\Ec=(\Lc\ltimes\Ic)\hotimes\Bc\ ,
$$
which contains $\Ic\hotimes\Bc$ as a closed ideal. Then $\Ec$ becomes automatically $\Rc$-admissible by taking $\Mc=(\Lc\ltimes\Ic)\hotimes\Rc$ and $\Nc=(\Lc\ltimes\Ic)\hotimes\Jc$ (remark that $(\Lc\ltimes\Ic)^n=\Lc\ltimes\Ic$ for any $n$ because $\Lc$ is unital, hence $\Nc^n$ is a direct summand in $\Mc$). The chain map $\Tr:F^{2n+1}_{\Ic\hotimes\Rc}X(\Mc) \to X(\Rc)$ is obtained by multiplying all the factors in $\Lc$ and $\Ic$, and taking the trace on $\Ic^{n+1}$, while the compatibility between the $\Nc$-adic and $\Jc$-adic filtrations is obvious. Although interesting examples arise under this form (see section \ref{sassem}), the algebra $\Ec$ cannot be decomposed into a tensor product with $\Bc$ in all situations.}
\end{example}

\begin{example}\label{ebott}\textup{ An important example where $\Ec$ cannot be taken in the previous form is provided by the Bott element of the real line. Here $\Ac=\cc$ and $\Bc=\cinf(0,1)$ is the algebra of smooth functions $f:[0,1]\to \cc$ such that $f$ and all its derivatives vanish at the endpoints $0$ and $1$. Take $\Ic=\cc$ as a 1-summable algebra, and $\Ec=\cinf[0,1]$ is the algebra of smooth functions $f:[0,1]\to\cc$ with the derivatives vanishing at the endpoints, while $f$ itself takes arbitrary values at $0$ and $1$. $\Bc$ and $\Ec$ provided with their usual Fr\'echet topology are $m$-algebras, and one has $\Ec\triangleright\Bc$. The Bott element is represented by the quasihomomorphism of \emph{odd degree} 
$$
\rho:\cc \to \Ec^s\triangleright \Ic^s\hotimes\Bc\ ,
$$
where $\Ic^s=M_2(C_1)$ and $\Ec^s\subset M_2(C_1)\hotimes \Ec$ by construction. The homomorphism $\rho:\cc\to\Ec^s_+$ is built from an arbitrary real-valued function $\xi\in \Ec$ such that $\xi(0)=0$, $\xi(1)=\pi/2$, and sends the unit $e\in\cc$ to the matrix
$$
\rho(e)=R^{-1}\left( \begin{matrix}
 e & 0 \\
 0 & 0 \end{matrix} \right)R\ ,\qquad R= \left( \begin{matrix}
 \cos\xi & \sin\xi \\
 -\sin\xi & \cos\xi \end{matrix} \right)\ .
$$
The algebra $\Bc$ is quasi-free hence we can choose $\Rc=\Bc$, $\Jc=0$ as quasi-free extension. The cyclic homology of $\Bc$ is therefore computed by $X(\Bc)$. Moreover, setting $\Mc=\Ec$ and $\Nc=0$ shows that $\Ec$ is $\Bc$-admissible. Indeed, $\Om^1\Bc_{\nat}$ is contained in the space $\Om^1(0,1)$ of ordinary (commutative) complex-valued smooth one-forms over $[0,1]$ vanishing at the endpoints with all their derivatives. The chain map $\Tr:F^{2n+1}_{\Bc}X(\Ec) \to X(\Bc)$ is thus well-defined for any $n\geq 1$, and just amounts to project noncommutative forms over $\Ec$ to ordinary (commutative) differential forms over $[0,1]$. }
\end{example}

We shall now construct the bivariant Chern character of a given $p$-summable quasihomomorphism $\rho: \Ac\to \Ec^s \triangleright \Ic^s\hotimes\Bc$. We take the universal free extension $0\to J\Ac \to T\Ac \to \Ac \to 0$ for $\Ac$, and choose some quasi-free extension $0\to \Jc \to \Rc \to \Bc \to 0$ for $\Bc$ with the property that the algebra $\Ec\triangleright \Ic\hotimes\Bc$ is $\Rc$-admissible. The bivariant Chern character should be represented by a chain map between the complexes $X(T\Ac)$ and $X(\Rc)$, compatible with the adic filtrations induced by the ideals $J\Ac$ and $\Jc$ (section \ref{scy}). Our task is thus to lift the quasihomomorphism to the quasi-free algebras $T\Ac$ and $\Rc$. First, the admissibility condition provides a diagram of extensions (\ref{adm}). From $\Mc\triangleright \Ic\hotimes\Rc$ define the $\zz_2$-graded algebra $\Mc^s$ in complete analogy with $\Ec^s$: depending on the degree of the quasihomomorphism, $\Mc^s$ is a subalgebra of $M_2(\cc)\hotimes\Mc$ (even case) or $M_2(C_1)\hotimes\Mc$ (odd case), with commutator property $[F,\Mc^s_+]\subset \Ic^s_-\hotimes\Rc$. Also, from $\Nc\triangleright \Ic\hotimes\Jc$ define $\Nc^s$ as the $\zz_2$-graded subalgebra of $M_2(\cc)\hotimes\Nc$ or $M_2(C_1)\hotimes\Nc$ with commutator $[F,\Nc^s_+]\subset \Ic^s_-\hotimes\Jc$. The algebras $\Ec^s$, $\Mc^s$ and $\Nc^s$ are gifted with a differential of odd degree induced by the graded commutator $[F,\ ]$ (its square vanishes because $F^2=1$). Then we get an extension of $\zz_2$-graded differential algebras
$$
0\to \Nc^s \to \Mc^s \to \Ec^s \to 0\ .
$$
The restriction to the even-degree subalgebras yields an extension of trivially graded algebras $0\to \Nc^s_+ \to \Mc^s_+ \to \Ec^s_+ \to 0$, split by a continuous linear map $\si:\Ec^s_+\to \Mc^s_+$ (recall the splitting is our basic hypothesis about extensions of $m$-algebras). The universal properties of the tensor algebra $T\Ac$ then allows to extend the homomorphism $\rho: \Ac\to \Ec^s_+$ to a continuous homomorphism $\rho_*:T\Ac\to \Mc^s_+$ by setting $\rho_*(a_1\otimes\ldots\otimes a_n)=\si\rho(a_1)\otimes\ldots\otimes \si\rho(a_n)$:
\be
\vcenter{\xymatrix{
0 \ar[r]  & J\Ac \ar[r] \ar[d]_{\rho_*} & T\Ac  \ar[r] \ar[d]_{\rho_*} & \Ac \ar[r]  \ar[d]^{\rho} & 0  \\
0 \ar[r] & \Nc^s_+ \ar[r] & \Mc^s_+ \ar[r] & \Ec^s_+ \ar[r] \ar@/_/@{.>}[l]_{\si} & 0 }} \label{uni}
\ee
A priori $\rho_*$ depends on the choice of linear splitting $\si$, but in a way which will not affect the cohomology class of the bivariant Chern character. This construction may be depicted in terms of a $p$-summable quasihomomorphism $\rho_*: T\Ac \to \Mc^s\triangleright \Ic^s\hotimes \Rc$, compatible with the adic filtration by the ideals in the sense that $J\Ac$ is mapped to $\Nc^s\triangleright \Ic^s\hotimes\Jc$. Hence, $\rho_*$ extends to a quasihomomorphism of pro-algebras
\be
\rho_* : \Th\Ac \to \Mch^s\triangleright \Ic^s\hotimes \Rch\ ,
\ee
where $\Th\Ac$, $\Mch^s$ and $\Rch$ are the adic completions of $T\Ac$, $\Mc^s$ and $\Rc$ with respect to the ideals $J\Ac$, $\Nc^s$ and $\Jc$. Next, depending on the degree of the quasihomomorphism, the even supertrace $\tr_s:M_2(\cc)\to \cc$ or the odd supertrace $\tr_s:M_2(C_1)\to\cc$ yields a chain map $X(\Mc^s)\to X(\Mc)$ by setting $\al x\mapsto \tr_s(\al)x$ and $\nat\, \al x\dd(\beta y)\mapsto \pm\tr_s(\al\beta)\nat x\dd y$ for any $x,y\in \Mc$ and $\al,\beta\in M_2(\cc)$ or $M_2(C_1)$. The sign $\pm$ is the parity of the matrix $\beta$, which has to move across the differential $\dd$. Hence composing with the chain map $\Tr:F^{2n+1}_{\Ic\hotimes\Rc}X(\Mc)\to X(\Rc)$ guaranteed by the admissibility condition, we obtain for any integer $n\geq \max(1,p-1)$ a {\it supertrace} chain map
\be
\tau: F^{2n+1}_{\Ic^s\hotimes\Rc}X(\Mc^s)\stackrel{\tr_s}{\to} F^{2n+1}_{\Ic\hotimes\Rc}X(\Mc) \stackrel{\Tr}{\to} X(\Rc)\ ,
\ee
of order zero with respect to the $\Nc^s$-adic filtration on $X(\Mc^s)$ and the $\Jc$-adic filtration on $X(\Rc)$. The parity of $\tau$ corresponds to the parity of the quasihomomorphism. This allows to construct a chain map $\chih^n:\Omh \Mc^s_+\to X(\Rc)$ from the $(b+B)$-complex of the algebra $\Mc^s_+$, in any degree $n\geq p$ having the same parity as the supertrace $\tau$. Observe that the linear map $x\in \Mc^s_+ \mapsto [F,x] \in\Ic^s_-\hotimes\Rc$ is continuous by construction.

\begin{proposition}\label{pchi}
Let $\rho:\Ac\to \Ec^s\triangleright \Ic^s\hotimes\Bc$ be a $p$-summable quasihomomorphism of parity $p\mod 2$, with $\Rc$-admissible algebra $\Ec$. Given any integer $n\geq p$ of the same parity, consider two linear maps $\chih^n_0:\Om^n \Mc^s_+\to \Rc$ and $\chih^n_1:\Om^{n+1}\Mc^s_+\to \Om^1\Rc_{\nat}$ defined by
\be
\chih^n_0(x_0\dd x_1\ldots\dd x_n) = (-)^n\frac{\Gamma(1+\frac{n}{2})}{(n+1)!} \sum_{\la\in S_{n+1}} \eps(\la)\, \tau(x_{\la(0)}[F,x_{\la(1)}]\ldots [F,x_{\la(n)}])\label{chi}
\ee
$$
\chih^n_1(x_0\dd x_1\ldots\dd x_{n+1}) = (-)^n\frac{\Gamma(1+\frac{n}{2})}{(n+1)!} \sum_{i=1}^{n+1}  \tau\nat(x_0[F,x_1]\ldots\dd x_i \ldots [F,x_{n+1}])
$$
where $S_{n+1}$ is the cyclic permutation group of $n+1$ elements and $\eps$ is the signature. Then $\chih^n_0$ and $\chih^n_1$ define together a chain map $\chih^n:\Omh \Mc^s_+\to X(\Rc)$ of parity $n \mod 2$, i.e. fulfill the relations
\be
\chih^n_0 B = 0\ ,\quad \nat\dd\chih^n_0-(-)^n\chih^n_1B=0\ ,\quad \bb\chih^n_1-(-)^n\chih^n_0b=0\ ,\quad \chih^n_1b=0\ . \label{rchi}
\ee
Moreover $\chih^n$ is invariant under the Karoubi operator $\kappa$ acting on $\Om^n\Mc^s_+$ and $\Om^{n+1}\Mc^s_+$.
\end{proposition}
{\it Proof:} This follows from purely algebraic manipulations, using the following general properties:\\
- The graded commutator $[F,\ ]$ is a differential and $\tau([F,\ ])=0$;\\
- $\dd F=0$ so that $[F,\ ]$ and $\dd$ are anticommuting differentials;\\
- $\tau\nat$ is a supertrace.\\
The computation is lengthy but straightforward. \cqfd \\

Thus we have attached to a $p$-summable quasihomomorphism $\rho:\Ac\to\Ec^s\triangleright\Ic^s\hotimes\Bc$ of parity $p\mod 2$ a sequence of cocycles $\chih^n$ ($n\geq p$) of the same parity in the $\zz_2$-graded complex $\hom(\Omh \Mc^s_+,X(\Rc))$. They are in fact all cohomologous, and the proposition below gives an explicit transgression formula in terms of the \emph{eta-cochain}:
\begin{proposition}\label{peta}
Let $\rho:\Ac\to \Ec^s\triangleright \Ic^s\hotimes\Bc$ be a $p$-summable quasihomomorphism of parity $p\mod 2$, with $\Rc$-admissible algebra $\Ec$. Given any integer $n\geq p+1$ of parity \emph{opposite} to $p$, consider two linar maps $\etah^n_0:\Om^n\Mc^s_+\to \Rc$ and $\etah^n_1:\Om^{n+1}\Mc^s_+\to \Om^1\Rc_{\nat}$ defined by
\beq
\etah^n_0(x_0\dd x_1\ldots\dd x_n) &=& \frac{\Gamma(\frac{n+1}{2})}{(n+1)!} \, \frac{1}{2}\tau\Big(F x_0[F,x_1]\ldots [F,x_n] + \non\\
&&\sum_{i=1}^n (-)^{ni} [F,x_i]\ldots [F,x_n]Fx_0 [F,x_1]\ldots [F,x_{i-1}] \Big) \non
\eeq
\beq
\lefteqn{\etah^n_1(x_0\dd x_1\ldots\dd x_{n+1}) =}\label{eta}\\
&& \frac{\Gamma(\frac{n+1}{2})}{(n+2)!} \sum_{i=1}^{n+1}  \frac{1}{2}\tau\nat(ix_0 F + (n+2-i)Fx_0)[F,x_1]\ldots\dd x_i \ldots [F,x_{n+1}]\ .\non
\eeq
Then $\etah^n_0$ and $\etah^n_1$ define together a cochain $\etah^n\in \hom(\Omh \Mc^s_+,X(\Rc))$ of parity $n \mod 2$, whose coboundary equals the difference of cocycles
$$
\chih^{n-1}-\chih^{n+1}= (\nat\dd\oplus \bb) \etah^n - (-)^n \etah^n (b+B)\ .
$$
Expressed in terms of components this amounts to the identities
\beq
\chih^{n-1}_0 &=& -(-)^n\etah^n_0 B\ ,\quad\qquad  \bb\etah^n_1 -(-)^n(\etah^n_0b+\etah^{n+2}_0 B)=0\ , \label{reta}\\
\chih^{n-1}_1 &=& \nat\dd \etah^n_0 -(-)^n\etah^n_1 B\ ,\quad \nat\dd\etah^{n+2}_0 -(-)^n(\etah^n_1b+\etah^{n+2}_1 B)=0\ .\non
\eeq
\end{proposition}
{\it Proof:}  Direct computation. \cqfd\\

\begin{remark}\label{rsum}\textup{Using a trick of Connes \cite{C0}, we may replace the chain map $\tau$ by $\tau'=\frac{1}{2}\tau(F[F,\ ])$. This allows to improve the summability condition by requiring the quasihomomorphism to be only $(p+1)$-summable instead of $p$-summable, while the condition on the degree remains $n\geq p$ for $\chih^n$ and $n\geq p+1$ for $\etah^n$. It is traightforward to write down the new formulas for $\chih^n$ and observe that it involves exactly $n+1$ commutators $[F,x]$. These formulas were actually obtained in \cite{P2} in a more general setting where we allow $\dd F \neq 0$.}
\end{remark}

\begin{definition}
The bivariant Chern character of the quasihomomorphism $\rho:\Ac\to\Ec^s\triangleright\Ic^s\hotimes\Bc$ is represented in any degree $n\geq p$ by the composition of chain maps 
\be
\ch^n(\rho): X(T\Ac) \stackrel{\gamma}{\longrightarrow} \Omh T\Ac \stackrel{\rho_*}{\longrightarrow} \Omh\Mc^s_+ \stackrel{\chih^n}{\longrightarrow} X(\Rc)\ ,
\ee
where $\gamma:X(T\Ac)\to \Omh T\Ac$ is the Goodwillie equivalence constructed in section \ref{scy} for any quasi-free algebra and $\rho_*:T\Ac\to \Mc^s_+$ is the classifying homomorphism. In the same way, define a transgressed cochain in $\hom(X(T\Ac),X(\Rc))$ by means of the eta-cochain in any degree:
\be
\tch^n(\rho) : X(T\Ac) \stackrel{\gamma}{\longrightarrow} \Omh T\Ac \stackrel{\rho_*}{\longrightarrow} \Omh\Mc^s_+ \stackrel{\etah^n}{\longrightarrow} X(\Rc)\ .
\ee
It fulfills the transgression property $\ch^n(\rho)-\ch^{n+2}(\rho)=[\d, \tch^{n+1}(\rho)]$ where $\d$ is the $X$-complex boundary map.
\end{definition}

Recall that $\gamma(x)=(1-\phi)^{-1}(x)$ and $\gamma(\nat x\dd y)=(1-\phi)^{-1}(x\dd y+b(x\phi(y)))$ for any $x,y\in T\Ac$, where the map $\phi: \Om^nT\Ac\to \Om^{n+2}T\Ac$ is uniquely defined from its restriction to zero-forms. Its existence is guaranteed by the fact that $T\Ac$ is a free algebra. Several choices are possible, but conventionally we always take $\phi:T\Ac\to \Om^2T\Ac$ as the canonical map obtained by setting $\phi(a)=0$ on the generators $a\in \Ac\subset T\Ac$, and then extended to all $T\Ac$ by the algebraic property $\phi(xy)=\phi(x)y+x\phi(y)+\dd x\dd y$.\\
Of course $\ch^n(\rho)$ and $\tch^n(\rho)$ are not very interesting a priori, because the $X$-complex of the non-completed tensor algebra $T\Ac$ is contractible. However, taking into account the adic filtrations induced by the ideals $J\Ac\subset T\Ac$ and $\Jc\subset \Rc$ yields non-trivial bivariant objects. By virtue of Remark \ref{rsum} we suppose from now on that $\Ic$ is $(p+1)$-summable.

\begin{proposition}\label{padic}
Let $\rho:\Ac\to\Ec^s\triangleright\Ic^s\hotimes\Bc$ be a $(p+1)$-summable quasihomomorphism of parity $p \mod 2$ with $\Rc$-admissible extension $\Ec$, and let $n\geq p$ be an integer of the same parity. The composites $\chih^n\rho_*\gamma$ and $\etah^{n+1}\rho_*\gamma$ are linear maps $X(T\Ac)\to X(\Rc)$ verifying the adic properties
\beq
\chih^n\rho_*\gamma &:& F^k_{J\Ac}X(T\Ac) \to F^{k-n}_{\Jc}X(\Rc)\ , \non\\
\etah^{n+1}\rho_*\gamma &:& F^k_{J\Ac}X(T\Ac) \to F^{k-n-2}_{\Jc}X(\Rc)\ ,\non
\eeq
for any $k\in\zz$. Consequently the composite $\ch^n(\rho)=\chih^n\rho_*\gamma$ defines a cocycle of parity $n \mod 2$ in the complex $\hom^n(\Xh(T\Ac,J\Ac),\Xh(\Rc,\Jc))$ and the Chern character is a bivariant cyclic cohomology class of degree $n$:
\be
\ch^n(\rho)\in HC^n(\Ac,\Bc)\ ,\quad \forall n\geq p\ .
\ee
Moreover, the transgression relation $\ch^n(\rho)-\ch^{n+2}(\rho)=[\d, \tch^{n+1}(\rho)]$ holds in the complex $\hom^{n+2}(\Xh(T\Ac,J\Ac),\Xh(\Rc,\Jc))$, which implies
\be
\ch^{n+2}(\rho) \equiv S\ch^n(\rho) \ \mbox{in}\ HC^{n+2}(\Ac,\Bc)\ .
\ee
In particular the cocycles $\ch^n(\rho)$ for different $n$ define the same periodic cyclic cohomology class $\ch(\rho)\in HP^n(\Ac,\Bc)$.
\end{proposition}
{\it Proof:} Let us denote by $0\to \Gc \to \Fc \to \Ac \to 0$ the universal extension $ 0\to J\Ac \to T\Ac \to \Ac \to 0$. Recall that $\Mc^s$ and its ideal $\Nc^s$ are $\zz_2$-graded differential algebras on which the graded commutator $[F,\ ]$ acts as a differential of odd degree. Moreover the commutation relations  $[F,\Mc^s_+]\subset \Ic^s_-\hotimes\Rc$ and $[F,\Nc^s_+]\subset \Ic^s_-\hotimes\Jc$ hold. Now, we have to investigate the adic behaviour of the Goodwillie equivalence $\gamma: X(\Fc)\to \Omh\Fc$ with respect to the filtration $F^k_{\Gc}X(\Fc)$. The first step in that direction was actually done in the proof of Proposition \ref{pgood}, where the following filtration of the subspaces $\Om^n\Fc$ was introduced:
$$
H^k\Om^n\Fc = \sum_{k_0+\ldots +k_n \geq k}\Gc^{k_0} \dd \Gc^{k_1} \ldots \dd\Gc^{k_n}\ \subset \Om^n\Fc\ .
$$
Let us look at the image of the latter filtration under the maps $\chih^n\rho_*$ and $\etah^n\rho_*:\Omh\Fc\to X(\Rc)$ given by Eqs. (\ref{chi}, \ref{eta}). We know that the homomorphism $\rho_*:\Fc\to \Mc^s_+$ respects the ideals $\Gc$ and $\Nc^s_+$. Hence if $x_0, \ldots,x_n$ denote $n+1$ elements in $\Gc^{k_0},\ldots,\Gc^{k_n}$ respectively, with $k_0+\ldots+k_n\geq k$, then $x_{\la(0)}[F,x_{\la(1)}]\ldots [F,x_{\la(n)}]\in (\Nc^s)^k$ for any permutation $\la\in S_{n+1}$. Hence applying the supertrace $\tau$, which is a chain map of order zero with respect to the $\Nc^s$-adic and $\Jc$-adic filtrations on $X(\Mc^s)$ and $X(\Rc)$, yields (from now on we omit to write the homomorphism $\rho_*$)
\be
\chih^n_0(H^k\Om^n\Fc)\subset \Jc^k \ . \label{1}
\ee
In the same way, for $n+1$ elements $x_0, \ldots,x_{n+1}$ in $\Gc^{k_0},\ldots,\Gc^{k_{n+1}}$, the one-form $\nat x_0[F,x_1]\ldots \dd x_i \ldots [F,x_{n+1}]$ involves $k_0+\ldots+k_{n+1}\geq k$ powers of the ideal $\Nc^s$, hence lies in the subspace $\nat((\Nc^s)^k \dd\Mc^s+(\Nc^s)^{k-1}\dd\Nc^s)$. Thus applying the supertrace $\tau$ one gets
\be
\chih^n_1(H^k\Om^{n+1}\Fc)\subset \nat(\Jc^k\dd\Rc+ \Jc^{k-1}\dd\Jc)\ .  \label{2}
\ee
Proceeding in exactly the same fashion with the maps $\etah^n_0:\Om^n\Fc\to \Rc$ and $\etah^n_1:\Om^{n+1}\Fc\to \Om^1\Rc_{\nat}$, it is clear that
\beq
\etah^n_0(H^k\Om^n\Fc) &\subset& \Jc^k \ . \label{3} \\
\etah^n_1(H^k\Om^{n+1}\Fc) &\subset& \nat(\Jc^k\dd\Rc+ \Jc^{k-1}\dd\Jc)\ .  \label{4}
\eeq
However these estimates are not optimal concerning the component $\chih^n_1$. We need a refinement of the $H$-filtration. For any $k\in\zz$, $n\geq 0$, let us define the subspaces 
$$
G^k\Om^n\Fc = \sum_{k_0+\ldots +k_n \geq k}\Gc^{k_0} (\dd\Fc) \Gc^{k_1} (\dd\Fc)\ldots \Gc^{k_{n-1}}(\dd \Fc) \Gc^{k_n} + H^{k+1}\Om^n\Fc\ .
$$
Then for fixed $n$, $G^*\Om^n\Fc$ is a decreasing filtration of $\Om^n\Fc$, and by convention $G^k\Om^n\Fc=\Om^n\Fc$ for $k\leq 0$. One has $G^k\Om^n\Fc \subset H^k\Om^n\Fc$. Now observe the following. Since $[F,\ ]$ and $\dd$ are derivations, the map $x_0\dd x_1\ldots \dd x_i \ldots\dd x_{n+1} \mapsto \nat x_0[F,x_1]\ldots \dd x_i \ldots [F,x_{n+1}]$ has the property that
\beq
\lefteqn{\Gc^{k_0} (\dd\Fc) \Gc^{k_1} (\dd\Fc) \ldots (\dd \Fc) \Gc^{k_i}\ldots (\dd \Fc) \Gc^{k_{n+1}} \to } \non\\
&&\nat (\Nc^s)^{k_0} [F,\Mc^s] (\Nc^s)^{k_1} [F,\Mc^s] \ldots (\dd\Mc^s)(\Nc^s)^{k_i}\ldots [F,\Mc^s] (\Nc^s)^{k_{n+1}} \non\\
&& \subset \nat(\Nc^s)^k\dd\Mc^s\ ,\non
\eeq
and because $\chih^n_1(H^{k+1}\Om^{n+1}\Fc)\subset \nat(\Jc^{k+1}\dd\Rc+ \Jc^k\dd\Jc)\subset \nat \Jc^k\dd\Rc$, one gets the crucial estimate
\be
\chih^n_1(G^k\Om^{n+1}\Fc)\subset \nat \Jc^k\dd\Rc\ .  \label{5}
\ee
Now we have to understand the way $\gamma$ sends the $X$-complex filtration
\beq 
F_{\Gc}^{2k}X(\Fc) &:& \Gc^{k+1}+[\Gc^k,\Fc] \ \rightleftarrows \ \nat \Gc^k\dd \Fc \non\\
F_{\Gc}^{2k+1}X(\Fc) &:& \Gc^{k+1}\ \rightleftarrows\ \nat(\Gc^{k+1}\dd \Fc + \Gc^k\dd \Gc)\ ,\non
\eeq 
to the filtration $G^*\Om^n\Fc$, in all degrees $n$. Recall that $\gamma(x)|_{\Om^{2n}\Fc}=\phi^n(x)$ and $\gamma(\nat x\dd y)|_{\Om^{2n+1}\Fc}=\phi^n(x\dd y + b(x\phi(y)))$ for any $x,y\in \Fc$, where the map $\phi: \Om^n\Fc\to \Om^{n+2}\Fc$ is obtained from its restriction to zero-forms as
$$
\phi(x_0\dd x_1\ldots \dd x_n)=\sum_{i=0}^n(-)^{ni}\phi(x_i)\dd x_{i+1}\ldots \dd x_n\dd x_0\ldots \dd x_{i-1}\ .
$$
Note the following important properties of $\phi$. Firstly, it is invariant under the Karoubi operator $\kappa:\Om^n\Fc\to\Om^n\Fc$ in the sense that $\phi\circ\kappa =\phi$, and vanishes on the image of the boudaries $d,B: \Om^n\Fc\to \Om^{n+1}\Fc$. Secondly, the relation $\phi b - b\phi=B$ holds on $\Om^n\Fc$ whenever $n\geq 1$ (see \cite{P1} \S 4). Since we have to apply successive powers of $\phi$ on the filtration $G^k\Om^n\Fc$, the computation will be greatly simplified by exploiting $\kappa$-invariance. Define the linear map $\phit:\Om^n\Fc\to\Om^{n+2}\Fc$ by 
$$
\phit(x_0\dd x_1\ldots \dd x_n)=\sum_{i=0}^n(-)^{i}\dd x_0 \ldots \dd x_{i-1}\phi(x_i)\dd x_{i+1}\ldots \dd x_n\ .
$$
Then $\phit$ coincides with $\phi$ modulo the image of $1-\kappa$. In particular the relation $\phi^n=\phi\circ \phit^{n-1}$ holds. The advantage of the map $\phit$ stems from the fact that it does not involve cyclic permutations of the elements $x_i$, and verifies the following optimal compatibility with the $G$-filtration
$$
\phit(G^k\Om^n\Fc)\subset G^{k-1}\Om^{n+2}\Fc \qquad \forall k,n \geq 0\ ,
$$
whereas the map $\phi$ is only compatible with the (coarser) $H$-filtration:
$$
\phi(H^k\Om^n\Fc)\subset H^{k-1}\Om^{n+2}\Fc \qquad \forall k,n \geq 0\ .
$$
We shall now evaluate the image of the filtration $F^k_{\Gc}X(\Fc)$ under the map $\gamma:X(\Fc)\to \Omh\Fc$. Firstly, one has $\gamma(\Gc^{k+1})\cap \Om^{2n}\Fc=\phi^n(\Gc^{k+1})$. But $\Gc^{k+1}\subset G^{k+1}\Om^0\Fc$ and $\phi^n=\phi\circ \phit^{n-1}$, hence
\be
\gamma(\Gc^{k+1})\cap \Om^{2n}\Fc \subset \phi(G^{k-n+2}\Om^{2n-2}\Fc)\ .\label{a}
\ee
Secondly, the image of $\nat\Gc^k\dd\Fc$ in $\Om^{2n+1}\Fc$ is given by $\phi^n(\Gc^k\dd\Fc+b(\Gc^k\phi(\Fc)))$. One has $b(\Gc^k\phi(\Fc))\subset [\Gc^k\dd\Fc,\Fc] \subset \Gc^k\dd\Fc$, hence we only need to compute $\phi^n(\Gc^k\dd\Fc)\subset \phi^n(G^k\Om^1\Fc)$, and
\be
\gamma(\nat \Gc^k\dd\Fc)\cap \Om^{2n+1}\Fc \subset \phi(G^{k-n+1}\Om^{2n-1}\Fc)\ . \label{b}
\ee
Thirdly, $[\Gc^k,\Fc]=\bb\nat\Gc^k\dd\Fc$ so that $\gamma([\Gc^k,\Fc])=(b+B)\gamma(\nat\Gc^k\dd\Fc)$ because $\gamma$ is a chain map. Therefore, the image of $[\Gc^k,\Fc]$ restricted to $\Om^{2n}\Fc$ is contained in $B\phi^{n-1}(\Gc^k\dd\Fc) + b\phi^n(\Gc^k\dd\Fc)$. We may estimate coarsly the first term as $B\phi^{n-1}(\Gc^k\dd\Fc) \subset B\phi^{n-1}(H^k\Om^1\Fc) \subset H^{k-n+1}\Om^{2n}\Fc$, and the second term as $b\phi \phit^{n-1}(\Gc^k\dd\Fc) \subset b\phi(G^{k-n+1}\Om^{2n-1}\Fc)$. Hence
\be
\gamma([\Gc^k,\Fc])\cap \Om^{2n}\Fc \subset H^{k-n+1}\Om^{2n}\Fc + b\phi(G^{k-n+1}\Om^{2n-1}\Fc)\ . \label{c}
\ee
Fourthly, the image of $\nat\Gc^k\dd\Gc$ in $\Om^{2n+1}\Fc$ is given by $\phi^n(\Gc^k\dd\Gc+b(\Gc^k\phi(\Gc)))$. We estimate coarsly $\phi^n(\Gc^k\dd\Gc)\subset \phi^n(H^{k+1}\Om^1\Fc)\subset H^{k-n+1}\Om^{2n+1}\Fc$. Then, one has $\Gc^k\phi(\Gc)\subset \Gc^k\dd\Fc\dd\Fc \subset G^k\Om^2\Fc$, and using repeatedly the relations $\phi b-b\phi=B$, $\phi B=0$ gives $\phi^nb(\Gc^k\phi(\Gc)) \subset b\phi^n(G^k\Om^2\Fc)+B\phi^{n-1}(G^k\Om^2\Fc) \subset b\phi(G^{k-n+1}\Om^{2n}\Fc) + BH^{k-n+1}\Om^{2n}\Fc$. Thus
\be
\gamma(\nat \Gc^k\dd\Gc)\cap \Om^{2n+1}\Fc \subset H^{k-n+1}\Om^{2n+1}\Fc + b\phi(G^{k-n+1}\Om^{2n}\Fc)\ . \label{d}
\ee
Now everything is set to evaluate the adic behaviour of the composites $\chih^n\gamma$ and $\etah^n\gamma$. We shall deal only with even degrees, the odd case is similar. Hence let us start with the map $\chih^{2n}_0\gamma: \Fc \to \Rc$. For any $k\in\zz$, Eq. (\ref{a}) gives $\chih^{2n}_0\gamma(\Gc^{k+1}) \subset \chih^{2n}_0\circ\phi(G^{k-n+2}\Om^{2n-2}\Fc)$. But $\chih^{2n}_0$ is $\kappa$-invariant, hence $\chih^{2n}_0\circ \phi = \chih^{2n}_0\circ \phit$. Therefore, $\chih^{2n}_0\gamma(\Gc^{k+1}) \subset \chih^{2n}_0(G^{k-n+1}\Om^{2n}\Fc) \subset \Jc^{k-n+1}$ using $G^{k-n+1}\Om \subset H^{k-n+1}\Om$ and (\ref{1}). Now we look at its companion $\chih^{2n}_1\gamma: \Om^1\Fc_{\nat} \to \Om^1\Rc_{\nat}$. From (\ref{b}) one gets $\chih^{2n}_1(\nat\Gc^k\dd\Fc) \subset \chih^{2n}_1\circ\phi(G^{k-n+1}\Om^{2n-1}\Fc)$. But $\chih^{2n}_1$ is also $\kappa$-invariant and $\chih^{2n}_1\circ \phi = \chih^{2n}_1\circ \phit$, thus $\chih^{2n}_1\gamma(\nat\Gc^k\dd\Fc) \subset \chih^{2n}_1(G^{k-n}\Om^{2n+1}\Fc) \subset \nat \Jc^{k-n}\dd\Rc$ by (\ref{5}). This allows to estimate the image of $[\Gc^k,\Fc]=\bb\nat\Gc^k\dd\Fc$ under the chain map $\chih^{2n}\gamma$. Indeed $\chih^{2n}_0\gamma(\bb\nat\Gc^k\dd\Fc)=\bb \chih^{2n}_1\gamma(\nat\Gc^k\dd\Fc)\subset \bb \nat \Jc^{k-n}\dd\Rc$, so that $\chih^{2n}_0\gamma([\Gc^k,\Fc]) \subset  [\Jc^{k-n},\Rc]$. Collecting these results shows the effect of the map $\chih^{2n}\gamma$ on the adic filtration in even degree:
$$
\left\{ \begin{array}{llcl}
\chih^{2n}_0\gamma : & \Gc^{k+1}+[\Gc^k,\Fc] & {\longrightarrow} & \Jc^{k-n+1}+[\Jc^{k-n},\Rc] \\
\chih^{2n}_1\gamma : & \nat \Gc^k\dd\Fc & {\longrightarrow} & \nat \Jc^{k-n}\dd\Rc \end{array} \right.
$$
hence $\chih^{2n}\gamma : F^{2k}_{\Gc}X(\Fc)\to F^{2k-2n}_{\Jc}X(\Rc)$. To understand the effect on the filtration in odd degree, one has to evaluate $\chih^{2n}_1\gamma$ on $\nat\Gc^k\dd\Gc$. From (\ref{d}), one gets $\chih^{2n}_1\gamma(\nat\Gc^k\dd\Gc)\subset \chih^{2n}_1(H^{k-n+1}\Om^{2n+1}\Fc + b\phi(G^{k-n+1}\Om^{2n}\Fc))$. But (\ref{rchi}) shows $\chih^{2n}_1\circ b=0$,  and (\ref{2}) implies $\chih^{2n}_1\gamma(\nat\Gc^k\dd\Gc) \subset \nat(\Jc^{k-n+1}\dd\Rc + \Jc^{k-n}\dd\Jc)$. One thus gets the adic behaviour of the chain map $\chih^{2n}\gamma$ on the filtration of odd degree:
$$
\left\{\begin{array}{llcl}
\chih^{2n}_0\gamma : & \Gc^{k+1} & {\longrightarrow} & \Jc^{k-n+1} \\
\chih^{2n}_1\gamma : & \nat (\Gc^{k+1}\dd\Fc + \Gc^k\dd\Gc) & {\longrightarrow} & \nat (\Jc^{k-n+1}\dd\Rc + \Jc^{k-n}\dd\Jc) \end{array}\right. 
$$
hence $\chih^{2n}\gamma : F^{2k+1}_{\Gc}X(\Fc)\to F^{2k-2n+1}_{\Jc}X(\Rc)$ and $\chih^{2n}\gamma$ is a map of order $2n$. Using similar methods, one shows that $\chih^{2n+1}\gamma$ is of order $2n+1$. \\
Now we investigate the eta-cochain. Consider $\etah^{2n}_0\gamma: \Fc\to \Rc$. (\ref{a}) gives $\etah^{2n}_0\gamma(\Gc^{k+1}) \subset \etah^{2n}_0\phi(G^{k-n+2}\Om^{2n-2}\Fc)$. However $\etah^{2n}$ is not $\kappa$-invariant, so that we cannot replace $\phi$ by $\phit$. We are forced to consider $\phi(G^{k-n+2}\Om^{2n-2}\Fc)\subset H^{k-n+1}\Om^{2n}\Fc$ and consequently $\etah^{2n}_0\gamma(\Gc^{k+1}) \subset \Jc^{k-n+1}$ by (\ref{3}). Similarly, (\ref{c}) implies $\etah^{2n}_0\gamma([\Gc^k,\Fc]) \subset \etah^{2n}_0(H^{k-n+1}\Om^{2n}\Fc)+\etah^{2n}_0b\phi(G^{k-n+1}\Om^{2n-1}\Fc) \subset \etah^{2n}_0(H^{k-n}\Om^{2n}\Fc)$ hence $\etah^{2n}_0\gamma([\Gc^k,\Fc]) \subset \Jc^{k-n}$. Its companion $\etah^{2n}_1\gamma:\Om^1\Fc_{\nat}\to \Om^1\Rc_{\nat}$ evaluated on $\nat \Gc^k\dd\Fc$ uses equation (\ref{b}) again with $\phi(G^{k-n+1}\Om^{2n-1}\Fc) \subset H^{k-n}\Om^{2n+1}\Fc$, so that $\etah^{2n}_1\gamma (\nat \Gc^k\dd\Fc) \subset \etah^{2n}_1(H^{k-n}\Om^{2n+1}\Fc) \subset \nat(\Jc^{k-n}\dd\Rc+\Jc^{k-n-1}\dd\Jc)$ by (\ref{4}). This shows the effect of $\etah^{2n}\gamma$ on the filtration of even degree
$$
\left\{\begin{array}{llcl}
\etah^{2n}_0\gamma : & \Gc^{k+1}+[\Gc^k,\Fc] & {\longrightarrow} & \Jc^{k-n} \\
\etah^{2n}_1\gamma : & \nat \Gc^k\dd\Fc & {\longrightarrow} & \nat (\Jc^{k-n}\dd\Rc + \Jc^{k-n-1}\dd\Jc) \end{array} \right.
$$
hence $\etah^{2n}\gamma : F^{2k}_{\Gc}X(\Fc)\to F^{2k-2n-1}_{\Jc}X(\Rc)$. For the odd filtration, let us compute from (\ref{d}) $\etah^{2n}_1\gamma(\nat \Gc^k\dd\Gc) \subset \etah^{2n}_1(H^{k-n+1}\Om^{2n+1}\Fc)+\etah^{2n}_1b\phi(G^{k-n+1}\Om^{2n}\Fc)$. But the identities (\ref{reta}) show that $\etah^{2n}_1b= \chih^{2n+1}_1$, hence using $\kappa$-invariance one gets $\etah^{2n}_1b\phi(G^{k-n+1}\Om^{2n}\Fc) \subset \chih^{2n+1}_1\phit(G^{k-n+1}\Om^{2n}\Fc) \subset \chih^{2n+1}_1(G^{k-n}\Om^{2n+2}\Fc)$. Therefore, (\ref{4}) and (\ref{5}) imply $\etah^{2n}_1\gamma(\nat \Gc^k\dd\Gc) \subset \nat \Jc^{k-n}\dd\Rc$. These results give the adic behaviour of $\etah^{2n}\gamma$ with respect to the odd filtration
$$
\left\{\begin{array}{llcl}
\etah^{2n}_0\gamma : & \Gc^{k+1} & {\longrightarrow} & \Jc^{k-n+1} \\
\etah^{2n}_1\gamma : & \nat (\Gc^{k+1}\dd\Fc + \Gc^k\dd\Gc) & {\longrightarrow} & \nat \Jc^{k-n}\dd\Rc  \end{array} \right.
$$
hence $\etah^{2n}\gamma : F^{2k+1}_{\Gc}X(\Fc)\to F^{2k-2n}_{\Jc}X(\Rc)$ and $\etah^{2n}\gamma$ is a map of order $2n+1$. Similarly, one shows that $\etah^{2n+1}\gamma$ is of order $2n+2$.   \cqfd\\

Note that the chain maps $\gamma$ and $\chih^n$ extend to the adic completions of all the algebras involved, so that from now on we will consider the bivariant Chern character $\ch^n(\rho)\in \hom^n(\Xh(T\Ac,J\Ac),\Xh(\Rc,\Jc))$ as a chain map of pro-complexes
\be
\ch^n(\rho): X(\Th\Ac) \stackrel{\gamma}{\longrightarrow} \Omh \Th\Ac \stackrel{\rho_*}{\longrightarrow} \Omh\Mch^s_+ \stackrel{\chih^n}{\longrightarrow} X(\Rch)\ .
\ee
We would like to introduce some equivalence relations among quasihomomorphisms, and discuss the corresponding invariance properties of the Chern character. The first equivalence relation is (smooth) homotopy. It involves the algebra $\cinf[0,1]$ of smooth functions $f:[0,1]\to\cc$, such that all the derivatives of order $\geq 1$ vanish at the enpoints $0$ and $1$, while the values of $f$ itself remain arbitrary. We have already seen that $\cinf[0,1]$ endowed with its usual Fr\'echet topology is an $m$-algebra. It is moreover nuclear \cite{Gr}, so that its projective tensor product $\Ac\hotimes\cinf[0,1]$ with any $m$-algebra $\Ac$ is isomorphic to the algebra of smooth $\Ac$-valued functions over $[0,1]$, with all derivatives of order $\geq 1$ vanishing at the endpoints. We will usually denote by $\Ac[0,1]$ this $m$-algebra. The second equivalence relation of interest among quasihomomorphisms is conjugation by an invertible element of the unitalized algebra $(\Ec^s_+)^+$.

\begin{definition}
Let $\rho_0:\Ac\to\Ec^s\triangleright\Ic^s\hotimes\Bc$ and $\rho_1:\Ac\to\Ec^s\triangleright\Ic^s\hotimes\Bc$ be two quasihomomorphisms with same parity. They are called\\

\noindent {\bf i) homotopic} if there exists a quasihomomorphism $\rho:\Ac\to\Ec[0,1]^s\triangleright\Ic^s\hotimes\Bc[0,1]$ such that evaluation at the endpoints gives $\rho_0$ and $\rho_1$; \\

\noindent {\bf ii) conjugate} if there exists an invertible element in the unitalized algebra $U\in(\Ec^s_+)^+$ with $U-1\in \Ec^s_+$, such that $\rho_1=U^{-1}\rho_0 U$ as a homomorphism $\Ac\to \Ec^s_+$.
\end{definition}

Remark that the commutators $[F,U]$ and $[F,U^{-1}]$ always lie in the ideal $\Ic^s\hotimes\Bc\subset \Ec^s$. When the algebra $\Ic$ is $M_2$-stable (i.e. $M_2(\Ic)\cong\Ic$), two conjugate quasihomomorphisms are also homotopic, but the converse is not true. Hence conjugation is strictly stronger than homotopy as an equivalence relation. The former is an analogue of ``compact perturbation'' of quasihomomorphisms in Kasparov's bivariant $K$-theory for $C^*$-algebras, see \cite{Bl}. \\
The proposition below describes the compatibility between these equivalence relations and the bivariant Chern character.

\begin{proposition}\label{pinv}
Let $\rho_0:\Ac\to\Ec^s\triangleright\Ic^s\hotimes\Bc$ and $\rho_1:\Ac\to\Ec^s\triangleright\Ic^s\hotimes\Bc$ be two $(p+1)$-summable quasihomomorphisms of parity $p\mod 2$, with $\Ec$ admissible with respect to a quasi-free extension $\Rc$ of $\Bc$. Let $n\geq p$ be any integer of the same parity.\\

\noindent i) If $\rho_0$ and $\rho_1$ are homotopic, then $S\ch^n(\rho_0) \equiv S\ch^n(\rho_1)$ in $HC^{n+2}(\Ac,\Bc)$. In particular $\ch^n(\rho_0)\equiv \ch^n(\rho_1)$ whenever $n\geq p+2$.\\

\noindent ii) If $\rho_0$ and $\rho_1$ are conjugate, then $\ch^n(\rho_0) \equiv \ch^n(\rho_1)$ in $HC^{n}(\Ac,\Bc)$ for all $n\geq p$.
\end{proposition}
{\it Proof:} First observe that if $\rho:\Ac\to \Ec^s\triangleright \Ic^s\hotimes\Bc$ is a quasihomomorphism with $\Rc$-admissible algebra $\Ec$, the lifting homomorphism $\rho_*:T\Ac\to \Mc^s_+$ factors through the tensor algebra $T\Ec^s_+$ by virtue of the commutative diagram
$$
\xymatrix{
0 \ar[r]  & J\Ac \ar[r] \ar[d]_{\varphi} & T\Ac  \ar[r] \ar[d]_{\varphi} & \Ac \ar[r]  \ar[d]_{\rho} & 0  \\
0 \ar[r]  & J\Ec^s_+ \ar[r] \ar[d] & T\Ec^s_+  \ar[r] \ar[d] & \Ec^s_+ \ar[r]  \ar@{=}[d] & 0  \\
0 \ar[r] & \Nc^s_+ \ar[r] & \Mc^s_+ \ar[r] & \Ec^s_+ \ar[r] \ar@/_/@{.>}[l]_{\si} & 0 }
$$
where the homomorphism $\varphi: T\Ac \to T\Ec^s_+$ is $\varphi(a_1\otimes\ldots\otimes a_n) = \rho(a_1)\otimes\ldots\otimes \rho(a_n)$, and the arrow $T\Ec^s_+\to\Mc^s_+$ maps a tensor product $e_1\otimes\ldots\otimes e_n\in T\Ec^s_+$ to the product $\si(e_1)\ldots\si(e_n)$. By the naturality of the Goodwillie equivalences $\gamma_{\Ac}:X(T\Ac)\to \Omh T\Ac$ and $\gamma_{\Ec^s_+}:X(T\Ec^s_+)\to \Omh T\Ec^s_+$, one immediately sees that the bivariant Chern character coincides with the composition of chain maps 
$$
\ch^n(\rho): X(T\Ac) \stackrel{X(\varphi)}{\longrightarrow} X(T\Ec^s_+) \stackrel{\gamma_{\Ec^s_+}}{\longrightarrow} \Omh T\Ec^s_+ \stackrel{\chih^n}{\longrightarrow} X(\Rc)\ .
$$
Hence, all the information about the homomorphism $\rho:\Ac\to\Ec^s_+$ is concentrated in the chain map $X(\varphi):X(T\Ac)\to X(T\Ec^s_+)$. This will simplify the comparison of Chern characters associated to homotopic or conjugate quasihomomorphisms.\\
 
\noindent i) Homotopy: the cocycles $\ch^n(\rho_0)$ and $\ch^n(\rho_1)$ differ only by the chain maps $X(\varphi_{i}): X(T\Ac)\to X(T\Ec^s_+)$, $i=0,1$. We view $\rho:\Ac\to \Ec^s_+[0,1]$ as a smooth family of homomorphisms $\rho_t:\Ac\to \Ec^s_+$ parametrized by $t\in [0,1]$, giving a homotopy between the two endpoints $\rho_0$ and $\rho_1$. Cuntz and Quillen prove in \cite{CQ1} a Cartan homotopy formula which provides a transgression between the chain maps $X(\varphi_{i})$. At any point $t\in[0,1]$, denote by $\dot\varphi=\frac{d}{dt}\varphi : T\Ac\to T\Ec^s_+$ the derivative of the homomorphism $\varphi_t$ with respect to $t$, and define a linear map $\iota: \Om^mT\Ac\to \Om^{m-1}T\Ec^s_+$ by
$$
\iota(x_0\dd x_1\ldots \dd x_m) = (\varphi x_0)(\dot\varphi x_1)\dd(\varphi x_2)\ldots \dd(\varphi x_m)\ .
$$
The tensor algebra $T\Ac$ is quasi-free, hence consider any $\phi:T\Ac\to\Om^2 T\Ac$ verifying $\phi(xy)=\phi(x)y+x\phi(y)+\dd x\dd y$, and let $h:X(T\Ac)\to X(T\Ec^s_+)$ be the linear map of odd degree
$$
h(x)=\nat \iota \phi(x)\ ,\qquad h(\nat x\dd y) = \iota(x\dd y + b(x\phi(y)))
$$
(the latter is well-defined on $\nat x\dd y$). Then Cuntz and Quillen show the following adic properties of $h$ for any $k\in\zz$,
\beq
h(F^k_{J\Ac}X(T\Ac)) \subset F^{k-1}_{J\Ec^s_+}X(T\Ec^s_+) &\mbox{if}& \dot\varphi(J\Ac) \subset J\Ec^s_+\ ,\non\\
h(F^k_{J\Ac}X(T\Ac)) \subset F^k_{J\Ec^s_+}X(T\Ec^s_+) &\mbox{if}& \dot\varphi(T\Ac) \subset J\Ec^s_+\ ,\non
\eeq
and moreover the transgression formula $\frac{d}{dt}X(\varphi)  = [\d, h]$ holds. Hence if we define by integration over $[0,1]$ the odd chain $H=\int_0^1dt\, h$, one has
$$
X(\varphi_1)-X(\varphi_0)=[\d,H]
$$
in the complex $\hom^1(\Xh(T\Ac,J\Ac),\Xh(T\Ec^s_+,J\Ec^s_+))$ in case $\dot\varphi(J\Ac) \subset J\Ec^s_+$, or in the complex $\hom^0(\Xh(T\Ac,J\Ac),\Xh(T\Ec^s_+,J\Ec^s_+))$ in case $\dot\varphi(T\Ac) \subset J\Ec^s_+$. For a general homotopy we are in the first case $\dot\varphi(J\Ac) \subset J\Ec^s_+$. After composition by the chain map $\chih^n\gamma_{\Ec^s_+} \in \hom^n(\Xh(T\Ec^s_+,J\Ec^s_+),\Xh(\Rc,\Jc))$, this shows the transgression relation
$$
\ch^n(\rho_1)-\ch^n(\rho_0) = (-)^n[\d, \chih^n\gamma_{\Ec^s_+} H]\in \hom^{n+1}(\Xh(T\Ac,J\Ac),\Xh(\Rc,\Jc))\ ,
$$
whence $S\ch^n(\rho_1)\equiv S\ch^n(\rho_0)$ in $HC^{n+2}(\Ac,\Bc)$. The sign $(-)^n$ comes from the parity of the chain map $\chih^n\gamma_{\Ec^s_+}$. \\

\noindent ii) Conjugation: now $\varphi_0,\varphi_1:T\Ac\to T\Ec^s_+$ are the homomorphism lifts of $\rho_0$ and $\rho_1=U^{-1}\rho_0 U$. Introduce the pro-algebra $\Th\Ec^s_+=\varprojlim_k T\Ec^s_+/(J\Ec^s_+)^k\cong \prod_{k\geq 0}\Om^{2k}\Ec^s_+$, and consider the invertible $U\in(\Ec^s_+)^+$ as an element $\Uh$ of the unitalization $(\Th\Ec^s_+)^+$, via the linear inclusion of zero-forms $\Ec^s_+\hookrightarrow  \Th\Ec^s_+$. By proceeding as in \cite{CQ1}, it turns out that $\Uh$ is invertible, with inverse given by the series
$$
\Uh^{-1}= \sum_{k\geq 0} U^{-1}(dU\,dU^{-1})^k\ \in (\Th\Ec^s_+)^+\ .
$$
Of course the image of $\Uh^{-1}$ under the multiplication map $(\Th\Ec^s_+)^+\to (\Ec^s_+)^+$ is $U^{-1}$. We will show that $\varphi_1$, viewed as a homomorphism $T\Ac\to \Th\Ec^s_+$, is homotopic to the homomorphism $\Uh^{-1}\varphi_0\Uh$. For any $t\in [0,1]$ define a linear map $\si_t:\Ac\to \Th\Ec^s_+$ by
$$
\si_t(a) = (1-t)\rho_1(a) + t\, \Uh^{-1}\rho_0(a)\Uh\ ,\quad \forall a\in\Ac\ ,
$$
where $\rho_0(a)$ and $\rho_1(a)$ are considered as elements of the subspace of zero-forms $\Ec^s_+\hookrightarrow  \Th\Ec^s_+$. Thus $\si_t$ is a linear lifting of the constant homomorphism $\rho_1:\Ac\to\Ec^s_+$. Then use the universal property of the tensor algebra $T\Ac$ to build a smooth family of homomorphisms $\varphi(t): T\Ac\to \Th\Ec^s_+$ by means of the commutative diagram with exact rows
$$
\xymatrix{
0 \ar[r] & J\Ac \ar[r] \ar[d]_{\varphi(t)} & T\Ac \ar[r] \ar[d]_{\varphi(t)} & \Ac \ar[r] \ar[d]^{\rho_1} \ar@{.>}[dl]_{\si_t} & 0 \\
0 \ar[r] & \Jh\Ec^s_+ \ar[r] & \Th\Ec^s_+ \ar[r] & \Ec^s_+ \ar[r] & 0 }
$$
By construction one has $\varphi(0)=\varphi_1$, $\varphi(1)=\Uh^{-1}\varphi_0\Uh$ and the derivative $\dot\varphi$ sends $T\Ac$ to the ideal $\Jh\Ec^s_+$. Hence from the Cartan homotopy formula of part i) we deduce that the chain maps $X(\varphi_1)$ and $X(\Uh^{-1}\varphi_0\Uh)$ are cohomologous in the complex $\hom^0(\Xh(T\Ac,J\Ac),\Xh(T\Ec^s_+,J\Ec^s_+))$. Then we have to show that $ X(\Uh^{-1}\varphi_0\Uh)$ and $X(\varphi_0)$ are cohomologous. Consider the following linear map of odd degree $h: X(T\Ac) \to X(\Th\Ec^s_+)\cong \Xh(T\Ec^s_+,J\Ec^s_+) $ defined by
$$
h(x)=\nat(\Uh^{-1}\varphi_0(x) \dd \Uh)\ ,\qquad h(\nat x\dd y)=0\ .
$$
It is easy to see that $h$ defines a cochain of order zero, i.e. lies in the complex $\hom^0(\Xh(T\Ac,J\Ac),\Xh(T\Ec^s_+,J\Ec^s_+))$. Moreover, one has the transgression relation $X(\Uh^{-1}\varphi_0\Uh)-X(\varphi_0) = [\d, h]$. Indeed (we replace $\varphi_0 (x)$ by $x$ for notational simplicity)
\beq
X(\Uh^{-1}\varphi_0\Uh)(x)-X(\varphi_0)(x) &=& \Uh^{-1}x\Uh - x\ =\ [\Uh^{-1}x,\Uh]\non\\
&=& \bb\nat(\Uh^{-1}x \dd \Uh)\ =\ \bb h(x)\ ,\non
\eeq
\beq
\lefteqn{X(\Uh^{-1}\varphi_0\Uh)(\nat x\dd y)-X(\varphi_0)(\nat x\dd y) = \nat\Uh^{-1}x\Uh\dd(\Uh^{-1}y\Uh) - \nat x\dd y} \non\\
&&\qquad = \nat(\Uh^{-1}x\Uh \dd \Uh^{-1}y\Uh + \Uh^{-1}x\dd y\Uh + \Uh^{-1}xy\dd\Uh -x\dd y) \non\\
&& \qquad = \nat(-yx \dd\Uh \Uh^{-1} + xy \dd \Uh\Uh^{-1})\ =\nat \Uh^{-1}[x,y]\dd \Uh \non\\
&&\qquad =h(\bb\nat x\dd y)\ ,\non
\eeq
where in the second computation we use the identity $\dd \Uh^{-1}= -\Uh^{-1}\dd\Uh\Uh^{-1}$ deduced from $\dd 1=0$. This shows the equality of bivariant cyclic cohomology classes
$$
X(\varphi_1)\equiv X(\varphi_0)\ \in HC^0(\Ac,\Ec^s_+)\ ,
$$
so that after composition with $\chih^n\gamma_{\Ec^s_+} \in HC^n(\Ec^s_+,\Bc)$, the equality $\ch^n(\rho_1)\equiv \ch^n(\rho_0)$ holds in $HC^{n}(\Ac,\Bc)$. \cqfd\\

Part ii) of the above proof also shows the independence of the cohomology class $\ch^n(\rho)\in HC^n(\Ac,\Bc)$ with respect to the choice of linear splitting $\si:\Ec^s_+\to \Mc^s_+$ used to lift the homomorphism $\rho$, two such splittings being always homotopic. Then, from section \ref{scy} we know that any class in $HC^n(\Ac,\Bc)$ induces linear maps of degree $-n$ between the cyclic homologies of $\Ac$ and $\Bc$, compatible with the $SBI$ exact sequence. Hence, if the quasihomomorphism is $(p+1)$-summable and with parity $p\mod 2$, the lowest degree representative of the Chern character $\ch^p(\rho)\in HC^p(\Ac,\Bc)$ carries the maximal information. We collect these results in a theorem:

\begin{theorem}\label{tbiv}
Let $\rho:\Ac\to \Ec^s_+\triangleright \Ic^s_+\hotimes\Bc$ be a $(p+1)$-summable quasihomomorphism of parity $p \mod 2$, with $\Ec$ admissible with respect to a quasi-free extension of $\Bc$. The bivariant Chern character $\ch^p(\rho)\in HC^p(\Ac,\Bc)$ induces a graded-commutative diagram
$$
\vcenter{\xymatrix{
HP_{n+1}(\Ac) \ar[r]^S \ar[d]^{\ch^p(\rho)} & HC_{n-1}(\Ac) \ar[r]^B \ar[d]^{\ch^p(\rho)} & HN_n(\Ac)  \ar[r]^I \ar[d]^{\ch^p(\rho)} & HP_n(\Ac)  \ar[d]^{\ch^p(\rho)} \\
HP_{n-p+1}(\Bc) \ar[r]^S & HC_{n-p-1}(\Bc) \ar[r]^B & HN_{n-p}(\Bc) \ar[r]^I & HP_{n-p}(\Bc) }}
$$
invariant under conjugation of quasihomomorphisms. Moreover the arrow in periodic cyclic homology $HP_n(\Ac)\to HP_{n-p}(\Bc)$ is invariant under homotopy of quasihomomorphisms.
\end{theorem}
{\it Proof:} The fact that $S\ch^p(\rho)\in HC^{p+2}(\Ac,\Bc)$ is homotopy invariant shows its image in the periodic theory $HP^p(\Ac,\Bc)$ is homotopy invariant. \cqfd\\

\begin{example}\textup{When $\Ac$ is arbitrary and $\Bc=\cc$, we saw in Example \ref{ek} that a $(p+1)$-summable quasihomomorphism $\rho:\Ac\to \Lc^s\triangleright \Ic^s$, represents a $K$-homology class of $\Ac$. By hypothesis, the degree of the quasihomomorphism is $p\mod 2$. The Chern character $\ch^p(\rho)\in HC^p(\Ac,\cc)\cong HC^p(\Ac)$ is a cyclic cohomology class of degree $p$ over $\Ac$, represented by a chain map $\Xh(T\Ac,J\Ac)\to\cc$ vanishing on the subcomplex $F^p\Xh(T\Ac,J\Ac)$. Using the pro-vector space isomorphism $\Xh(T\Ac,J\Ac)\cong \Omh\Ac$, one finds that $\ch^p(\rho)$ is non-zero only on the subspace of $p$-forms $\Om^p\Ac$, explicitly
$$
\ch^p(\rho)(a_0da_1\ldots da_p)= \frac{c_p}{2}\, \Tr_s(F[F,a_0]\ldots [F,a_p])\ ,
$$
where $\Tr_s:(\Ic^s)^{p+1}\to \cc$ is the supertrace of the $(p+1)$-summable algebra $\Ic$ and $c_p$ is a constant depending on the degree. One has $c_p=(-)^n(n!)^2/p!$ when $p=2n$ is even, and $c_p=\sqrt{2\pi i}\, (-)^n/2^p$ when $p=2n+1$ is odd. This coincides with the Chern-Connes character \cite{C0, C1}, up to a scaling factor accounting for the homotopy equivalence between the $X$-complex $\Xh(T\Ac,J\Ac)$ and the $(b+B)$-complex $\Omh\Ac$.
}
\end{example}

\begin{example}\textup{When $\Ac=\cc$ and $\Bc=\cinf(0,1)$, the Bott element (Example \ref{ebott}) represented by the odd $1$-summable quasihomomorphism $\rho:\cc\to\Ec^s\triangleright \Ic^s\hotimes\Bc$ with $\Bc$-admissible extension $\Ec=\cinf[0,1]$, has a Chern character in $HC^1(\cc,\Bc)\cong HP_1(\Bc)$. The periodic cyclic homology of $\Bc$ is isomorphic to the de Rham cohomology of the open interval $(0,1)$, hence $HP_0(\Bc)=0$ and $HP_1(\Bc)=\cc$. Consequently, the Chern character $\ch^1(\rho)$ may be represented by a smooth one-form over $[0,1]$ vanishing at the endpoints. It involves a real-valued function $\xi\in \Ec$, with $\xi(0)=0$ and $\xi(1)=\pi/2$, used in the construction of the homomorphism $\rho:\cc\to\Ec^s_+$. One explicitly finds
$$
\ch(\rho)= \sqrt{2\pi i}\, \dd(\sin^2\xi)\ ,
$$
so that its integral over the interval $[0,1]$ is normalized to $\sqrt{2\pi i}$, and of course does no depend on the chosen function $\xi$. This is due to the fact that quasihomomorphisms associated to different choices of $\xi$ are homotopic. }
\end{example}

\section{Topological $K$-theory}\label{st}

We review here the topological $K$-theory of Fr\'echet $m$-algebras following Phillips \cite{Ph}, and construct various Chern character maps with value in cyclic homology. Topological $K$-theory for Fr\'echet $m$-algebras is defined in analogy with Banach algebras and fulfills the same properties of homotopy invariance, Bott periodicity and excision \cite{Ph}. For our purposes, only homotopy invariance and Bott periodicity are needed. A basic example of Fr\'echet $m$-algebra is provided by the algebra $\Kc$ of ``smooth compact operators''. $\Kc$ is the space of infinite matrices $(A_{ij})_{i,j\in\nn}$ with entries in $\cc$ and rapid decay, endowed with the family of submultiplicative norms
$$
\|A\|_n= \sup_{(i,j)\in\nn^2} (1+i+j)^nA_{ij} < \infty\qquad \forall n\in\nn\ .
$$
The multiplication of matrices makes $\Kc$ a Fr\'echet $m$-algebra. Moreover $\Kc$ is nuclear as a locally convex vector space \cite{Gr}. If $\Ac$ is any Fr\'echet $m$-algebra, the completed tensor product $\Kc\hotimes\Ac$ is the \emph{smooth stabilization} of $\Ac$. Other important examples are the algebras $\cinf[0,1]$, resp. $\cinf(0,1)$, of smooth $\cc$-valued functions over the interval, with all derivatives of order $\geq 1$, resp. $\geq 0$, vanishing at the endpoints. As already mentioned in section \ref{sbiv}, these are again nuclear Fr\'echet $m$-algebras and the completed tensor products $\Ac[0,1]=\Ac\hotimes \cinf[0,1]$ and $\Ac(0,1)=\Ac\hotimes \cinf(0,1)$ are isomorphic to the algebras of smooth $\Ac$-valued functions over the interval, with the appropriate vanishing boundary conditions. In particular $S\Ac:=\Ac(0,1)$ is the \emph{smooth suspension} of $\Ac$. We say that two idempotents $e_0$, $e_1$ of an algebra $\Ac$ are smoothly homotopic if there exists an idempotent $e\in \Ac[0,1]$ whose evaluation at the endpoints gives $e_0$ and $e_1$. Similarly for invertible elements.\\ 
The definition of topological $K$-theory involves idempotents and invertibles of the unitalized algebra $(\Kc\hotimes\Ac)^+$. Choosing an isomorphism $M_2(\Kc)\cong\Kc$ makes $(\Kc\hotimes\Ac)^+$ a semigroup for the direct sum $a\oplus b=\bigl(
\begin{smallmatrix}
a & 0 \\
0 & b \end{smallmatrix} \bigr)$. We denote by $p_0$ the idempotent $\bigl(
\begin{smallmatrix}
1 & 0 \\
0 & 0 \end{smallmatrix} \bigr)$ of the matrix algebra $M_2(\Kc\hotimes\Ac)^+$.

\begin{definition}[Phillips \cite{Ph}]
Let $\Ac$ be a Fr\'echet $m$-algebra. The topological $K$-theory of $\Ac$ in degree zero and one is defined by
\beq
\Kt_0(\Ac) &=& \{\mbox{set of smooth homotopy classes of idempotents $e\in M_2(\Kc\hotimes\Ac)^+$}\non\\
 && \qquad \mbox{such that}\ e- p_0 \in M_2(\Kc\hotimes\Ac)\ \}\non\\
\Kt_1(\Ac) &=& \{\mbox{set of smooth homotopy classes of invertibles $g\in (\Kc\hotimes\Ac)^+$}\non\\
 && \qquad  \mbox{such that}\ g-1 \in \Kc\hotimes\Ac\ \}\non
\eeq
\end{definition}

$\Kt_0(\Ac)$ and $\Kt_1(\Ac)$ are semigroups for the direct sum of idempotents and invertibles; in the case of idempotents, the direct sum $e\oplus e' \in M_4(\Kc\hotimes\Ac)^+$ has to be conjugated by the invertible matrix
\be
c=\left( \begin{matrix}
1 & 0 & 0 & 0 \\
0 & 0 & 1 & 0 \\
0 & 1 & 0 & 0 \\
0 & 0 & 0 & 1 \end{matrix} \right)\ \in M_4(\Kc\hotimes\Ac)^+\ ,\quad c^{-1}=c\ ,\label{mat}
\ee
in order to preserve the condition $c(e\oplus e')c - \pt_0\in M_4(\Kc\hotimes\Ac)$, with $\pt_0$ the diagonal matrix diag$(1,1,0,0)$. The proof that $\Kt_0(\Ac)$ and $\Kt_1(\Ac)$ are actually abelian groups will be recalled in Lemma \ref{lab}. The unit of $\Kt_0(\Ac)$ is the class of the idempotent $p_0\in M_2(\Kc\hotimes\Ac)^+$, whereas the unit of $\Kt_1(\Ac)$ is represented by $1\in (\Kc\hotimes\Ac)^+$. \\
The fundamental property of topological $K$-theory is Bott periodicity \cite{Ph}. Let $S\Ac=\Ac(0,1)$ be the smooth suspension of $\Ac$. Define two additive maps
\be
\al:\Kt_1(\Ac)\to \Kt_0(S\Ac)\ ,\qquad \beta:\Kt_0(\Ac)\to \Kt_1(S\Ac)
\ee
as follows. First choose a real-valued function $\xi\in \cinf[0,1]$ such that $\xi(0)=0$ and $\xi(1)=\pi/2$ (we recall that all the derivatives of $\xi$ vanish at the endpoints). Let $g\in (\Kc\hotimes\Ac)^+$ represent an element of $\Kt_1(\Ac)$. Then the idempotent
$$
\al(g)=G^{-1}p_0G\ ,\qquad \al(g)-p_0\in M_2(\Kc\hotimes S\Ac)
$$
defines an element of $\Kt_0({S}\Ac)$, where $G:[0,1]\to M_2(\Kc\hotimes\Ac)^+$ is the matrix function over $[0,1]$
$$
G=R^{-1}\left(
\begin{array}{cc}
1 & 0 \\
0 & g \end{array} \right)R\quad \mbox{with}\quad R=\left(
\begin{array}{cc}
\cos\xi & \sin\xi \\
-\sin\xi & \cos\xi \end{array} \right)\ . 
$$
Now $z=\exp(4i\xi)$ is a complex-valued invertible function over $[0,1]$ with winding number 1. The functions $z-1$ and $z^{-1}-1$ lie in $\cinf(0,1)$. Then for any idempotent $e\in M_2(\Kc\hotimes\Ac)^+$ representing a class in $\Kt_0(\Ac)$, we define the invertible element
$$
\beta(e)=(1+(z-1)e)(1+(z-1)p_0)^{-1}\ .
$$
One has $(1+(z-1)p_0)^{-1}=(1+(z^{-1}-1)p_0)$, and the idempotent relations $e^2=e$, $p_0^2=p_0$ imply $\beta(e)=1+(z-1)e(e-p_0)+(z^{-1}-1)(p_0-e)p_0$, which shows that $\beta(e)-1$ is an element of the algebra $M_2(\Kc\hotimes{S}\Ac)\cong \Kc\hotimes{S}\Ac$, hence $\beta(e)$ defines a class in $\Kt_1({S}\Ac)$.

\begin{proposition}[Bott periodicity \cite{Ph}]
The two maps defined above $\al:\Kt_1(\Ac)\to \Kt_0({S}\Ac)$ and $\beta:\Kt_0(\Ac)\to \Kt_1({S}\Ac)$ are isomophisms of abelian groups. \cqfd
\end{proposition}

Hence Bott periodicity implies $\Kt_i({S}^2\Ac)=\Kt_i(\Ac)$ for $i=0,1$, so that we may define topological $K$-theory groups in any degree $n\in\zz$:
\be
\Kt_n(\Ac) = \left\{\begin{array}{l} 
\Kt_0(\Ac)\quad n\ \mbox{even}\\
\Kt_1(\Ac)\quad n\ \mbox{odd}\ .
\end{array} \right.
\ee
Following Cuntz and Quillen \cite{CQ1}, we construct Chern characters with values in periodic cyclic homology $\Kt_n(\Ac)\to HP_n(\Ac)$. Recall (section \ref{scy}) that periodic cyclic homology is computed from any quasi-free extension $0 \to \Jc \to \Rc \to \Ac \to 0$ by the pro-complex
$$
\Xh(\Rc,\Jc)= X(\Rch)\ :\ \Rch \rightleftarrows \Om^1\Rch_{\nat}\ ,
$$
where the pro-algebra $\Rch=\varprojlim_n \Rc/\Jc^n$ is the $\Jc$-adic completion of the quasi-free algebra $\Rc$. In particular, the universal free extension $0 \to J\Ac \to T\Ac \to \Ac \to 0$ is quasi-free and the universal property of the tensor algebra leads to a classifying homomorphism $T\Ac\to \Rc$ compatible with the ideals $J\Ac$ and $\Jc$ by means of the commutative diagram
$$
\vcenter{\xymatrix{
0 \ar[r]  & J\Ac \ar[r] \ar[d] & T\Ac  \ar[r] \ar[d] & \Ac \ar[r]  \ar@{=}[d] & 0  \\
0 \ar[r] & \Jc \ar[r] & \Rc \ar[r] & \Ac \ar[r] \ar@/_/@{.>}[l]_{\si} & 0 }}
$$
for any choice of continuous linear section $\si:\Ac\to\Rc$. The homomorphism $T\Ac\to\Rc$ thus extends to a homomorphism of pro-algebras $\Th\Ac\to \Rch$ and the induced morphism of complexes $X(\Th\Ac)\to X(\Rch)$ is a homotopy equivalence. The Chern character on topological $K$-theory requires to lift idempotents and invertible elements from the algebra $\Kc\hotimes\Ac$ to the pro-algebra 
$$
\Kc\hotimes \Rch = \varprojlim_n \Kc\hotimes (\Rc/\Jc^n)\ .
$$ 
If $e\in M_2(\Kc\hotimes\Ac)^+$ is an idempotent such that $e-p_0\in M_2(\Kc\hotimes\Ac)$, there always exists an idempotent lift $\eh\in M_2(\Kc\hotimes\Rch)^+$ with $\eh-p_0\in M_2(\Kc\hotimes\Rch)$, and two such liftings are always conjugate \cite{CQ1}. A concrete way to construct an idempotent lift is to work first with the tensor algebra and then push forward by the homomorphism $\Kc\hotimes\Th\Ac\to \Kc\hotimes\Rch$. Using the isomorphism of pro-vector spaces $\Th\Ac\cong \Omh^+\Ac$, the following differential form of even degree defines an idempotent \cite{CQ1}
\be
\eh= e + \sum_{k\geq 1} \frac{(2k)!}{(k!)^2}(e-\frac{1}{2})(dede)^k\ \in M_2(\Kc\hotimes \Th\Ac)^+ \ ,\label{elif}
\ee
where concatenation products over $M_2(\Kc)$ are taken. We will refer to (\ref{elif}) as the \emph{canonical lift} of $e$, but it should be stressed that other choices are possible. Denoting also by $\eh$ its image in $M_2(\Kc\hotimes\Rch)^+$, the Chern character of $e$ is represented by the cycle of even degree
\be
\ch_0(\eh)=\Tr(\eh-p_0)\ \in \Rch \ ,\label{ch0}
\ee
where the partial trace $\Tr:M_2(\Kc\hotimes\Rch)\to\Rch$ comes from the usual trace of matrices with rapid decay. We will show below that the cyclic homology class of $\ch_0(\eh)$ is invariant under smooth homotopies of $\eh$. Moreover, the invariance of the trace under similarity implies that $\ch_0$ is additive. Next, if $g\in (\Kc\hotimes\Ac)^+$ is an invertible element such that $g-1 \in \Kc\hotimes\Ac$, we have again to choose an invertible lift $\gh\in (\Kc\hotimes\Rch)^+$ with $\gh-1\in \Kc\hotimes\Rch$. It turns out that any lifting of $g$ is invertible, and two such liftings are always homotopic \cite{CQ1}. A concrete way to construct an invertible lift is to use the linear inclusion of zero-forms $\Kc\hotimes\Ac\hookrightarrow \Kc\hotimes\Th\Ac\cong \Kc\hotimes\Omh^+\Ac$ and consider $g$ as an element $\gh=g\in (\Kc\hotimes\Th\Ac)^+$. A simple computation shows that it is invertible, with inverse
\be
\gh^{-1}=\sum_{k\geq 0} g^{-1}(dg\,dg^{-1})^k\ \in (\Kc\hotimes\Th\Ac)^+\ . \label{glif}
\ee
Here again we shall refer to the above $\gh$ as the \emph{canonical lift} of $g$, but other choices are possible. Then denoting also by $\gh$ its image in $(\Kc\hotimes\Rch)^+$, the Chern character of $g$ is represented by the cycle of odd degree
\be
\ch_1(\gh)=\frac{1}{\sqrt{2\pi i}}\, \Tr \nat \gh^{-1}\dd \gh\ \in \Om^1\Rch_{\nat} \ ,\label{ch1}
\ee
with the trace map $\Tr:\Om^1(\Kc\hotimes\Rch)_{\nat}\to\Om^1\Rch_{\nat}$. In this case also we will show that the cyclic homology class of $\ch_1(\gh)$ is invariant under smooth homotopies of $\gh$. Clearly $\ch_1$ is additive. The factor $1/\sqrt{2\pi i}$ is chosen for consistency with the bivariant Chern character. \\
Note the following important property of idempotents and invertibles: two idempotents $\eh_0,\eh_1 \in M_2(\Kc\hotimes\Rch)^+$ are homotopic if and only if their projections $e_0,e_1\in M_2(\Kc\hotimes\Ac)^+$ are homotopic, and similarly with invertibles \cite{CQ1}. Since the cyclic homology classes of the Chern characters $\ch_0(\eh)$ and $\ch_1(\gh)$ are homotopy invariant with respect to $\eh$ and $\gh$, one gets well-defined additive maps $\ch_0: \Kt_0(\Ac)\to HP_0(\Ac)$ and $\ch_1: \Kt_1(\Ac)\to HP_1(\Ac)$ on the topological $K$-theory groups. They do not depend on the quasi-free extension $\Rc$ since we know that the classifying homomorphism $\Th\Ac\to \Rch$ induces a homotopy equivalence of pro-complexes $X(\Th\Ac)\stackrel{\sim}{\to} X(\Rch)$.\\ 
To show the homotopy invariance of the Chern characters, we introduce the Cherns-Simons transgressions. Let $\Rch[0,1]$ be the tensor product $\Rch\hotimes \cinf[0,1]$, and let $\eh$ be any idempotent of $M_2(\Kc\hotimes\Rch[0,1])^+$ with $\eh-p_0\in M_2(\Kc\hotimes\Rch[0,1])$. Denote by $s:\cinf[0,1]\to \Om^1[0,1]$ the de Rham coboundary map with values in ordinary (commutative) one-forms over the interval. We then define the Chern-Simons form associated to $\eh$ as the chain of odd degree
\be
\cs_1(\eh)=\int_0^1\Tr\nat (-2\eh+1)\,s\eh\,\dd\eh\ \in \Om^1\Rch_{\nat}\ ,\label{cs1}
\ee
with obvious notations. Now let $\gh\in (\Kc\hotimes\Rch[0,1])^+$ be any invertible element such that $\gh-1\in \Kc\hotimes\Rch[0,1]$. The Chern-Simons form associated to $\gh$ is the chain of even degree
\be
\cs_0(\gh)=\frac{1}{\sqrt{2\pi i}}\int_0^1\Tr(\gh^{-1}s\gh)\ \in \Rch\ .\label{cs0}
\ee

\begin{lemma}\label{ltrans}
Let $\eh$ be an idempotent of the algebra $M_2(\Kc\hotimes\Rch[0,1])^+$ with $\eh-p_0\in M_2(\Kc\hotimes\Rch[0,1])$. Denote by $\eh_0$ and $\eh_1$ the idempotents of $M_2(\Kc\hotimes\Rch)^+$ obtained by evaluation at $0$ and $1$. Then one has
\be
\bb\cs_1(\eh)=\ch_0(\eh_1)-\ch_0(\eh_0)\ \in \Rch\ .
\ee
Let $\gh\in (\Kc\hotimes\Rch[0,1])^+$ be an invertible element such that $\gh-1\in \Kc\hotimes\Rch[0,1]$. Denote by $\gh_0$ and $\gh_1$ the invertibles of $(\Kc\hotimes\Rch)^+$ obtained by evaluation at $0$ and $1$. Then one has
\be
\nat\dd\cs_0(\gh)=\ch_1(\gh_1)-\ch_1(\gh_0)\ \in \Om^1\Rch_{\nat}\ .
\ee
\end{lemma}
{\it Proof:} First notice that the current $\int_0^1$ is odd, so that 
$$
\bb \cs_1(\eh)= - \int_0^1 \bb\Tr\nat (-2\eh+1)\,s\eh\,\dd\eh\ ,
$$
and taking into account the fact that $s\eh$ is also odd, one has
$$
\bb \Tr\nat (-2\eh+1)\,s\eh\,\dd\eh=-\Tr [ (-2\eh+1)\,s\eh ,\eh]= \Tr((2\eh-1)s\eh\, \eh - \eh s\eh)\ ,
$$
where we use the idempotent property $\eh^2=\eh$ for the last equality. Since $s$ is a derivation, one has $s\eh= s(\eh^2)=s\eh\, \eh+\eh s\eh$ and $\eh s\eh\, \eh=0$, whence
$$
\bb\cs_1(e)=\int_0^1 \Tr (s\eh\, \eh + \eh s\eh)=\int_0^1 s \Tr \eh = \Tr(\eh_1-\eh_0)= \ch_0(\eh_1)-\ch_0(\eh_0)\ ,
$$
because $\eh_0$ and $\eh_1$ are the evaluations of $\eh$ respectively at $0$ and $1$. Let us proceed now with invertibles:
$$
\nat\dd \cs_0(\gh)=\frac{-1}{\sqrt{2\pi i}}\int_0^1 \Tr\nat \dd(\gh^{-1}s\gh)\ ,
$$
and because $\dd$ is an odd derivation anticommuting with $s$, one has
$$
\Tr\nat\dd(\gh^{-1}s\gh)= \Tr\nat(\dd \gh^{-1} s\gh - \gh^{-1}s\dd\gh)\ .
$$
Then, $\dd 1=0=s1$ implies $\dd \gh^{-1}=-\gh^{-1}\dd\gh \gh^{-1}$ and $s\gh^{-1}=-\gh^{-1}s\gh \gh^{-1}$. But $\Tr\nat$ is a supertrace, hence
$$
\Tr\nat(-\gh^{-1}\dd\gh \gh^{-1}s\gh-\gh^{-1}s\dd\gh)= \Tr\nat(\gh^{-1}s\gh \gh^{-1}\dd\gh-\gh^{-1}s\dd\gh)=-s \Tr\nat(\gh^{-1}\dd\gh)\ .
$$
By integration over the interval $[0,1]$, one gets
$$
\nat\dd \cs_0(\gh)=\frac{1}{\sqrt{2\pi i}}\Tr\nat(\gh_1^{-1}\dd\gh_1-\gh_0^{-1}\dd\gh_0)=\ch_1(\gh_1)-\ch_1(\gh_0)
$$
as wanted. \cqfd\\

Hence the cyclic homology classes of the Chern characters are homotopy invariant as claimed. There is another consequence of the above lemma. Observe that the suspensions ${S}\Ac=\Ac(0,1)$ and $S\Rch=\Rch(0,1)$ are subalgebras of $\Ac[0,1]$ and $\Rch[0,1]$. If $e\in M_2(\Kc\hotimes{S}\Ac)^+$ is an idempotent representing a class in $\Kt_0({S}\Ac)$, and $g\in(\Kc\hotimes{S}\Ac)^+$ an invertible representing a class in $\Kt_1({S}\Ac)$, we choose some lifts $\eh\in M_2(\Kc\hotimes{S}\Rch)^+$ and $\gh\in (\Kc\hotimes{S}\Rch)^+$. Then $\cs_1(\eh)$ and $\cs_0(\gh)$ are closed and define homology classes in $HP_1(\Ac)$ and $HP_0(\Ac)$ respectively. The following lemma shows the compatibility with Bott periodicity:

\begin{lemma}\label{lbott1}
The Chern-Simons forms define additive maps $\cs_1:\Kt_0({S}\Ac)\to HP_1(\Ac)$ and $\cs_0:\Kt_1({S}\Ac)\to HP_0(\Ac)$. Moreover they are compatible with the Bott isomorphisms $\al:\Kt_1(\Ac)\to\Kt_0({S}\Ac)$ and $\beta:\Kt_0(\Ac)\to\Kt_1({S}\Ac)$ and Chern characters, up to multiplication by a factor $\sqrt{2\pi i}$:
\beq
\cs_1\circ\al \equiv \sqrt{2\pi i}\, \ch_1 &:&\Kt_1(\Ac)\to HP_1(\Ac)\ , \non\\
\cs_0\circ\beta  \equiv \sqrt{2\pi i}\, \ch_0 &:& \Kt_0(\Ac)\to HP_0(\Ac)\ . \non
\eeq
\end{lemma}
{\it Proof:} Let $\eh$ be any idempotent of the pro-algebra $M_2(\Kc\hotimes{S}\Rch)^+$. We have to prove the homotopy invariance of the cyclic homology class determined by the cycle
$$
\cs_1(\eh)=\int_0^1\Tr\nat (-2\eh+1)\,s\eh\,\dd\eh\ .
$$
To this end, consider a smooth family of idempotents $\eh_t\in M_2(\Kc\hotimes{S}\Rch)^+$ parametrized by $t\in \rr$, such that $\eh_t-p_0\in M_2(\Kc\hotimes{S}\Rch), \forall t$. Denote by $\dot{\eh}$ the derivative $\d \eh/\d t$. The idempotent property of the family $\eh$ implies the following identity:
$$
\frac{\d}{\d t}\Tr\nat(-2\eh+1)s\eh\dd\eh=-\nat\dd\Tr(\eh(\dot{\eh} s\eh-s\eh\dot{\eh}))- s\Tr\nat\eh(\dot{\eh} \dd\eh-\dd\eh\dot{\eh})\ .
$$
Since for any fixed $t$, the idempotent $\eh_t$ equals $p_0$ at the boundaries of the suspended algebra $M_2(\Kc\hotimes{S}\Rch)^+$, one gets $\int_0^1s \Tr\nat\eh(\dot{\eh} \dd\eh-\dd\eh\dot{\eh})=0$ and
$$
\frac{\d}{\d t}\int_0^1\Tr\nat (-2\eh+1)\,s\eh\,\dd\eh= \nat\dd\int_0^1 \Tr(\eh(\dot{\eh} s\eh-s\eh\dot{\eh}))\ .
$$
This implies that the cyclic homology class of $\cs_1(\eh)$ is a homotopy invariant of $\eh$. Hence if $\eh$ lifts an idempotent $e\in M_2(\Kc\hotimes{S}\Ac)^+$, the cyclic homology class of $\cs_1(\eh)$ only depends on the homotopy class of $e$, and  the map $\cs_1:\Kt_0({S}\Ac)\to HP_1(\Ac)$ is well-defined. We now have to show the compatibility with Bott periodicity. Thus let $g\in (\Kc\hotimes\Ac)^+$ be an invertible such that $g-1\in \Kc\hotimes\Ac$, and let $\al(g)$ be the idempotent $G^{-1}p_0 G\in M_2(\Kc\hotimes{S}\Ac)^+$ constructed by means of a rotation matrix
$$
R=\left(
\begin{matrix}
\cos\xi & \sin\xi \\
-\sin\xi & \cos\xi \end{matrix} \right)\ ,\qquad G=R^{-1}\left(
\begin{matrix}
1 & 0 \\
0 & g \end{matrix} \right)R\ ,
$$
where $\xi\in \cinf[0,1]$ is a real function with $\xi(0)=0$ and $\xi(1)=\pi/2$. Then, it is clear that the idempotent $\eh=\Gh^{-1}p_0 \Gh \in M_2(\Kc\hotimes{S}\Rch)^+$ is a lifting of $\al(g)$, where the matrix $\Gh=R^{-1}\bigl(
\begin{smallmatrix}
1 & 0 \\
0 & \gh \end{smallmatrix} \bigr)R$ is built from any lifting $\gh\in (\Kc\hotimes\Rch)^+$ of $g$. Hence, the cyclic homology class of $\cs_1(\widehat{\al(g)})$ is represented by $\cs_1(\eh)$. A direct computation shows the equality
$$
\Tr\nat (-2\eh+1)\,s\eh\,\dd\eh= s(\cos\xi)\Tr\nat(-\gh^{-1}\dd\gh + \frac{1}{2}\dd\gh - \frac{1}{2}\dd\gh^{-1})\ ,
$$
so that after integration over $[0,1]$ one gets, modulo boundaries $\nat\dd(\cdot)$
$$
\cs_1(\widehat{\al(g)})\equiv\Tr\nat(\gh^{-1}\dd\gh)\mod\nat\dd \equiv \sqrt{2\pi i}\, \ch_1(\gh)\mod\nat\dd\ .
$$
Next, we turn to the map $\cs_0$. If $\gh\in (\Kc\hotimes{S}\Rch)^+$ is any invertible, one has
$$
\cs_0(\gh)=\frac{1}{\sqrt{2\pi i}}\int_0^1\Tr(\gh^{-1}s\gh)\ .
$$
We have to show that the cyclic homology class of $\cs_0(\gh)$ is a homotopy invariant of $\gh$. To this end, consider a smooth one-parameter family of invertibles $\gh_t\in (\Kc\hotimes{S}\Rch)^+$. One has, with $\dot{\gh}=\d \gh/\d t$,
$$
\frac{\d}{\d t}(\gh^{-1}s\gh)=-\gh^{-1}\dot{\gh}\gh^{-1}s\gh + \gh^{-1}s\dot{\gh} = [\gh^{-1}s\gh, \gh^{-1}\dot{\gh}] + s(\gh^{-1}\dot{\gh})\ .
$$
Since $\Tr [\gh^{-1}s\gh, \gh^{-1}\dot{\gh}]=- \bb\Tr\nat \gh^{-1}s\gh\dd(\gh^{-1}\dot{\gh})$, we get
$$
\frac{\d}{\d t}\int_0^1\Tr(\gh^{-1}s\gh)=\bb\int_0^1 \Tr\nat \gh^{-1}s\gh\dd(\gh^{-1}\dot{\gh})\ .
$$
Hence the cyclic homology class of $\cs_0(\gh)$ is homotopy invariant. In particular if $\gh$ lifts an invertible $g\in (\Kc\hotimes{S}\Ac)^+$, the cyclic homology class of $\cs_0(\gh)$ is a homotopy invariant of $g$ and the map $\cs_0:\Kt_1({S}\Ac)\to HP_0(\Ac)$ is well-defined. Now let $e\in M_2(\Kc\hotimes\Ac)^+$ be an idempotent, with $e-p_0\in M_2(\Kc\hotimes\Ac)$. Its image under the Bott map $\beta$ is the invertible element  $\beta(e)\in (\Kc\hotimes{S}\Ac)^+$ given by
$$
\beta(e)=(1+(z-1)e)(1+(z-1)p_0)^{-1}\ ,
$$
where $z=\exp(4i\xi)$. If $\eh\in M_2(\Kc\hotimes\Rch)^+$ is any idempotent lift of $e$, it is clear that the invertible 
$$
\gh=(1+(z-1)\eh)(1+(z-1)p_0)^{-1}\ \in (\Kc\hotimes{S}\Rch)^+
$$
is a lifting of $\beta(e)$. Hence the cyclic homology class of $\cs_0(\widehat{\beta(e)})$ is represented by $\cs_0(\gh)$. Let us compute explicitly $\Tr(\gh^{-1}s\gh)$:
\beq
\lefteqn{\Tr\big( (1+(z-1)p_0)(1+(z-1)\eh)^{-1}s\big((1+(z-1)\eh)(1+(z-1)p_0)^{-1}\big)\big)}\non\\
&&\qquad\qquad =\Tr\big( (1+(z^{-1}-1)\eh)sz\,\eh - (1+(z^{-1}-1)p_0)sz\, p_0 \big)\non\\
&& \qquad\qquad =\Tr(\eh-p_0)\, z^{-1}sz\ .\non
\eeq
Since the integration of $z^{-1}sz$ over the interval $[0,1]$ yields a factor $2\pi i$, one is left with equivalences modulo boundaries $\bb(\cdot)$
$$
\cs_0(\widehat{\beta(e)})\equiv \sqrt{2\pi i}\, \Tr(\eh-p_0)\mod\bb \equiv \sqrt{2\pi i}\, \ch_0(\eh)\mod\bb
$$
as wanted.  \cqfd\\

Since our main motivation is index theory we will have to consider the stabilization of $\Ac$ by a $p$-summable Fr\'echet $m$-algebra $\Ic$, that is, $\Ic$ is provided with a continuous trace $\Tr:\Ic^p\to\cc$ as in section \ref{sbiv}. Hence it will be convenient to define a Chern character $\Kt_n(\Ic\hotimes\Ac)\to HP_n(\Ac)$. The difficulty of course is that the trace is not defined on the algebra $\Kc\hotimes\Ic$ but only on its $p$-th power. To cope with this problem, we shall construct higher analogues of the Chern characters and Chern-Simons forms associated to idempotents and invertibles. Consider the following $p$-summable quasihomomorphism of even degree, from the algebra $\Ic\Ac := \Ic\hotimes\Ac$ to $\Ac$:
\be
\rho:\Ic\Ac\to \Ec^s\triangleright \Ic^s\Ac\ ,\qquad \Ec=\Ic^+\hotimes\Ac\ ,
\ee
where $\Ic^+$ is the unitalization of $\Ic$. Because $\rho$ is of even degree, it is entirely specified by a pair of homomorphisms $(\rho_+,\rho_-): \Ic\Ac\rightrightarrows\Ec$ which agree modulo the ideal $\Ic\Ac\subset\Ec$. Equivalently if we represent $\Ec^s$ in the $\zz_2$-graded matrix algebra $M_2(\Ec)$ we can write $\rho =\bigl( \begin{smallmatrix} \rho_+ & 0 \\ 0 & \rho_- \end{smallmatrix} \bigr)$. By definition we set
$$
\rho_+=\Id:\Ic\Ac\to\Ic\Ac\subset\Ec\ ,\qquad \rho_-=0\ .
$$
Let $0\to \Jc \to \Rc \to \Ac \to 0$ be any quasi-free extension of $\Ac$, with continuous linear splitting $\si:\Ac\to\Rc$. Then choosing $\Mc=\Ic^+\hotimes\Rc$ and $\Nc=\Ic^+\hotimes\Jc$, one gets a commutative diagram
$$
\vcenter{\xymatrix{
0\ar[r] & \Nc \ar[r] & \Mc \ar[r] & \Ec \ar[r] & 0 \\
0 \ar[r] & \Ic \Jc \ar[r] \ar[u] & \Ic\Rc \ar[r] \ar[u] & \Ic\Ac \ar[r] \ar[u] & 0}}
$$
Since $(\Ic^+)^n=\Ic^+$ for any integer $n$, one sees that $\Nc^n=\Ic^+\hotimes\Jc^n$ is a direct summand in $\Mc$. Moreover, there is an obvious chain map $\Tr: F^{2n+1}_{\Ic\Rc}X(\Mc)\to X(\Rc)$ for any $n\geq p-1$, obtained by taking the trace $\Ic^{n+1}\to \cc$. It follows that the algebra $\Ec\triangleright \Ic\Ac$ is $\Rc$-admissible (Definition \ref{dadm}), hence in all degrees $2n+1\geq p$ the bivariant Chern characters $\ch^{2n}(\rho)\in HC^{2n}(\Ic\Ac,\Ac)$ are defined and related by the $S$-operation in bivariant cyclic cohomology $\ch^{2n+2}(\rho) \equiv S\ch^{2n}(\rho)$. We recall briefly the construction of $\ch^{2n}(\rho)$. By the universal properties of the tensor algebra $T(\Ic\Ac)$, the homomorphism $\rho:\Ic\Ac\to \bigl( \begin{smallmatrix} \Ic\Ac & 0 \\ 0 & 0 \end{smallmatrix} \bigr) \subset \Ec^s_+$ lifts to a classifying homomorphism $\rho_*$ through the commutative diagram (\ref{uni})
$$
\vcenter{\xymatrix{
0 \ar[r]  & J(\Ic\Ac) \ar[r] \ar[d]_{\rho_*} & T(\Ic\Ac)  \ar[r] \ar[d]_{\rho_*} & \Ic\Ac \ar[r]  \ar[d]^{\rho} & 0  \\
0 \ar[r] & \Nc^s_+ \ar[r] & \Mc^s_+ \ar[r] & \Ec^s_+ \ar[r] \ar@/_/@{.>}[l]_{\Id\otimes\si} & 0 }}
$$
induced by the linear splitting. It extends to a homomorphism of pro-algebras $\rho_*:\Th(\Ic\Ac)\to \bigl(\begin{smallmatrix} \Ic\Rch & 0 \\ 0 & 0 \end{smallmatrix} \bigr) \subset \Mch^s_+$. The bivariant Chern character $\ch^{2n}(\rho)$ is the composite of the Goodwillie equivalence $\gamma: X(\Th(\Ic\Ac))\to\Omh\Th(\Ic\Ac)$ with the chain maps $\rho_*:\Omh\Th(\Ic\Ac) \to \Omh \Mch^s_+$ and $\chih^{2n}:\Omh \Mch^s_+\to X(\Rch)$. The two non-zero components of $\chih^{2n}$ are given by Eqs. (\ref{chi}) and defined on $2n$ and $(2n+1)$-forms respectively:
$$
\chih^{2n}_0:\Om^{2n}\Mch^s_+ \to \Rch\ ,\qquad \chih^{2n}_1:\Om^{2n+1}\Mch^s_+ \to \Om^1\Rch_{\nat}\ .
$$
The bivariant Chern character is designed to improve the summability degree and can be used to define the higher Chern characters of idempotents and invertibles via the composition
\be
\ch^{2n}_i : \Kt_i(\Ic\hotimes\Ac) \to  HP_i(\Ic\hotimes\Ac) \xrightarrow{\ch^{2n}_i(\rho)}  HP_i(\Ac)\ , \qquad i=0,1\ . \label{high}
\ee
We shall now establish very explicit formulas for these higher characters. Let $\eh$ be an idempotent of the algebra $M_2(\Kc\hotimes\Ic\Rch)^+$, such that $\eh-p_0\in M_2(\Kc\hotimes\Ic\Rch)$. It is well-known (see for example \cite{C1}) that the differential forms
\beq
\ch_{2n}(\eh) &=& (-)^n\frac{(2n)!}{n!}\, \Tr \big((\eh-\frac{1}{2})(\dd \eh \dd \eh)^n\big)\ \in \Om^{2n}(\Ic\Rch) \quad \mbox{for} \ n\geq 1 \non\\
\ch_{0}(\eh) &=& \Tr(\eh-p_0) \ \in \Om^{0}(\Ic\Rch) \label{e}
\eeq
are the components of a $(b+B)$-cycle of even degree over $\Ic\Rch$, i.e. fulfill the relations $B\ch_{2n}(\eh)+ b\ch_{2n+2}(\eh)=0$ for any $n$. Here $\Tr$ is the trace over $\Kc$. In the odd case, any invertible element $\gh\in (\Kc\hotimes\Ic\Rch)^+$ such that $\gh-1\in \Kc\hotimes\Ic\Rch$ gives rise to a $(b+B)$-cycle of odd degree with components
\be
\ch_{2n+1}(\gh)= \frac{(-)^n}{\sqrt{2\pi i}}\,  n!\, \Tr\big(\gh^{-1}\dd \gh(\dd \gh^{-1}\dd \gh)^n\big)\ \in \Om^{2n+1}(\Ic\Rch)\ .\label{g}
\ee
The homology classes of these cycles are of course homotopy invariant. If $\eh\in M_2(\Kc\hotimes\Ic\Rch[0,1])^+$ is a smooth path of idempotents, we define the components of the associated Chern-Simons form as
\be
\cs_{2n+1}(\eh)=(-)^n\frac{(2n)!}{n!}\int_0^1 \Tr \big((-2\eh+1)\sum_{i=0}^{2n}(\dd \eh)^i s\eh (\dd \eh)^{2n+1-i}\big)\label{f}
\ee
in $\Om^{2n+1}(\Ic\Rch)$. Similarly if $\gh\in (\Kc\hotimes\Ic\Rch[0,1])^+$ is a smooth path of invertibles, the components of the Chern-Simons form are for $n\geq 1$
\be
\cs_{2n}(\gh)=\frac{(-)^n}{\sqrt{2\pi i}}\,  (n-1)!\int_0^1 \Tr\big(\gh^{-1}\dd \gh \sum_{i=0}^{n-1}(\dd \gh^{-1}\dd \gh)^i\dd\om (\dd \gh^{-1}\dd \gh)^{n-1-i}\big)\label{u}
\ee
in $\Om^{2n}(\Ic\Rch)$, where $\om=\gh^{-1}s\gh$, and for $n=0$ we set as before $\cs_0(\gh)=\frac{1}{\sqrt{2\pi i}}\int_0^1\Tr(\om)$. Simple algebraic manipulations show that the $(b+B)$-boundaries of the Chern-Simons forms yield the difference of evaluations of the Chern characters at the endpoints:
\beq
&& B\cs_{2n-1}(\eh)+ b\cs_{2n+1}(\eh) = \ch_{2n}(\eh_1)-\ch_{2n}(\eh_0)\ ,\label{tran}\\
&& B\cs_{2n}(\gh)+ b\cs_{2n+2}(\gh) =  \ch_{2n+1}(\gh_1)-\ch_{2n+1}(\gh_0)\ .\non
\eeq
The higher Chern characters (\ref{high}) and their associated Chern-Simons forms are obtained by evaluation of these $(b+B)$-chains on the inclusion homomorphism $\iota_*:\Ic\Rch\hookrightarrow \bigl(\begin{smallmatrix} \Ic\Rch & 0 \\ 0 & 0 \end{smallmatrix} \bigr) \subset \Mch^s_+$ followed by the chain map $\chih^{2n}$ whenever $2n+1\geq p$. 

\begin{lemma}\label{lchcs}
Let $\Ic$ be $p$-summable and $2n+1\geq p$. For any idempotent $\eh\in M_2(\Kc\hotimes\Ic\Rch)^+$ such that $\eh-p_0\in M_2(\Kc\hotimes\Ic\Rch)$, and any invertible $\gh\in (\Kc\hotimes\Ic\Rch)^+$ such that $\gh-1\in \Kc\hotimes\Ic\Rch$, we define the higher Chern characters by the explicit formulas
\beq
\ch_0^{2n}(\eh) &=& \Tr\,(\eh-p_0)^{2n+1}\ ,\label{hch}\\
\ch_1^{2n}(\gh) &=& \frac{1}{\sqrt{2\pi i}}\, \frac{(n!)^2}{(2n)!} \, \Tr\nat\, \gh^{-1}[(\gh-1)(\gh^{-1}-1)]^n\dd \gh\ ,\non
\eeq
where we take concatenation products over $\Ic$ and $\Tr$ is the trace over the $p$-th power of $\Kc\hotimes\Ic$. Then one has $\ch^{2n}_0(\eh) = \chih^{2n}_0\iota_*\ch_{2n}(\eh)$ in $\Rch$ and $\ch^{2n}_1(\gh) = \chih^{2n}_1 \iota_*\ch_{2n+1}(\gh)$ in $\Om^1\Rch_{\nat}$. \\
Similarly, for any idempotent $\eh\in M_2(\Kc\hotimes\Ic\Rch[0,1])^+$ and any invertible $\gh\in (\Kc\hotimes\Ic\Rch[0,1])^+$, we define the higher Chern-Simons forms by the explicit formulas
\beq
&&\cs_1^{2n}(\eh) = \int_0^1\Tr\nat(-2\eh+1)\sum_{i=0}^n(\eh-p_0)^{2i}s\eh (\eh-p_0)^{2(n-i)}\dd\eh\ ,\label{hcs}\\
&&\cs_0^{2n}(\gh) = \frac{1}{\sqrt{2\pi i}}\, \frac{(n!)^2}{(2n)!} \int_0^1 \Tr\, \gh^{-1}[(\gh-1)(\gh^{-1}-1)]^n s\gh \non
\eeq
Then $\cs^{2n}_1(\eh) = \chih^{2n}_1 \iota_*\cs_{2n+1}(\eh)$ in $\Om^1\Rch_{\nat}$ and $\cs^{2n}_0(\gh) \equiv \chih^{2n}_0 \iota_*\cs_{2n}(\gh) \mod \bb$ in $\Rch$. Moreover the transgression relations hold:
$$
\bb\cs^{2n}_1(\eh)=\ch^{2n}_0(\eh_1)-\ch^{2n}_0(\eh_0)\ ,\qquad \nat\dd\cs^{2n}_0(\gh)=\ch^{2n}_1(\gh_1)-\ch^{2n}_1(\gh_0)\ .
$$
\end{lemma}
{\it Proof:} Let us brielfy explain the computation of the cycles $\chih^{2n}_0\iota_*\ch_{2n}(\eh)$ and $\chih^{2n}_1 \iota_*\ch_{2n+1}(\gh)$ associated to idempotents $\eh\in M_2(\Kc\hotimes\Ic\Rch)^+$ and invertibles $\gh\in (\Kc\hotimes\Ic\Rch)^+$. The upper left corner inclusion $\iota_*:\Ic\Rch\hookrightarrow \Mch^s_+$ canonically extends to a unital homomorphism $(\Kc\hotimes\Ic\Rch)^+\to (\Kc\hotimes\Mch^s_+)^+$, and in matrix form we can write
$$
\iota_*\eh=\left(\begin{matrix}
\eh & 0 \\
0 & p_0 \end{matrix} \right)\ ,\qquad \iota_*\gh=\left(\begin{matrix}
\gh & 0 \\
0 & 1 \end{matrix} \right)\ .
$$
Consequently the commutators with the odd multiplier $F=\bigl(\begin{smallmatrix} 0 & 1 \\ 1 & 0 \end{smallmatrix} \bigr)$ read
$$
[F,\iota_*\eh]=\left(\begin{matrix}
0 & p_0-\eh \\
\eh-p_0 & 0 \end{matrix} \right)\ ,\qquad [F,\iota_*\gh]=\left(\begin{matrix}
0 & 1-\gh \\
\gh-1 & 0 \end{matrix} \right)\ .
$$
It is therefore straightforward to evaluate the differential forms $\ch_{2n}(\eh)$ and $\ch_{2n+1}(\gh)$ on the chain map $\chih^{2n}$ given by (\ref{chi}). One finds $\chih^{2n}_0\iota_*\ch_{2n}(\eh)=\ch^{2n}_0(\eh) $ and $\chih^{2n}_1 \iota_*\ch_{2n+1}(\gh)=\ch^{2n}_1(\gh)$. Similarly with the Chern-Simons forms one finds $\chih^{2n}_1 \iota_*\cs_{2n+1}(\eh)=\cs^{2n}_1(\eh)$, whereas by setting $\om = \gh^{-1}s\gh$
\beq
\lefteqn{ \chih^{2n}_0 \iota_*\cs_{2n}(\gh) = \frac{1}{\sqrt{2\pi i}}\, \frac{(n!)^2}{(2n+1)!} \int_0^1 \Tr \big( \om [(\gh^{-1}-1)(\gh-1)]^n + } \non\\
&& (\gh-1)\om[(\gh^{-1}-1)(\gh-1)]^{n-1}(\gh^{-1}-1) + \ldots + [(\gh^{-1}-1)(\gh-1)]^n\om \big) \non
\eeq
coincides with $\cs^{2n}_0(\gh)$ only modulo commutators.\\
Finally, the transgression relations are an immediate consequence of Eqs. (\ref{tran}) and the fact that $\chih^{2n}$ is a chain map from the $(b+B)$-complex over $\Mch^s_+$ to the complex $X(\Rch)$.  \cqfd\\

In any degree $2n+1\geq p$ the Chern characters $\ch^{2n}_0:\Kt_0(\Ic\hotimes\Ac)\to HP_0(\Ac)$ and $\ch^{2n}_1:\Kt_1(\Ic\hotimes\Ac)\to HP_1(\Ac)$ are thus obtained by first lifting idempotents $e\in M_2(\Kc\hotimes\Ic\Ac)^+$ and invertibles $g\in (\Kc\hotimes\Ic\Ac)^+$ to some $\eh\in M_2(\Kc\hotimes\Ic\Rch)^+$ and $\gh\in (\Kc\hotimes\Ic\Rch)^+$, and then taking the cyclic homology classes of $\ch^{2n}_0(\eh)\in \Rch$ and $\ch^{2n}_1(\gh)\in \Om^1\Rch_{\nat}$. Although $\eh$ and $\gh$ are only defined up to homotopy, the above lemma shows these higher Chern characters are well-defined, and moreover independent of the degree $2n$ because the cocycles $\chih^{2n}$ are all related by the transgression relations
$$
\chih^{2n}-\chih^{2n+2}=[\d,\etah^{2n+1}]\ \in \hom(\Omh\Mch^s_+,X(\Rch))\ .
$$
Passing to suspensions, lifting any idempotent $e\in M_2(\Kc\hotimes \Ic S\Ac)^+$ or invertible $g\in (\Kc\hotimes \Ic S\Ac)^+$ gives rise to an odd cycle $\cs^{2n}_1(\eh)\in \Om^1\Rch_{\nat}$ or an even cycle $\cs^{2n}_0(\gh)\in \Rch$. As expected, this is well-defined at the $K$-theory level and compatible with Bott periodicity:

\begin{lemma}\label{lbott2}
Let $\Ic$ be $p$-summable. In any degre $2n+1\geq p$, the Chern-Simons forms define additive maps $\cs^{2n}_1:\Kt_0(\Ic\hotimes{S}\Ac)\to HP_1(\Ac)$ and $\cs^{2n}_0:\Kt_1(\Ic\hotimes{S}\Ac)\to HP_0(\Ac)$, independent of $n$, and compatible with the Bott isomorphisms:
\beq
\cs^{2n}_1\circ\al \equiv \sqrt{2\pi i}\, \ch^{2n}_1 &:&\Kt_1(\Ic\hotimes\Ac)\to HP_1(\Ac)\ , \non\\
\cs^{2n}_0\circ\beta \equiv \sqrt{2\pi i}\, \ch^{2n}_0 &:& \Kt_0(\Ic\hotimes\Ac)\to HP_0(\Ac)\ . \non
\eeq
\end{lemma}
{\it Proof:} Consider an idempotent $\eh\in M_2(\Kc\hotimes \Ic S\Rch)^+$. We have to show the homotopy invariance of the cyclic homology class $\cs^{2n}_1(\eh)$ with respect to $\eh$. This can be shown by direct computation from Formulas (\ref{hcs}). Define the matrix idempotent $\fh=\bigl(\begin{smallmatrix} \eh & 0 \\ 0 & p_0 \end{smallmatrix} \bigr)$. Then one has $s\fh=\bigl(\begin{smallmatrix} s\eh & 0 \\ 0 & 0 \end{smallmatrix} \bigr)$ and $\dd\fh=\bigl(\begin{smallmatrix} \dd\eh & 0 \\ 0 & 0 \end{smallmatrix} \bigr)$, and $\cs^{2n}_1(\eh)$ can be rewritten by means of the operator $F=\bigl(\begin{smallmatrix} 0 & 1 \\ 1 & 0 \end{smallmatrix} \bigr)$ and the supertrace $\tau$:
$$
\cs_1^{2n}(\eh)=(-)^n\int_0^1 \tau\nat (-2\fh+1)\sum_{i=0}^n[F,\fh]^{2i}s\fh [F,\fh]^{2(n-i)}\dd\fh\ .
$$
Now suppose that $\eh$ depends smoothly on an additional parameter $t$. The above integrand may be expressed in terms of the odd differential $\delta=s+\dd+ dt\frac{\d}{\d t} + [F,\ ]$ as
$$
\tau\nat (-2\fh+1)\sum_{i=0}^n[F,\fh]^{2i}s\fh [F,\fh]^{2(n-i)}\dd\fh=\frac{-1}{n+1}\, \tau\nat \, \fh(\delta\fh)^{2n+1}|_{s,\dd}
$$
where $|_{s,\dd}$ means that we select the terms containing only $s$, $\dd$ and not $dt$. Because $\tau\nat$ is a supertrace, the cocycle property $(s+\dd+dt\frac{\d}{\d t})\tau\nat \, \fh(\delta\fh)^{2n+1}=\tau\nat \, (\delta\fh)^{2n+2}=0$ holds, and projecting this relation on $s,\dd,dt$ yields
$$
s\, \tau\nat \, \fh(\delta\fh)^{2n+1}|_{\dd,dt}+\tau\nat \,\dd( \fh(\delta\fh)^{2n+1})|_{s,dt} + dt \frac{\d}{\d t}\tau\nat \, \fh(\delta\fh)^{2n+1}|_{s,\dd}=0\ .
$$
This may be rephrased as
\beq
\lefteqn{\frac{\d}{\d t}\big(\Tr\nat (-2\eh+1)\sum_{i=0}^n(\eh-p_0)^{2i}s\eh (\eh-p_0)^{2(n-i)}\dd\eh \big) }\non\\
&\equiv& s \big( \Tr\nat (-2\eh+1)\sum_{i=0}^n(\eh-p_0)^{2i}\frac{\d\eh}{\d t} (\eh-p_0)^{2(n-i)}\dd\eh\big) \mod \nat\dd\ ,\non
\eeq
and integration over the current $\int_0^1$ shows the homotopy invariance of the class $\cs_1^{2n}(\eh)$. Hence the map $\cs_1^{2n}:\Kt_0(\Ic\hotimes{S}\Ac)\to HP_1(\Ac)$ is well-defined. Its compatibility with Bott periodicity can be established, without computation, as follows. Let $g\in (\Kc\hotimes \Ic\Ac)^+$ be an invertible element and  $e=\al(g)\in M_2(\Kc\hotimes \Ic S\Ac)^+$ its idempotent image under the Bott isomorphism. Choose an invertible lift $\gt\in (\Kc\hotimes S\Th(\Ic\Ac))^+$ of $g$ and an idempotent lift $\et\in M_2(\Kc\hotimes S\Th(\Ic\Ac))^+$ of $e$. The differential forms $\ch_{2n+1}(\gt)$ and $\cs_{2n+1}(\et)$ in $\Om^{2n+1}\Th(\Ic\Ac)$ defined by (\ref{g}) and (\ref{f}) are the components of two $(b+B)$-cycles $\ch_*(\gt)$ and $\cs_*(\et)$, whose projections on the odd part of the complex $X(\Th(\Ic\Ac))$ are
$$
\nat\ch_1(\gt)=\frac{1}{\sqrt{2\pi i}} \Tr\nat\gt^{-1}\dd\gt\ ,\qquad \nat\cs_1(\et)=\int_0^1\Tr\nat(-2\et+1)s\et\dd\et)\ .
$$
By Lemma \ref{lbott1}, the cycles $\sqrt{2\pi i}\nat\ch_1(\gt)$ and $\nat\cs_1(\et)$ are homologous. But we know that the projection $\Omh\Th(\Ic\Ac)\to X(\Th(\Ic\Ac))$ is a homotopy equivalence, with inverse the Goodwillie map $\gamma$. Hence $\sqrt{2\pi i}\, \ch_*(\gt)$ and $\cs_*(\et)$ are $(b+B)$-homologous in $\Omh\Th(\Ic\Ac)$. Finally, it remains to observe that under the homomorphism $\rho_*:\Th\Ac\to \Mc^s_+$, the invertible $\rho_*\gt=\bigl(\begin{smallmatrix} \gh & 0 \\ 0 & 1 \end{smallmatrix} \bigr)$ gives a choice of lifting $\gh\in (\Kc\hotimes \Ic\Rch)^+$ and the idempotent $\rho_*(\et)=\bigl(\begin{smallmatrix} \eh & 0 \\ 0 & p_0 \end{smallmatrix} \bigr)$ gives a choice of lifting $\eh\in M_2(\Kc\hotimes \Ic\Rch)^+$. Because the cycles $\sqrt{2\pi i}\,\chih^{2n}\rho_*\ch_*(\gt)$ and $\chih^{2n}\rho_*\cs_*(\et)$ are homologous in $X(\Rch)$, we have $\sqrt{2\pi i}\, \ch_1^{2n}(\gh)\equiv \cs_1^{2n}(\eh)$ in $HP_1(\Ac)$.\\
We proceed similarly with the map $\cs_0^{2n}$. Let $\gh\in (\Kc\hotimes \Ic S\Rch)^+$ be an invertible. We have to show the homotopy invariance of the cyclic homology class $\cs^{2n}_0(\gh)$. Define the invertible matrix $\uh=\bigl(\begin{smallmatrix} \gh & 0 \\ 0 & 1 \end{smallmatrix} \bigr)$. Then $s\uh=\bigl(\begin{smallmatrix} s\gh & 0 \\ 0 & 0 \end{smallmatrix} \bigr)$ and one has
$$
\cs_0^{2n}(\gh) = \frac{(-)^n}{\sqrt{2\pi i}}\frac{(n!)^2}{(2n)!}\int_0^1 \tau (\uh^{-1}([F,\uh][F,\uh^{-1}])^n s\uh)\ .
$$
Now suppose that $\gh$ is a smooth family of invertibles depending on an additional parameter $t$. The above integrand may be expressed in terms of the odd derivation $\delta=s+ dt\frac{\d}{\d t} + [F,\ ]$ as
$$
\tau\, \uh^{-1}([F,\uh][F,\uh^{-1}])^n s\uh \equiv \frac{(-)^n}{2n+1}\, \tau  (\uh^{-1}\delta\uh)^{2n+1}|_{s} \mod \bb\ .
$$
One has the relation $(s+dt\frac{\d}{\d t})\tau(\uh^{-1}\delta\uh)^{2n+1} =-\tau (\uh^{-1}\delta\uh)^{2n+2}\equiv 0\mod\bb$, hence by projection on $s,dt$
$$
s\, \tau  (\uh^{-1}\delta\uh)^{2n+1}|_{dt}+ dt \frac{\d}{\d t}\tau (\uh^{-1}\delta\uh)^{2n+1}|_{s} \equiv 0 \mod \bb\ .
$$
This may be rephrased as
$$
\frac{\d}{\d t}\big(\Tr\, \gh^{-1}[(\gh-1)(\gh^{-1}-1)]^ns\gh \big) \equiv s \big( \Tr\, \gh^{-1}[(\gh-1)(\gh^{-1}-1)]^n \frac{\d\gh}{\d t} \big) \mod \bb\ ,\non
$$
and integration over the current $\int_0^1$ shows the homotopy invariance of the class $\cs_0^{2n}(\gh)$. Hence the map $\cs_0^{2n}:\Kt_1(\Ic\hotimes{S}\Ac)\to HP_0(\Ac)$ is well-defined. Its compatibility with Bott periodicity is be established as before, replacing invertibles by idempotents and conversely. \cqfd

\section{Multiplicative $K$-theory}\label{sm}

Let $\Ac$ and $\Ic$ be Fr\'echet $m$-algebras, $\Ic$ being $p$-summable. We shall define the multiplicative $K$-theory groups $MK^{\Ic}_n(\Ac)$ in any degree $n\in\zz$. They are intermediate between the topological $K$-theory $\Kt_n(\Ic\hotimes\Ac)$ and the non-periodic cyclic homology $HC_n(\Ac)$. Recall from section \ref{scy} that if $0\to \Jc \to \Rc \to \Ac \to 0$ is any quasi-free extension, $HC_n(\Ac)$ is computed by the quotient complex $X_n(\Rc,\Jc)=X(\Rc)/F^n_{\Jc}X(\Rc)$ induced by the $\Jc$-adic filtration:
$$
HC_n(\Ac)=H_{n+2\zz}(X_n(\Rc,\Jc))\ ,\qquad \forall n\in\zz\ .
$$
Of course $HC_n(\Ac)$ vanishes whenever $n<0$. Multiplicative $K$-theory classes are represented by idempotents or invertibles whose higher Chern characters (Lemma \ref{lchcs}) can be transgressed up to a certain order. As before we adopt the notation $\Ic\Ac=\Ic\hotimes\Ac$ and $\Ic\Rch=\Ic\hotimes\Rch$, where $\Rch$ is the $\Ic$-adic completion of $\Rc$.

\begin{definition}
Let $0\to \Jc \to \Rc \to \Ac \to 0$ be any quasi-free extension of Fr\'echet $m$-algebras, and let $\Ic$ be a $p$-summable Fr\'echet $m$-algebra. Choose an integer $q$ such that $2q+1\geq p$. We define the multiplicative $K$-theory $MK^{\Ic}_n(\Ac)$, in any even degree $n=2k\in \zz$, as the set of equivalence classes of pairs $(\eh,\te)$ such that $\eh\in M_2(\Kc\hotimes\Ic\Rch)^+$ is an idempotent and $\te\in X_{n-1}(\Rc,\Jc)$ is a chain of odd degree related by the transgression formula
$$
\ch^{2q}_0(\eh)=\bb \te\ \in X_{n-1}(\Rc,\Jc)\ .
$$
Two pairs $(\eh_0,\te_0)$ and $(\eh_1,\te_1)$ are equivalent if and only if there exists an idempotent $\eh\in M_2(\Kc\hotimes\Ic\Rch[0,1])^+$ whose evaluation yields $\eh_0$ and $\eh_1$ at the endpoints, and a chain $\la\in X_{n-1}(\Rc,\Jc)$ of even degree such that
$$
\te_1-\te_0=\cs^{2q}_1(\eh)+ \nat\dd \la\ \in X_{n-1}(\Rc,\Jc)\ .
$$
In the same way, we define the multiplicative $K$-theory $MK^{\Ic}_{n}(\Ac)$, in any odd degree $n=2k+1\in \zz$, as the set of equivalence classes of pairs $(\gh,\te)$ such that $\gh\in (\Kc\hotimes\Ic\Rch)^+$ is an invertible and $\te\in X_{n-1}(\Rc,\Jc)$ is a chain of even degree related by the transgression formula
$$
\ch^{2q}_1(\gh)=\nat\dd\te\ \in X_{n-1}(\Rc,\Jc)\ .
$$
Two pairs $(\gh_0,\te_0)$ and $(\gh_1,\te_1)$ are equivalent if and only if there exists an invertible $\gh\in (\Kc\hotimes\Ic\Rch[0,1])^+$ whose evaluation yields $\gh_0$ and $\gh_1$ at the endpoints, and a chain $\la\in X_{n-1}(\Rc,\Jc)$ of odd degree such that
$$
\te_1-\te_0=\cs^{2q}_0(\gh)+\bb \la\ \in X_{n-1}(\Rc,\Jc)\ .
$$
We will prove in Proposition \ref{pmul} that for any $n\in\zz$ the set $MK_n^{\Ic}(\Ac)$ depends neither on the degree $2q+1\geq p$ chosen to represent the Chern characters nor on the quasi-free extension $\Rc$.
\end{definition}

Recall from section \ref{scy} that the homology $H_{n-1+ 2\zz}(X_{n}(\Rc,\Jc))$ is the non-commutative de Rham homology $HD_{n-1}(\Ac)$. Hence the transgression relation $\ch^{2q}_0(\eh)=\bb \te\in X_{n-1}(\Rc,\Jc)$ exactly means that the class of $\ch^{2q}_0(\eh)$ vanishes in $HD_{n-2}(\Ac)$. Similarly in the odd case, $\ch^{2q}_1(\gh)=\nat\dd\te \in X_{n-1}(\Rc,\Jc)$ is equivalent to $\ch^{2q}_1(\gh)\equiv 0$ in $HD_{n-2}(\Ac)$. Since the complexes $X_{n-1}(\Rc,\Jc)$ vanish for $n\leq 0$, we immediately deduce that $MK^{\Ic}_n(\Ac)$ coincides with the topological $K$-theory $\Kt_n(\Ic\hotimes\Ac)$ whenever $n\leq 0$.\\
As in the case of topological $K$-theory, define an addition on $MK^{\Ic}_n(\Ac)$ by direct sum of idempotents and invertibles as follows ($c$ is the permutation matrix (\ref{mat})):
\beq
\mbox{even case:} && (\eh,\te) + (\eh',\te') = (c(\eh\oplus \eh')c, \te+\te')\ ,\non\\
\mbox{odd case:} && (\gh,\te) + (\gh',\te') = (\gh\oplus \gh', \te+\te')\ . \non
\eeq
This turns $MK^{\Ic}_n(\Ac)$ into a semigroup, the unit being represented by $(p_0,0)$ in the even case and $(1,0)$ in the odd case.

\begin{lemma}\label{lab}
$MK_n^{\Ic}(\Ac)$ is an abelian group for any $n\in\zz$.
\end{lemma}
{\it Proof:} We first need to recall the proof that $\Kt_0(\Ic\hotimes\Ac)$ is a group. Let $e\in M_2(\Kc\hotimes\Ic\Ac)^+$ be an idempotent, with $e-p_0\in M_2(\Kc\hotimes\Ic\Ac)$. The idempotent $1-e$ is orthogonal to $e$, as $e(1-e)=(1-e)e=0$. If $X\in M_2(\cc)$ is the permutation matrix $\bigl( \begin{smallmatrix}
          0 & 1 \\
          1 & 0 \end{smallmatrix} \bigr)$, we claim that the idempotent 
$$
X(1-e)X\in M_2(\Kc\hotimes\Ic\Ac)^+\ ,\quad\mbox{with}\quad X(1-e)X-p_0\in M_2(\Kc\hotimes\Ic\Ac)\ ,
$$
represents the inverse class of $e$. Indeed, we shall construct a homotopy between the direct sum $c(e\oplus X(1-e)X)c\in M_4(\Kc\hotimes\Ic \Ac)^+$ and the unit
$$
\tilde{p}_0= \left( \begin{matrix}
1 & 0 & 0 & 0 \\
0 & 1 & 0 & 0 \\
0 & 0 & 0 & 0 \\
0 & 0 & 0 & 0 \end{matrix} \right)\in M_4(\Kc\hotimes\Ic\Ac)^+\ .
$$
Choose a smooth real-valued function $\xi\in\cinf[0,1]$ ranging from $\xi(0)=0$ to $\xi(1)=\pi/2$, and consider the paths of invertible matrices
$$
R_{23}=\left( \begin{matrix}
1 & 0 & 0 & 0 \\
0 & \cos\xi & \sin\xi & 0 \\
0 & \sin\xi & -\cos\xi & 0 \\
0 & 0 & 0 & 1 \end{matrix} \right)\ ,\quad R_{14}=\left( \begin{matrix}
\cos\xi & 0 & 0 & \sin\xi \\
0 & 1 & 0 & 0 \\
0 & 0 & 1 & 0 \\
\sin\xi & 0 & 0 & -\cos\xi \end{matrix} \right)\ .
$$
A direct computation shows that the idempotents $R_{23}(t)^{-1}\bigl(\begin{smallmatrix}
          e & 0 \\
          0 & 0 \end{smallmatrix} \bigr)R_{23}(t)$ and $R_{14}(t)^{-1}\bigl(\begin{smallmatrix} 1-e & 0 \\ 0 & 0 \end{smallmatrix} \bigr)R_{14}(t)$ are orthogonal for any $t\in[0,1]$. Hence the sum
$$
f= R_{23}^{-1}\left(\begin{matrix}
          e & 0 \\
          0 & 0 \end{matrix} \right)R_{23}+R_{14}^{-1}\left(\begin{matrix}
          1-e & 0 \\
          0 & 0 \end{matrix} \right)R_{14}\ \in M_4(\Kc\hotimes\Ic\Ac[0,1])^+
$$
is an idempotent path such that $f-\tilde{p}_0\in M_4(\Kc\hotimes\Ic\Ac[0,1])$, and interpolates $f_0=\tilde{p}_0$ and $f_1=c(e\oplus X(1-e)X)c$. This shows that $\Kt_0(\Ic\hotimes\Ac)$ is a group. It is abelian because a direct sum $e\oplus e'$ can be connected via a smooth path (by conjugation with respect to rotation matrices) to $e'\oplus e$.\\
Now fix an integer $2q+1\geq p$ and let $(\eh,\te)$ represent an element of $MK_{n}^{\Ic}(\Ac)$ of even degree $n=2k$. Hence $\eh\in M_2(\Kc\hotimes\Ic\Rch)^+$ is an idempotent such that $\ch^{2q}_0(\eh)=\bb\te$ in the quotient complex $X_{n-1}(\Rc,\Jc)$. Consider the smooth idempotent path $\fh\in M_4(\Kc\hotimes\Ic\Rch[0,1])^+$ constructed as above replacing $\Ac$ by its extension $\Rch$ and $e$ by $\eh$. It provides an interpolation between $\fh_0=\tilde{p}_0$ and $\fh_1=c(\eh\oplus X(1-\eh)X)c$. We guess that the inverse of $(\eh,\te)$ is represented by the pair $(X(1-\eh)X, \cs^{2q}_1(\fh)-\te)$. Indeed, one has
$$
\bb \cs^{2q}_1(\fh)= \ch^{2q}_0(\fh_1)-\ch^{2q}_0(\fh_0)=\ch^{2q}_0(\eh)+\ch^{2q}_0(X(1-\eh)X)\ ,
$$
so that $\ch^{2q}_0(X(1-\eh)X)=\bb(\cs^{2q}_1(\fh)-\te)$ in the complex $X_{n-1}(\Rc,\Jc)$ and the pair $(X(1-\eh)X, \cs^{2q}_1(\fh)-\te)$ represents a class in $MK_{n}^{\Ic}(\Ac)$. Moreover, the sum
$$
(\eh,\te)+(X(1-\eh)X, \cs^{2q}_1(\fh)-\te)=(c(\eh\oplus X(1-\eh)X)c,\cs^{2q}_1(\fh))
$$
is equivalent to the unit $(\tilde{p}_0,0)$ because $\fh$ provides the interpolating idempotent. Hence $MK_{n}^{\Ic}(\Ac)$, $n=2k$ is a group. Abelianity is shown as for topological $K$-theory, by means of another interpolation between the idempotents $c(\eh\oplus \eh')c$ and $c(\eh'\oplus \eh)c$ with the property that its Chern-Simons form $\cs_1$ vanishes.\\
One proceeds similarly in the odd case. Let $(\gh,\te)$ represent an element of $MK^{\Ic}_{n}(\Ac)$ of odd degree $n=2k+1$. Hence $\ch^{2q}_1(\gh)=\nat\dd\te$ in $X_{n-1}(\Rc,\Jc)$. Define an invertible path $\uh\in M_2(\Kc\hotimes\Ic\Rch[0,1])^+$ by means of the rotation matrix $R=\bigl( \begin{smallmatrix}
          \cos\xi & \sin\xi \\
          -\sin\xi & \cos\xi \end{smallmatrix} \bigr)$:
$$
\uh=\left( \begin{matrix}
          \gh & 0 \\
          0 & 1 \end{matrix} \right) R^{-1} \left( \begin{matrix}
          \gh^{-1} & 0 \\
          0 & 1 \end{matrix} \right) R\ .
$$
Then $\uh-1\in M_2(\Kc\hotimes\Ic\Rch[0,1])$, and $\uh$ provides a smooth homotopy between the invertibles $\uh_0=1$ and $\uh_1 = \bigl( \begin{smallmatrix}
          \gh & 0 \\
          0 & \gh^{-1} \end{smallmatrix} \bigr)$ (the same argument shows that $\Kt_1(\Ic\hotimes\Ac)$ is an abelian group). We guess that the inverse class of $(\gh,\te)$ is represented by the pair $(\gh^{-1},\cs^{2q}_0(\uh)-\te)$. Indeed, one has
$$
\nat\dd\cs^{2q}_0(\uh) = \ch^{2q}_1(\uh_1) - \ch^{2q}_1(\uh_0) = \ch^{2q}_1(\gh) + \ch^{2q}_1(\gh^{-1}) 
$$
so that $\ch_1^{2q}(\gh^{-1})= \nat\dd(\cs_0^{2q}(\uh) - \te)$ in $X_{n-1}(\Rc,\Jc)$ and $(\gh^{-1},\cs^{2q}_0(\uh)-\te)$ represents a class in $MK^{\Ic}_{n}(\Ac)$. Moreover, the sum
$$
(\gh,\te) + (\gh^{-1},\cs^{2q}_0(\uh)-\te) = (\gh\oplus \gh^{-1} , \cs^{2q}_0(\uh))
$$
is equivalent to the unit $(1,0)$ through the interpolating invertible $\uh$. Hence $MK^{\Ic}_{n}(\Ac)$, $n=2k+1$ is a group as claimed. Abelianity is shown once again by means of rotation matrices.\cqfd\\

\begin{remark}\label{rcan}\textup{We know that two different liftings of a given idempotent $e\in M_2(\Kc\hotimes\Ic\Ac)^+$ are always homotopic in $M_2(\Kc\hotimes\Ic\Rch)^+$. Hence choosing the universal free extension $\Rc=T\Ac$ allows to represent any multiplicative $K$-theory class of even degree by a pair $(\eh,\te)$ where $\eh\in M_2(\Kc\hotimes\Ic\Th\Ac)^+$ is the \emph{canonical lift} of some idempotent $e$. Moreover, the transgression formula established in the proof of Lemma \ref{lbott2} shows that two such pairs $(\eh_0,\te_0)$ and $(\eh_1,\te_1)$ are equivalent if and only if $e_0$ and $e_1$ can be joined by an idempotent $e\in M_2(\Kc\hotimes\Ic\Ac[0,1])^+$ such that $\te_1-\te_0\equiv \cs^{2q}_1(\eh)\mod \nat\dd$, where $\eh\in M_2(\Kc\hotimes\Ic\Th\Ac[0,1])^+$ is the canonical lift of $e$. The same is true with invertibles: any multiplicative $K$-theory class of odd degree may be represented by a pair $(\gh,\te)$ where $\gh\in (\Kc\hotimes\Ic\Th\Ac)^+$ is the \emph{canonical lift} of some invertible $g\in (\Kc\hotimes\Ic\Ac)^+$. Two such pairs $(\gh_0,\te_0)$ and $(\gh_1,\te_1)$ are equivalent if and only if $g_0$ and $g_1$ can be joined by an invertible $g\in (\Kc\hotimes\Ic\Ac[0,1])^+$ such that $\te_1-\te_0\equiv \cs^{2q}_0(\gh)\mod \bb$, where $\gh\in (\Kc\hotimes\Ic\Th\Ac[0,1])^+$ is the canonical lift of $g$.  }
\end{remark}

The particular case $\Ic=\cc$ is essentially equivalent Karoubi's definition of multiplicative $K$-theory \cite{K1,K2}. The groups $MK^{\Ic}_n(\Ac)$ are designed to fit in a long exact sequence
\be
\Kt_{n+1}(\Ic\hotimes\Ac) \to HC_{n-1}(\Ac) \stackrel{\delta}{\to} MK^{\Ic}_n(\Ac)  \to \Kt_n(\Ic\hotimes\Ac)  \to HC_{n-2}(\Ac)  \label{ex}
\ee
The map $\Kt_n(\Ic\hotimes \Ac)\to HC_{n-2}(\Ac)$ corresponds to the composition of the Chern character $\Kt_n(\Ic\hotimes \Ac)\to HP_n(\Ac)$ with the natural map $HP_n(\Ac)\to HC_{n-2}(\Ac)$ induced by the projection $\Xh(\Rc,\Jc)\to X_{n-2}(\Rc,\Jc)$. The map $MK^{\Ic}_n(\Ac)\to \Kt_n(\Ac)$ is the forgetful map, which sends a pair $(\eh,\te)$ or $(\gh,\te)$ respectively on its image $e$ or $g$ under the projection homomorphism $\Rch\to\Ac$. The connecting map $\delta:HC_{n-1}(\Ac)\to MK^{\Ic}_n(\Ac)$ sends a cycle $\te\in X_{n-1}(\Rc,\Jc)$ to
\be
\delta(\te) = \left\{ \begin{array}{ll} 
(p_0, \sqrt{2\pi i}\,\te) & \mbox{$n$ even,}\\
(1, \sqrt{2\pi i}\,\te) & \mbox{$n$ odd.} \end{array} \right.
\ee
There is also an additive Chern character map $\ch_n: MK^{\Ic}_n(\Ac)\to HN_n(\Ac)$ defined in all degrees $n\in\zz$, with values in \emph{negative} cyclic homology. Recall that the latter is the homology in degree $n$ mod $2$ of the subcomplex $F^{n-1}\Xh(\Rc,\Jc)=\ker (\Xh(\Rc,\Jc)\to X_{n-1}(\Rc,\Jc))$:
$$
HN_n(\Ac)=H_{n+2\zz}(F^{n-1}\Xh(\Rc,\Jc))\ .
$$
Hence in particular $HN_n(\Ac)=HP_n(\Ac)$ whenever $n\leq 0$, and $HP_*$, $HC_*$ $HN_*$ are related by the $SBI$ long exact sequence (section \ref{scy}). To define the Chern character $\ch_{n}: MK^{\Ic}_{n}(\Ac)\to HN_{n}(\Ac)$ in even degree $n=2k$, we first have to choose an integer $2q+1\geq p$. Then, a multiplicative $K$-theory class of degree $n$ is represented by a pair $(\eh,\te)$, such that the transgression formula $\ch^{2q}_0(\eh)=\bb \te$ holds in $X_{n-1}(\Rc,\Jc)$. Choose an arbitrary lifting $\tilde{\te}\in \Xh(\Rc,\Jc)$ of $\te$, and define the negative Chern character as
\be
\ch_{n}(\eh,\te)= \ch^{2q}_0(\eh)- \bb \tilde{\te}\ \in F^{n-1}\Xh(\Rc,\Jc)\ .
\ee
It is clearly closed, and its negative cyclic homology class does not depend on the choice of lifting $\tilde{\te}$, since the difference of two such liftings lies in the subcomplex $F^{n-1}\Xh(\Rc,\Jc)$. We will show in the proposition below that it does not depend on the representative $(\eh,\te)$ of the $K$-theory class, nor on the integer $2q+1\geq p$. In odd degree $n=2k+1$, the Chern character $\ch_{n}: MK^{\Ic}_{n}(\Ac)\to HN_{n}(\Ac)$ is defined exactly in the same way: take a representative $(\gh,\te)$ of a multiplicative $K$-theory class, with $\ch^{2q}_1(\gh)=\nat\dd \te$ in $X_{n-1}(\Rc,\Jc)$. Then if $\tilde{\te}\in \Xh(\Rc,\Jc)$ denotes an arbitrary lifting of $\te$, the cycle
\be
\ch_{n}(\gh,\te)= \ch^{2q}_1(\gh)- \nat\dd \tilde{\te}\ \in F^{n-1}\Xh(\Rc,\Jc)
\ee
defines a negative cyclic homology class. The following proposition shows the compatibility between the $K$-theory exact sequence (\ref{ex}) and the $SBI$ exact sequence (\ref{SBI}), through the various Chern character maps.

\begin{proposition}\label{pmul}
Let $\Ac$ and $\Ic$ be Fr\'echet $m$-algebras, such that $\Ic$ is $p$-summable. Then one has a commutative diagram with long exact rows
\be
\vcenter{\xymatrix{
\Kt_{n+1}(\Ic\hotimes\Ac) \ar[r] \ar[d] & HC_{n-1}(\Ac) \ar[r]^{\delta} \ar@{=}[d] & MK^{\Ic}_n(\Ac)  \ar[r] \ar[d] & \Kt_n(\Ic\hotimes\Ac)  \ar[d] \\
HP_{n+1}(\Ac) \ar[r]^S & HC_{n-1}(\Ac) \ar[r]^{\widetilde{B}} & HN_n(\Ac) \ar[r]^I & HP_n(\Ac) }} \label{mul}
\ee
where $\widetilde{B}$ is the connecting map of the $SBI$ sequence rescaled by a factor $-\sqrt{2\pi i}$.
\end{proposition}
{\it Proof:} We show the exactness of the sequence (\ref{ex}) in the case of even degree $n=2k$ (the odd case is completely similar):
$$
\Kt_{1}(\Ic\hotimes\Ac) \stackrel{\ch_1}{\longrightarrow} HC_{n-1}(\Ac) \stackrel{\delta}{\to} MK^{\Ic}_{n}(\Ac)  \stackrel{\iota}{\to} \Kt_0(\Ic\hotimes\Ac)  \stackrel{\ch_0}{\longrightarrow} HC_{n-2}(\Ac)\ .
$$
Fix once and for all an integer $2q+1\geq p$ to represent the Chern characters. We first have to check that the maps $\delta$ and $\iota$ are well-defined. Let $\te_0$ and $\te_1=\te_0+\nat\dd \la$ be two homologous odd cycles in $X_{n-1}(\Rc,\Jc)$ representing the same cyclic homology class $[\te]\in HC_{n-1}(\Ac)$. Their images by $\delta$ are respectively $(p_0,\sqrt{2\pi i}\, \te_0)$ and $(p_0,\sqrt{2\pi i}\, \te_1)$, which obviously represent the same class in $MK^{\Ic}_{n}(\Ac)$ by virtue of the equivalence relation $\te_1-\te_0=\nat\dd\la$ (take the constant idempotent $p_0\in M_2(\Kc\hotimes\Ic\Rch[0,1])^+$ as interpolation, with $\cs^{2q}_1(p_0)=0$). Hence $\delta$ is well-defined.\\
Now take two equivalent pairs $(\eh_0,\te_0)$ and $(\eh_1,\te_1)$ representing the same element in $MK^{\Ic}_{n}(\Ac)$. In particular, the idempotents $\eh_0$ and $\eh_1$ are smoothly homotopic and their projections $e_0,e_1\in M_2(\Kc\hotimes\Ic\Ac)^+$ define the same class in $\Kt_0(\Ic\hotimes\Ac)$. Since $\iota(\eh_0,\te_0)=e_0$ and $\iota(\eh_1,\te_1)=e_1$, the map $\iota$ is well-defined.\\

\noindent \emph{Exactness at $HC_{n-1}(\Ac)$:} Let $[g]\in \Kt_{1}(\Ic\hotimes\Ac)$ be a class represented by an invertible element $g\in (\Kc\hotimes\Ic\Ac)^+$, and consider its idempotent image $\al(g)\in M_2(\Kc\hotimes \Ic S\Ac)^+$ under the Bott isomorphism $\al:\Kt_{1}(\Ic\hotimes\Ac)\to \Kt_{0}(\Ic\hotimes{S}\Ac)$. Choose any invertible lift $\gh\in (\Kc\hotimes\Ic\Rch)^+$ and any idempotent lift $\widehat{\al(g)} \in M_2(\Kc\hotimes \Ic S\Rch)^+$. By Lemma \ref{lbott2}, we have the equality of periodic cyclic homology classes
$$
\ch^{2q}_1(\gh)\equiv \frac{1}{\sqrt{2\pi i}} \cs^{2q}_1(\widehat{\al(g)})\  \in HP_1(\Ac)\ ,
$$
hence this equality also holds in $HC_{n-1}(\Ac)$. It follows that $\delta(\ch_1(g))$ is represented by 
$$
\delta\Big(\frac{1}{\sqrt{2\pi i}} \cs^{2q}_1(\widehat{\al(g)})\Big)= \big(p_0,\cs^{2q}_1(\widehat{\al(g)})\big)\ .
$$
But the idempotent path $\widehat{\al(g)}$ evaluated at the endpoints is $p_0$, so that the pairs $(p_0,\cs^{2q}_1(\widehat{\al(g)}))$ and $(p_0,0)$ are equivalent in $MK^{\Ic}_n(\Ac)$. Hence $\delta\circ \ch_1=0$.\\
Now let a class $[\te]\in HC_{n-1}(\Ac)$ be in the kernel of $\delta$. It means that the pair $\delta(\te)=(p_0,\sqrt{2\pi i}\, \te)$ is equivalent to $(p_0,0)$. Hence there exists an idempotent $\eh\in M_2(\Kc\hotimes \Ic S\Rch)^+$ and a chain $\la$ such that $\sqrt{2\pi i}\, \te= \cs^{2q}_1(\eh)+\nat\dd\la$ in $X_{n-1}(\Rc,\Jc)$. By Bott periodicity, there exists an element $[g]\in \Kt_{1}(\Ic\hotimes\Ac)$ such that $\cs^{2q}_1(\eh)\equiv \sqrt{2\pi i}\, \ch^{2q}_1(\gh)$ in $HC_{n-1}(\Ac)$, where $\gh\in (\Kc\hotimes \Ic\Rch)^+$ is any invertible lift of $g$. Whence the equality of cylic homology classes $[\te]\equiv \ch^{2q}_1(\gh)$. It follows that $\ker\delta=\im \ch_1$.\\

\noindent \emph{Exactness at $MK^{\Ic}_{n}(\Ac)$:} Let $[\te]\in HC_{n-1}(\Ac)$ be any cyclic homology class. Then $\delta(\te)=(p_0,\sqrt{2\pi i}\, \te)$, and $\iota(\delta(\te))=p_0$ is the zero-class in topological $K$-theory. Therefore $\iota\circ \delta=0$.\\
Now let $(\eh,\te)\in MK^{\Ic}_n(\Ac)$ be in the kernel of $\iota$: it means that $\eh$ is smoothly homotopic to $p_0$. Hence there exists an idempotent $\fh\in M_2(\Kc\hotimes\Ic\Rch[0,1])^+$, with evaluations $\fh_0=\eh$ and $\fh_1=p_0$, and the pair $(\eh,\te)$ is equivalent to $(p_0,\te+\cs^{2q}_1(\fh))$. Remark that the odd chain $\te+\cs^{2q}_1(\fh)\in X_{n-1}(\Rc,\Jc)$ is closed (indeed, $\bb\te=\ch^{2q}_0(\eh)$ and $\bb\cs^{2q}_1(\fh)=-\ch^{2q}_0(\eh)$), and we can write
$$
(p_0,\te+\cs^{2q}_1(\fh))=\delta\Big(\frac{1}{\sqrt{2\pi i}}(\te+ \cs^{2q}_1(\fh))\Big)\ .
$$
It follows that $\ker\iota  = \im \delta$.\\

\noindent \emph{Exactness at $\Kt_{0}(\Ic\hotimes\Ac)$:} Let $(\eh,\te)\in MK^{\Ic}_{n}(\Ac)$ represent any multiplicative $K$-theory class. Then $\ch^{2q}_0(\eh)\equiv 0$ in non-commutative de Rham homology $HD_{n-2}(\Ac)$, and therefore also in $HC_{n-2}(\Ac)$. Thus, the Chern character of $\iota(\eh,\te)=e$ vanishes and $\ch_0\circ\iota=0$.\\
Now let $[e]\in \Kt_{0}(\Ic\hotimes\Ac)$ be in the kernel of $\ch_0$. We know from section \ref{scy} that the natural map $HD_{n-2}(\Ac)\to HC_{n-2}(\Ac)$ is injective, so that $\ch^{2q}_0(\eh)\equiv 0$ in $HD_{n-2}(\Ac)$ for any idempotent lift $\eh\in M_2(\Kc\hotimes \Ic\Rch)^+$. Hence there exists an odd chain $\te\in X_{n-1}(\Rc,\Jc)$ such that $\ch^{2q}_0(\eh)= \bb\te$, and $e=\iota(\eh,\te)$. This shows that $\ker \ch_0= \im\iota$.\\

\noindent Let us now show the independence of multiplicative $K$-theory upon the choice of degree $2q+1\geq p$. To this end, write $(\eh, \te)^q\in MK^{\Ic}_{n}(\Ac)^q$ for a representative of a class obtained using the higher Chern character $\ch_0^{2q}(\eh)=\bb\te$ of degree $2q$. We shall construct a map $MK^{\Ic}_{n}(\Ac)^q\to MK^{\Ic}_{n}(\Ac)^{q+1}$ which turns out to be an isomorphism. Let $\rho:\Ic\Ac\to \Ec^s\triangleright \Ic^s\Ac$ be the canonical $p$-summable quasihomomorphism of even degree considered in section \ref{st}, for the construction of the higher Chern characters. Recall that $\Ec=\Ic^+\hotimes\Ac$ with extension $\Mc=\Ic^+\hotimes\Rc$. From Proposition \ref{padic}, we know that the chain maps $\chi^{2q}:\Omh(\Mch^s_+)\to X(\Rch)$ associated to $\rho$ are related in successive degrees by the transgression formula involving the eta-cochain $\chih^{2q}-\chih^{2q+2}=[\d, \etah^{2q+1}]$. More precisely:
\beq
\chih^{2q}_0-\chih^{2q+2}_0 &=& \bb \etah^{2q+1}_1 + \etah^{2q+1}_0(b+B)\ ,\non\\
\chih^{2q}_1-\chih^{2q+2}_1 &=& \nat\dd \etah^{2q+1}_0 + \etah^{2q+1}_1(b+B)\ .\non
\eeq
The evaluation of the first equation on the $(b+B)$-cycle $\ch_*(\eh)\in \Omh^+(\Ic\Rch)$ yields (see section \ref{st}; we also omit reference to the inclusion homomorphism $\iota_*:\Ic\Rch\hookrightarrow \bigl( \begin{smallmatrix}
          \Ic\Rch & 0 \\
          0 & 0 \end{smallmatrix} \bigr) \subset \Mch^s_+$)
$$
\ch^{2q}_0(\eh)-\ch^{2q+2}_0(\eh) = \bb \big(\etah^{2q+1}_1\ch_{2q+2}(\eh)\big)\ ,
$$
with $\ch^{2q}_0(\eh)=\chih^{2q}_0\ch_{2q}(\eh)$ by Lemma \ref{lchcs}. Therefore, we guess that the map $MK^{\Ic}_{n}(\Ac)^q\to MK^{\Ic}_{n}(\Ac)^{q+1}$ should send a pair $(\eh,\te)^q$ to the pair $(\eh,\te')^{q+1}$ with $\te'=\te -\etah^{2q+1}_1\ch_{2q+2}(\eh)$. Indeed, one has the correct transgression relation
$$
\ch^{2q+2}_0(\eh) = \ch^{2q}_0(\eh) - \bb(\etah^{2q+1}_1\ch_{2q+2}(\eh)) = \bb\te'
$$
in the complex $X_{n-1}(\Rc,\Jc)$. Moreover, this map is well-defined at the level of equivalence classes: let $(\eh_0,\te_0)^q$ and $(\eh_1,\te_1)^q$ be two equivalent pairs. Then there exists an interpolating idempotent $\eh\in M_2(\Kc\hotimes\Ic\Rch[0,1])^+$ and a chain $\la$ such that $\te_1-\te_0=\cs^{2q}_1(\eh)+\nat\dd\la$. Hence the respective images $(\eh_0,\te_0')^{q+1}$ and $(\eh_1,\te_1')^{q+1}$ verify
$$
\te_1'-\te_0' = \te_1-\te_0 -\etah^{2q+1}_1 \, (\ch_{2q+2}(\eh_1)-\ch_{2q+2}(\eh_0)) \ .
$$
But we know the transgression relation (\ref{tran})
$$
\ch_{2q+2}(\eh_1)-\ch_{2q+2}(\eh_0)=B\cs_{2q+1}(\eh)+b\cs_{2q+3}(\eh)\ ,
$$
so that, using the identities $\chih^{2q}_1\cs_{2q+1}(\eh)=(\nat\dd\etah^{2q+1}_0+\etah^{2q+1}_1B)\cs_{2q+1}(\eh)$ and $-\chih^{2q+2}_1\cs_{2q+3}(\eh)=\etah^{2q+1}_1b \cs_{2q+3}(\eh)$ one gets (recall $\cs_1^{2q}(\eh)=\chih_1^{2q}\cs_{2q+1}(\eh)$)
\beq
\te_1'-\te_0' &=& \cs^{2q}_1(\eh)+\nat\dd\la\non\\
&& - \chih_1^{2q}\cs_{2q+1}(\eh) + \chih_1^{2q+2}\cs_{2q+3}(\eh) + \nat\dd \big(\etah^{2q+1}_0\cs_{2q+1}(\eh)\big) \non\\
&=& \cs^{2q+2}_1(\eh)+\nat\dd\big(\la+\etah^{2q+1}_0\cs_{2q+1}(\eh)\big)\ .\non
\eeq 
Hence $(\eh_0,\te_0')^{q+1}$ and $(\eh_1,\te_1')^{q+1}$ are equivalent in $MK^{\Ic}_{n}(\Ac)^{q+1}$. It remains to show that the map $MK^{\Ic}_{n}(\Ac)^q\to MK^{\Ic}_{n}(\Ac)^{q+1}$ is an isomorphism. Consider the following diagram:
$$
\xymatrix{
 HC_{n-1}(\Ac) \ar[r]^{\delta} \ar@{=}[d] & MK^{\Ic}_n(\Ac)^q  \ar[r] \ar[d] & \Kt_0(\Ic\hotimes\Ac)  \ar@{=}[d] \\
 HC_{n-1}(\Ac) \ar[r]^{\delta} & MK^{\Ic}_n(\Ac)^{q+1} \ar[r] & \Kt_0(\Ic\hotimes\Ac) }
$$
For any cyclic homology class $[\te]\in HC_{n-1}(\Ac)$ represented by a closed chain $\te$, one has $\delta(\te)=(p_0, \sqrt{2\pi i}\, \te)^q$ in $MK^{\Ic}_{n}(\Ac)^q$. But observe that $\etah^{2q+1}_1\ch_{2q+2}(p_0)=0$, so that $(p_0, \sqrt{2\pi i}\, \te)^q$ is mapped to $(p_0, \sqrt{2\pi i}\, \te)^{q+1}$ in $MK^{\Ic}_{n}(\Ac)^{q+1}$. Hence the left square is commutative. Moreover the right square is obviously commutative. The isomorphism $MK^{\Ic}_n(\Ac)^q\cong MK^{\Ic}_n(\Ac)^{q+1}$ then follows from the five-lemma.\\

\noindent The negative Chern character $\ch_{n}: MK^{\Ic}_{n}(\Ac)^q\to HN_{n}(\Ac)$ is also independent of $q$, and compatible with the $SBI$ exact sequence. Indeed if $(\eh,\te)^q$ is a representative of a class in $MK^{\Ic}_{n}(\Ac)^q$, one has by definition
$$
\ch_{n}(\eh,\te)^q = \ch^{2q}_0(\eh) - \bb\tilde{\te}\ ,
$$
where $\tilde{\te}\in \Xh(\Rc,\Jc)$ is an arbitrary lifting of $\te$. First remark that $\ch_{n}$ is well-defined at the level of equivalence classes: if $(\eh_0,\te_0)^q$ and $(\eh_1,\te_1)^q$ are equivalent, there exists an interpolation $\eh\in M_2(\Kc\hotimes\Ic\Rch[0,1])^+$ and a chain $\la$ such that $\te_1-\te_0=\cs^{2q}_1(\eh)+ \nat\dd\la$ in $X_{n-1}(\Rc,\Jc)$. Let $\tilde\te_0$, $\tilde\te_1$ and $\tilde\la$ be arbitrary liftings; then there exists a chain $\mu \in F^{n-1}\Xh(\Rc,\Jc)$ such that $\tilde\te_1-\tilde\te_0=\cs^{2q}_1(\eh)+\nat\dd\tilde\la +\mu$ in $\Xh(\Rc,\Jc)$. Hence the difference
\beq
\lefteqn{\ch_{n}(\eh_1,\te_1)^q - \ch_{n}(\eh_0,\te_0)^q = \ch^{2q}_0(\eh_1) - \ch^{2q}_0(\eh_0) - \bb(\tilde\te_1-\tilde\te_0)} \non\\
&& \qquad\qquad\qquad\qquad  = \bb\cs^{2q}_1(\eh)- \bb(\cs^{2q}_1(\eh)+\nat\dd\tilde\la +\mu)  = -\bb\mu\non
\eeq
is a coboundary of the subcomplex $F^{n-1}\Xh(\Rc,\Jc)$, and the Chern character $\ch_{n}: MK^{\Ic}_{n}(\Ac)^q\to HN_{n}(\Ac)$ is well-defined. Now if $(\eh,\te)^{q}\in MK^{\Ic}_{n}(\Ac)^q$ is any class, its image in $MK^{\Ic}_{n}(\Ac)^{q+1}$ is represented by $(\eh,\te')^{q+1}$ with $\te'=\te - \etah^{2q+1}_1\ch_{2q+2}(\eh)$.  One has
\beq
\ch_{n}(\eh,\te')^{q+1} &=& \ch^{2q+2}_0(\eh) - \bb(\tilde\te - \etah^{2q+1}_1\ch_{2q+2}(\eh)) \non\\
&=& \ch^{2q+2}_0(\eh) - \bb(\tilde\te) + \chih^{2q}_0\ch_{2q}(\eh) - \chih^{2q+2}_0\ch_{2q+2}(\eh) \non\\
&=&  \ch^{2q}_0(\eh) - \bb(\tilde\te) \ =\ \ch_{n}(\eh,\te)^q\ ,\non
\eeq
and the negative Chern character does not depend on the degree $q$. Finally, for any cyclic homology class $[\te]\in HC_{n-1}(\Ac)$, one has
$$
\ch_{n}(\delta(\te)) = \ch_{n}(p_0, \sqrt{2\pi i}\, \te) = -\sqrt{2\pi i}\, \bb(\tilde\te)\ ,
$$
which shows the commutativity of the middle square (\ref{mul}). The compatibility between the negative Chern character and the periodic Chern character on topological $K$-theory is obvious, whence the commutativity of the right square (\ref{mul}). \\

\noindent Concerning the independence of $MK^{\Ic}_n(\Ac)$ with respect to the choice of quasi-free extension $\Rc$, it suffices to consider the universal extension $0\to J\Ac \to T\Ac \to \Ac \to 0$ together with the classifying homomorphisms $T\Ac\to\Rc$ and $J\Ac\to\Jc$. The various Chern characters and Chern-Simons forms constructed in $\Xh(\Rc,\Jc)$ are obtained from the universal ones in $\Xh(T\Ac,J\Ac)$ by applying the chain map $\Xh(T\Ac,J\Ac)\to \Xh(\Rc,\Jc)$, which we know is a homotopy equivalence compatible with the adic filtrations. Once again the conclusion follows from the five-lemma.\\
The case of odd degree $n=2k+1$ is established along the same lines, replacing idempotents by invertibles. \cqfd\\

Before ending this section we need to establish the invariance of topological and multiplicative $K$-theory with respect to adjoint actions of multipliers on the $p$-summable Fr\'echet $m$-algebra $\Ic$. We say that $U$ is a multiplier if it defines continuous linear maps (left and right multiplications) $x\mapsto Ux$ and $x\mapsto xU$ on $\Ic$, which commute and fulfill 
\beq
&i)& U(xy)=(Ux)y\ ,\ (xU)y=x(Uy)\ ,\ (xy)U=x(yU)\qquad \forall x,y\in \Ic\ ,\non\\
&ii)& \Tr([U,\Ic^p])=0 \ . \non
\eeq 
$U$ is invertible if there exists a multiplier $U^{-1}$ such that the compositions $UU^{-1}$ and $U^{-1}U$ induce the identity on $\Ic$, while left and right multiplications by $U$ and $U^{-1}$ commute. In this case the adjoint action of $U$ defined by $\ad U (x)= U^{-1}x U$ is a continuous automorphism of $\Ic$ preserving the trace on $\Ic^p$. If $\Ac$ is any Fr\'echet $m$-algebra, the adjoint action of $U$ extends to the tensor product $\Kc\hotimes\Ic\Ac$ by acting trivially on the factors $\Kc$ and $\Ac$, thus defines an automorphism of $\Kt_n(\Ic\hotimes\Ac)$. Similarly if $0\to \Jc\to \Rc \to \Ac \to 0$ is a quasi-free extension, and $(\eh,\te)$ (resp. $(\gh,\te)$) represents a multiplicative $K$-theory class of even (resp. odd) degree, the adjoint action of $U$ extends to an automorphism of the pro-algebra $\Kc\hotimes\Ic\Rch$ and define maps
\be
\ad U:\quad (\eh,\te)\mapsto (U^{-1}\eh U,\te)\ ,  \quad (\gh,\te)\mapsto (U^{-1}\gh U,\te)\ .
\ee
The images represent multiplicative $K$-theory classes because the invariance of the trace implies $\ch^{2q}_0(U^{-1}\eh U)=\ch^{2q}_0(\eh)=\bb\te$ and $\ch^{2q}_0(U^{-1}\gh U)=\ch^{2q}_0(\gh)=\nat\dd\te$. The adjoint action is actually well-defined at the level of $K$-theory:

\begin{lemma}\label{linv}
Let $U$ be an invertible multiplier of $\Ic$. Then the adjoint action $\ad U$ induces the identity on $\Kt_n(\Ic\hotimes\Ac)$ and $MK^{\Ic}_n(\Ac)$. 
\end{lemma}
{\it Proof:} First we show that an idempotent $e\in M_2(\Kc\hotimes\Ic\Ac)^+$, with $e-p_0\in M_2(\Kc\hotimes\Ic\Ac)$, is smoothly homotopic to its adjoint $U^{-1}eU$. Introduce the idempotent $f_0=\bigl( \begin{smallmatrix} e & 0 \\ 0 & p_0 \end{smallmatrix} \bigr) \in M_4(\Kc\hotimes\Ic\Ac)^+$, and choose a smooth real-valued function $\xi\in\cinf[0,1]$ such that $\xi(0)=0$ and $\xi(1)=\pi/2$. We define a path $W$ of invertible multipliers of $M_4(\Kc\hotimes\Ic\Ac)$ by means of the formula
$$
W=R^{-1}\left( \begin{matrix}
1 & 0 \\ 
0 & U \end{matrix} \right)R \ ,\qquad R=\left( \begin{matrix}
\cos\xi & \sin\xi \\ 
-\sin\xi & \cos\xi \end{matrix} \right)\ ,
$$
where each entry should be viewed as a $2\times 2$ block matrix. Hence, $W$ commutes with the matrix $\pt_0=\bigl( \begin{smallmatrix} p_0 & 0 \\ 0 & p_0 \end{smallmatrix} \bigr)$. The smooth path of idempotents $f=W^{-1}f_0W$ thus provides an interpolation between $f_0$ and $f_1=\bigl( \begin{smallmatrix} U^{-1}eU & 0 \\ 0 & p_0 \end{smallmatrix} \bigr)$. Put in another way, $cfc$ interpolates the $K$-theoretic sums $e+p_0$ and $U^{-1}eU+p_0$. This shows that $e$ and $U^{-1}eU$ define the same topological $K$-theory class. \\
Now suppose that $(\eh,\te)\in MK^{\Ic}_n(\Ac)$ represents a multiplicative $K$-theory class, with $\eh\in M_2(\Kc\hotimes\Ic\Rch)^+$, $\te\in X_{n-1}(\Rc,\Jc)$ and $\ch^{2q}_0(\eh)=\bb\te$. We define as before $\fh_0=\bigl( \begin{smallmatrix} \eh & 0 \\ 0 & p_0 \end{smallmatrix} \bigr)$, and $\fh=W^{-1}\fh_0W$ provides an interpolation between $\fh_0$ and $\fh_1=\bigl( \begin{smallmatrix} U^{-1}\eh U & 0 \\ 0 & p_0 \end{smallmatrix} \bigr)$. If $s: \cinf[0,1]\to \Om^1[0,1]$ denotes the differential over $[0,1]$ and $\dd: \Rch\to \Om^1\Rch$ the noncommutative differential, the Chern-Simons form (\ref{hcs}) associated to $c\fh c$ reads
$$
\cs_1^{2q}(c\fh c) = \int_0^1\Tr\nat(-2\fh+1)\sum_{i=0}^q(\fh-\pt_0)^{2i}s\fh (\fh-\pt_0)^{2(q-i)}\dd\fh\ .
$$
One has $\dd\fh=W^{-1}\dd\fh_0W$ and $s\fh=W^{-1}(-sWW^{-1}\fh_0+\fh_0sWW^{-1})W$, hence
\beq
\lefteqn{\Tr\nat(-2\fh+1)\sum_{i=0}^q(\fh-\pt_0)^{2i}s\fh (\fh-\pt_0)^{2(q-i)}\dd\fh} \non\\
&=& - \Tr\nat(-2\fh_0+1)\sum_{i=0}^q(\fh_0-\pt_0)^{2i}sWW^{-1}\fh_0 (\fh_0-\pt_0)^{2(q-i)}\dd\fh_0 \non\\
&& +  \Tr\nat(-2\fh_0+1)\sum_{i=0}^q(\fh_0-\pt_0)^{2i}\fh_0 sWW^{-1} (\fh_0-\pt_0)^{2(q-i)}\dd\fh_0\non
\eeq
Observe that $\Tr\nat$ is a trace. In the first term of the r.h.s., we can use the identity $\fh_0 (\fh_0-\pt_0)^{2(q-i)}=(\fh_0-\pt_0)^{2(q-i)}\fh_0$ which holds for any two idempotents $\fh_0$ and $\pt_0$, and then $\fh_0\dd\fh_0 (-2\fh_0+1)=\fh_0\dd\fh_0$. In the second term of the r.h.s., we simply write $(-2\fh_0+1)(\fh_0-\pt_0)^{2i}\fh_0=(-2\fh_0+1)\fh_0(\fh_0-\pt_0)^{2i}=-\fh_0(\fh_0-\pt_0)^{2i}$. Hence we arrive at
\beq
\lefteqn{\Tr\nat(-2\fh+1)\sum_{i=0}^q(\fh-\pt_0)^{2i}s\fh (\fh-\pt_0)^{2(q-i)}\dd\fh} \non\\
&=& - \sum_{i=0}^q \Tr\nat\, (\fh_0-\pt_0)^{2i}sWW^{-1} (\fh_0-\pt_0)^{2(q-i)}\fh_0 \dd\fh_0 \non\\
&& - \sum_{i=0}^q\Tr\nat\, (\fh_0-\pt_0)^{2i} sWW^{-1} (\fh_0-\pt_0)^{2(q-i)}\dd\fh_0\, \fh_0\non \\
&=& - \sum_{i=0}^q\Tr\nat\, sWW^{-1} (\fh_0-\pt_0)^{2i} \dd\fh_0 (\fh_0-\pt_0)^{2(q-i)} \non
\eeq
by the idempotent property $\fh_0\dd\fh_0+\dd\fh_0\, \fh_0=\dd\fh_0$. It remains to show that the latter sum is a boundary:
$$
-\sum_{i=0}^q\Tr\nat\, sWW^{-1} (\fh_0-\pt_0)^{2i} \dd\fh_0 (\fh_0-\pt_0)^{2(q-i)} = \nat\dd \big( \Tr\, sWW^{-1}(\fh_0-\pt_0)^{2q+1}\big)\ .
$$
Indeed $\dd$ anticommutes with $sWW^{-1}$, and $\dd(\fh_0-\pt_0)=\dd\fh_0$. Hence if we can show that the terms $\Tr\nat\, sWW^{-1} (\fh_0-\pt_0)^{2i+1} \dd\fh_0 (\fh_0-\pt_0)^{2j+1}$ vanish, the conclusion follows. Since $\fh_0\dd\fh_0\fh_0=0$, we can write
\beq
\lefteqn{\Tr\nat\, sWW^{-1} (\fh_0-\pt_0)^{2i+1} \dd\fh_0 (\fh_0-\pt_0)^{2j+1}}\non\\
&=& \Tr\nat\, sWW^{-1} (\fh_0-\pt_0)^{2i}( -\fh_0\dd\fh_0 \pt_0 -\pt_0\dd\fh_0\fh_0 + \pt_0\dd\fh_0\pt_0)(\fh_0-\pt_0)^{2j}\non\\
&=& \Tr\nat\, sWW^{-1} (\fh_0-\pt_0)^{2i}( -\fh_0\dd\fh_0 \pt_0 -\dd\fh_0\fh_0\pt_0 + \dd\fh_0\pt_0)(\fh_0-\pt_0)^{2j}\non\\
&=& 0 \ ,\non
\eeq
where we used the fact that $\pt_0$ commutes with $sWW^{-1}$ and the even powers of $\fh_0-\pt_0$. Hence $\cs^{2q}_1(c\fh c)\equiv 0 \mod \nat\dd$, which shows that the pairs $(\eh,\te)$ and $(\Uh^{-1}\eh U,\te)$ are equivalent. The adjoint action of $U$ on multiplicative $K$-theory groups in even degrees is thus the identity.\\
One proceeds in the same fashion with odd groups. Let $g\in (\Kc\hotimes\Ic\Ac)^+$ be an invertible such that $g-1\in \Kc\hotimes\Ic\Ac$. Introduce $u_0=\bigl( \begin{smallmatrix} g & 0 \\ 0 & 1 \end{smallmatrix} \bigr)$ and the invertible path $u=W^{-1}u_0W$, where $W=R^{-1}\bigl( \begin{smallmatrix} 1 & 0 \\ 0 & U \end{smallmatrix} \bigr)R$ is now viewed as a path of invertible multipliers of $M_2(\Kc\hotimes\Ic\Ac)$. Hence $u$ interpolates between $u_0$ and $u_1=\bigl( \begin{smallmatrix} U^{-1}gU & 0 \\ 0 & 1 \end{smallmatrix} \bigr)$. This shows that $g$ and $U^{-1}gU$ define the same topological $K$-theory class.\\
Now suppose that $(\gh,\te)\in MK^{\Ic}_n(\Ac)$ represents a multiplicative $K$-theory class, with $\gh\in (\Kc\hotimes\Ic\Rch)^+$, $\te\in X_{n-1}(\Rc,\Jc)$ and $\ch^{2q}_1(\gh)=\nat\dd\te$. We define $\uh_0=\bigl( \begin{smallmatrix} \gh & 0 \\ 0 & 1 \end{smallmatrix} \bigr)$, and $\uh=W^{-1}\uh_0W$ provides an interpolation between $\uh_0$ and $\uh_1=\bigl( \begin{smallmatrix} U^{-1}\gh U & 0 \\ 0 & 1 \end{smallmatrix} \bigr)$. The Chern-Simons form (\ref{hcs}) associated to $\uh$ reads
$$
\cs^{2q}_0(\uh) = \frac{1}{\sqrt{2\pi i}} \frac{(q!)^2}{(2q)!}\int_0^1 \Tr\, \uh^{-1}[(\uh-1)(\uh^{-1}-1)]^q s\uh\ .
$$ 
Using $s\uh= W^{-1}(-sWW^{-1}\uh_0+\uh_0sWW^{-1})W$, one gets
\beq
\Tr\, \uh^{-1}[(\uh-1)(\uh^{-1}-1)]^q s\uh &=& -\Tr\, \uh_0^{-1}[(\uh_0-1)(\uh_0^{-1}-1)]^q sWW^{-1}\uh_0 \non\\
&& + \Tr\, \uh_0^{-1}[(\uh_0-1)(\uh_0^{-1}-1)]^q \uh_0sWW^{-1}\non\\
&\equiv& 0 \mod \bb\non
\eeq
Hence $\cs^{2q}_0(\uh) \equiv 0\mod\bb$ and the pair $(\Uh^{-1}\gh U,\te)$ is equivalent to $(\gh,\te)$. The adjoint action of $U$ on the odd multiplicative $K$-theory groups is therefore the identity. \cqfd\\

\begin{example}\label{ec}\textup{Take $\Ac=\cc$ and $\Ic=\Lc^p(H)$ a Schatten ideal. It is known that $\Kt_0(\Ic)=\zz$ and $\Kt_1(\Ic)=0$. Furthermore $HC_n(\cc)=\cc$ for $n=2k\geq 0$ and vanishes otherwise. Hence the exact sequence yields
$$
MK_n^{\Ic}(\cc)=\left\{ \begin{array}{ll}
\zz & n\leq 0\ \mbox{even} \\
\cc^{\times} & n>0\ \mbox{odd} \\
0 & \mbox{otherwise} \end{array} \right.
$$ 
The multiplicative $K$-theory of $\cc$ is the natural target for index maps in even degree, and for regulator maps in odd degree (see  \cite{CK} and Example \ref{ereg}). }
\end{example}

Multiplicative $K$-theory has close connections with higher algebraic $K$-theory \cite{K1,R}. In fact there exists a morphism $\Ka_n(\Ac)\to MK^{\Ic}_{n}(\Ac)$ in all degrees, and composition with the negative Chern character coincides with the Jones-Goodwillie map $\Ka_n(\Ac)\to HN_n(\Ac)$ \cite{J}. See \cite{CT} for an exact sequence relating topological and algebraic $K$-theories of locally convex algebras stabilized by operator ideals.

\section{Riemann-Roch-Grothendieck theorem}\label{srr}

In this section we construct direct images of topological and multiplicative $K$-theory under quasihomomorphisms and show their compatibility with the $K$-theory and cyclic homology exact sequences. This provides a noncommutative version of the Riemann-Roch-Grothendieck theorem. If $\Ic$ is a $p$-summable Fr\'echet $m$-algebra, with trace $\Tr:\Ic^p\to \cc$, the tensor product $\Ic\hotimes\Ic$ is in a natural way a $p$-sumable algebra with trace $\Tr\hotimes\Tr$. We demand that $\Ic$ is provided with an external product as follows.
\begin{definition}
A $p$-summable Fr\'echet $m$-algebra $\Ic$ is \emph{multiplicative} if there exists a continuous algebra homomorphism (external product)
$$
\boxtimes:\Ic\hotimes\Ic\to \Ic
$$
such that the composition $\Tr\circ\boxtimes$ coincides with the trace $\Tr\hotimes\Tr$ on $(\Ic\hotimes\Ic)^p$. Two external products $\boxtimes$ and $\boxtimes'$ are \emph{equivalent} if there exists an invertible multiplier $U$ of $\Ic$ such that $\boxtimes'=\ad U\circ\boxtimes$ on $\Ic\hotimes\Ic$. 
\end{definition}

Hence if $\Ic$ is multiplicative the homomorphism $\boxtimes$ induces additive maps $\Kt_n(\Ic\hotimes\Ic\hotimes\Ac)\to \Kt_n(\Ic\hotimes\Ac)$ and $MK^{\Ic\hotimes\Ic}_n(\Ac)\to MK^{\Ic}_n(\Ac)$, clearly compatible with the commutative diagram of Proposition \ref{pmul}. Moreover two equivalent products induce the same maps, the adjoint action $\ad U$ being trivial on $K$-theory by Lemma \ref{linv}. In practice the algebra $\Ic$ often arises with external products defined only modulo equivalence:

\begin{example}\textup{Let $\Ic=\Lc^p(H)$ be the Schatten ideal of $p$-summable operators on a separable infinite-dimensional Hilbert space $H$, provided with the operator trace. The tensor product $\Lc^p(H)\hotimes \Lc^p(H)$ is naturally mapped to the algebra $\Lc^p(H\otimes H)$, and choosing an isomorphism of Hilbert spaces $H\otimes H \cong H$ allows to identify $\Lc^p(H\otimes H)$ with $\Lc^p(H)$ modulo the adjoint action of unitary operators $U\in \Lc(H)$. The product $\boxtimes:\Lc^p(H)\hotimes\Lc^p(H)\to\Lc^p(H)$ is thus compatible with the traces, and canonically defined modulo the adjoint action of unitary operators. }
\end{example}

Let $\Ac$ and $\Bc$ be any Fr\'echet $m$-algebras. Let $\rho:\Ac\to\Ec^s\triangleright \Ic^s\hotimes\Bc$ be a quasihomomorphism of parity $p$ mod $2$, and suppose that $\Ic$ is finitely summable (the exact summability degree is irrelevant for the moment). We want to show that $\rho$ induces an additive map
\be
\rho_!:\Kt_n(\Ic\hotimes \Ac)\to \Kt_{n-p}(\Ic\hotimes\Bc)\qquad \forall n\in\zz
\ee
provided $\Ic$ is multiplicative. This is has nothing to do with cyclic homology and we don't need to assume $\Ec$ admissible. Thanks to Bott periodicity, it is sufficient to define $\rho_!$ on $\Kt_1(\Ic\hotimes\Ac)$, where it is given by very explicit formulas. Suppose first that $p$ is even. Then $\rho$ is described by a pair of homomorphisms $(\rho_+,\rho_-):\Ac\rightrightarrows \Ec$ which coincide modulo $\Ic\hotimes\Bc$. For any invertible element $g\in (\Kc\hotimes\Ic\hotimes\Ac)^+$ with $g-1\in \Kc\hotimes\Ic\hotimes\Ac$ representing a $K$-theory class in $\Kt_1(\Ic\hotimes\Ac)$, one has $\rho_{\pm}(g) \in (\Kc\hotimes\Ic\hotimes\Ec)^+$ and $\rho_+(g)-\rho_-(g)\in \Kc\hotimes\Ic\hotimes\Ic\hotimes\Bc$, where the homomorphisms $\rho_+$ and $\rho_-$ are extended to the unitalized algebra $(\Kc\hotimes\Ic\hotimes\Ac)^+$ by acting trivially on the factor $\Kc\hotimes\Ic$ and preserving the unit. set
\be
\rho_!(g) = \rho_+(g)\rho_-(g)^{-1}\ \in (\Kc\hotimes\Ic\hotimes\Ic\hotimes\Bc)^+\ .
\ee
Using the homomorphism $\boxtimes:\Ic\hotimes\Ic\to\Ic$, we may therefore consider $\rho_!(g)$ as an invertible element of $(\Kc\hotimes\Ic\hotimes\Bc)^+$ such that $\rho_!(g)-1\in \Kc\hotimes\Ic\hotimes\Bc$. It is clear that the homotopy class of $\rho_!(g)$ only depends on the homotopy class of $g$, hence the map $\rho_!:\Kt_1(\Ic\hotimes\Ac)\to \Kt_1(\Ic\hotimes\Bc)$ is well-defined.\\
When $p$ is odd, $\rho$ is a homomorphism $\Ac\to M_2(\Ec)$ such that the off-diagonal terms lie in $\Ic\hotimes\Bc$. For any invertible element $g\in (\Kc\hotimes\Ic\hotimes\Ac)^+$ as above, one has $\rho(g)\in M_2(\Kc\hotimes\Ic\hotimes\Ec)^+$ with off-diagonal elements in $\Kc\hotimes\Ic\hotimes\Ic\hotimes\Bc$. Set
\be
\rho_!(g)=\rho(g)^{-1}p_0 \rho(g)\in M_2(\Kc\hotimes\Ic\hotimes\Ic\hotimes\Bc)^+\ ,
\ee
where $p_0=\bigl(\begin{smallmatrix} 1 & 0 \\ 0 & 0 \end{smallmatrix}\bigr)$ is the trivial matrix idempotent. Again applying the external product $\boxtimes$ we may consider $\rho_!(g)$ as an idempotent of $M_2(\Kc\hotimes\Ic\hotimes\Bc)^+$ such that $\rho_!(g)-p_0\in M_2(\Kc\hotimes\Ic\hotimes\Bc)$. The homotopy class only depends on the homotopy class of $g$ and we thus obtain $\rho_!: \Kt_1(\Ic\hotimes\Ac)\to \Kt_0(\Ic\hotimes\Bc)$.\\
To define $\rho_!$ on $\Kt_0(\Ic\hotimes\Ac)$ it suffices to pass to the suspensions $S\Ac = \Ac(0,1)$ and $S\Bc=\Bc(0,1)$, then apply the pushforward map constructed above $\rho_! : \Kt_1(\Ic\hotimes{S}\Ac) \to \Kt_{1-p}(\Ic\hotimes{S}\Bc)$ with trivial action on the factor $\cinf(0,1)$. The Bott isomorphisms $\Kt_n(\Ic\hotimes\cdot)\cong \Kt_{n+1}(\Ic\hotimes{S}\cdot)$ allow to define $\rho_!$ for the original algebras through a \emph{graded}-commutative diagram
\be
\vcenter{\xymatrix{
\Kt_0(\Ic\hotimes\Ac) \ar[r]^{\sim} \ar[d]_{\rho_!} & \Kt_1(\Ic\hotimes{S}\Ac) \ar[d]^{\rho_!} \\
\Kt_{-p}(\Ic\hotimes\Bc) \ar[r]^{\sim}  & \Kt_{1-p}(\Ic\hotimes{S}\Bc) }}
\ee 
Note the following subtlety concerning graduations: since $\Kt_n$ has parity $n$ mod 2 by definition, the Bott isomorphisms are \emph{odd}. As a consequence, when $p$ is also odd, the above diagram must be \emph{anti}-commutative. These conventions are necessary if we want to avoid sign problems with the theorem below. \\

Now choose a quasi-free extension $0\to \Jc \to \Rc \to \Bc \to 0$ for $\Bc$, and suppose that the algebra $\Ec\triangleright \Ic\hotimes\Bc$ is $\Rc$-admissible. We impose the following compatibility between the parity of the quasihomomorphism $\rho$ and the summability degree of $\Ic$: in the even case $\Ic$ is $(p+1)$-summable with $p$ even, and in the odd case $\Ic$ is $p$-summable with $p$ odd (this complicated choice is dictated by the theorem below). In both cases the bivariant Chern character $\ch^p(\rho)\in HC^p(\Ac,\Bc)$ constructed in section \ref{sbiv} induces a map
\be
\ch^p(\rho): HC_n(\Ac)\to HC_{n-p}(\Bc)\qquad \forall n\in \zz\ .
\ee
Combining $\rho_!$ with the bivariant Chern character yields a transformation in multiplicative $K$-theory, compatible with the diagram (\ref{mul}) of Proposition \ref{pmul}. This will be detailed in the proof of the following noncommutative version of the Riemann-Roch-Grothendieck theorem:

\begin{theorem}\label{trr}
Let $\Ac$, $\Bc$ be Fr\'echet $m$-algebras, and choose a quasi-free extension $0 \to \Jc \to \Rc \to \Bc \to 0$. Let $\rho:\Ac\to \Ec^s\triangleright \Ic^s\hotimes\Bc$ be a quasihomomorphism of parity $p \mod 2$, where $\Ic$ is multiplicative and $(p+1)$-summable in the even case, $p$-summable in the odd case. Suppose that $\Ec\triangleright \Ic\hotimes\Bc$ is $\Rc$-admissible. Then $\rho$ defines a transformation in multiplicative $K$-theory $\rho_!:MK^{\Ic}_n(\Ac)\to MK^{\Ic}_{n-p}(\Bc)$ compatible with the $K$-theory exact sequences for $\Ac$ and $\Bc$, whence a graded-commutative diagram 
\be
\vcenter{\xymatrix{
\Kt_{n+1}(\Ic\hotimes\Ac) \ar[r] \ar[d]^{\rho_!} & HC_{n-1}(\Ac) \ar[r] \ar[d]^{\ch^p(\rho)} & MK^{\Ic}_n(\Ac)  \ar[r] \ar[d]^{\rho_!}  & \Kt_n(\Ic\hotimes\Ac)  \ar[d]^{\rho_!}  \\
\Kt_{n+1-p}(\Ic\hotimes\Bc) \ar[r]  & HC_{n-1-p}(\Bc) \ar[r]  & MK^{\Ic}_{n-p}(\Bc)  \ar[r]  & \Kt_{n-p}(\Ic\hotimes\Bc) }}  \label{rr}
\ee
The vertical arrows are invariant under conjugation of quasihomomorphisms; the arrow in topological $K$-theory $\Kt_n(\Ic\hotimes\Ac)\to\Kt_{n-p}(\Ic\hotimes\Bc)$ is also invariant under homotopy of quasihomomorphisms. Moreover (\ref{rr}) is compatible with the commutative diagram of Theorem \ref{tbiv} (with connecting map $B$ rescaled by a factor $-\sqrt{2\pi i}$) after taking the Chern characters $MK^{\Ic}_n\to HN_n$ and $\Kt_n(\Ic\hotimes\ .\ )\to HP_n$.
\end{theorem}
{\it Proof:} As a general rule, the bivariant cyclic cohomomology $HC^p(\Ac,\Bc)$ is described as the cohomology of the complex $\hom^p(\Xh(T\Ac,J\Ac), \Xh(\Rc,\Jc))$ of linear maps of order $\leq p$, where we choose the universal free extension $0 \to J\Ac \to T\Ac \to \Ac \to 0$ for $\Ac$ and the quasi-free extension $0 \to \Jc \to \Rc \to \Bc \to 0$ for $\Bc$. By hypothesis, the algebra $\Ec \triangleright \Ic\hotimes\Bc$ is $\Rc$-admissible (Definition \ref{dadm}), i.e. one has a commutative diagram 
$$
\vcenter{\xymatrix{
0\ar[r] & \Nc \ar[r] & \Mc \ar[r] & \Ec \ar[r] & 0 \\
0 \ar[r] & \Ic\hotimes \Jc \ar[r] \ar[u] & \Ic\hotimes \Rc \ar[r] \ar[u] & \Ic\hotimes\Bc \ar[r] \ar[u] & 0}}
$$
verifying adequate properties with respect to the trace over $\Ic$. The detailed construction of the pushforward map in multiplicative $K$-theory $\rho_!:MK^{\Ic}_n(\Ac)\to MK^{\Ic}_{n-p}(\Bc)$ depends on the respective parities of $n$ and $p$. \\

\noindent  {\bf i) $n=2k+1$ is odd and $p=2q$ is even.} Our first task is to understand the composition of the topological Chern character $\ch^p_1:\Kt_1(\Ic\hotimes\Ac)\to HP_1(\Ac)$ with the bivariant Chern character $\ch^p(\rho) \in HC^p(\Ac,\Bc)$. For notational simplicity, we write as usual $\Ic\Ac$ for the tensor product $\Ic\hotimes\Ac$. Without loss of generality, we may suppose that an element of $\Kt_1(\Ic\Ac)$ is represented by an invertible $g\in (\Ic\Ac)^+$ such that $g-1\in \Ic\Ac$ (indeed the algebra $\Kc$ of smooth compact operators plays a trivial role in what follows). Since the universal free extension $T\Ac$ is chosen, we can take the \emph{canonical lift} $\gh\in (\Ic\Th\Ac)^+$ which corresponds to the image of $g$ under the canonical linear inclusion $\Ac\hookrightarrow \Omh^+\Ac\cong\Th\Ac$ as the the subspace of zero-forms. Its inverse is given by the series (\ref{glif}), with $\Kc$ replaced by $\Ic$. The $p$-th higher Chern character of $\gh$ is then represented by the cycle (\ref{hch})
$$
\ch^p_1(\gh) = \frac{1}{\sqrt{2\pi i}}\, \frac{(q!)^2}{p!}\, \Tr\nat \gh^{-1}[(\gh-1)(\gh^{-1}-1)]^q\dd\gh\ \in \Om^1\Th\Ac_{\nat}\ ,
$$
Observe that $(p+1)$ powers of $\Ic$ appear in the products of $(\gh-1)$, $(\gh^{-1}-1)$ and $\dd\gh$ hence the trace $\Tr:\Ic^{p+1}\to \cc$ is well-defined. On the other hand, the bivariant Chern character of the quasihomomorphism $\rho$ (section \ref{sbiv}) is represented by the composition of chain maps $\ch^p(\rho)=\chih^p\rho_*\gamma : X(\Th\Ac)\to \Omh\Th\Ac \to \Omh\Mch^s_+ \to X(\Rch)$, hence the composite $\ch^p(\rho)\cdot\ch^p_1(\gh)$ requires to compute first the image of $\ch^p_1(\gh)$ under the Goodwillie equivalence $\gamma:X(\Th\Ac)\to \Omh \Th\Ac$. This tricky computation can be simplified as follows. We use the isomorphism $\Ic\Ac\cong \cc\hotimes\Ic\Ac$ to identify $\gh$ with the invertible element
$$
\uh = 1+ e\otimes(\gh-1)  \ \in (\cc\hotimes\Ic\Ac)^+\hookrightarrow (\Th\cc\hotimes \Ic\Th\Ac)^+\ ,
$$
where $e$ is the unit of $\cc$. As usual we regard $\cc\hotimes\Ic\Ac$ as the subspace of zero-forms of the algebra $\Th\cc\hotimes \Ic\Th\Ac\cong \Omh^+\cc \hotimes\Ic\Omh^+\Ac$. It is not hard to calculate that the inverse of $\uh$ is given by the series
$$
\uh^{-1}=\sum_{i\geq 0}\big( (dede)^i\otimes[(\gh^{-1}-1)(\gh-1)]^i+e(dede)^i\otimes [(\gh^{-1}-1)(\gh-1)]^i(\gh^{-1}-1)\big)
$$
with the convention $(dede)^0=1$. Observe that the power of $\Ic$ is equal to the power of $e$ in each term of this series. Also, recall that the canonical lift of $e$ is the idempotent
$$
\eh=e+\sum_{i\geq 1}\frac{(2i)!}{(i!)^2} (e-\frac{1}{2})(dede)^i\ \in \Th\cc\ .
$$
We define the \emph{fundamental class} of degree $p=2q$ as the trace $[2q]:\Th\cc\to\cc$ vanishing on all the differential forms $e(dede)^i$ and $(dede)^i$ except $e(dede)^q$, and normalized so that $[2q]\, \eh=1$. One thus have
$$
[2q]\, e(dede)^q = \frac{(q!)^2}{p!}\ ,\qquad [2q]\, (\mbox{anything else})=0\ .
$$
Of course, $[2q]$ is the generator in degree $p$ of the cyclic cohomology of $\cc$. The fact that it is a trace over $\Th\cc$ is crucial. Indeed, one finds the identity
$$
\Tr\nat [2q]\,  \uh^{-1}\dd\uh = \frac{(q!)^2}{p!}\, \Tr\nat \gh^{-1}[(\gh-1)(\gh^{-1}-1)]^q\dd\gh\ \in \Om^1\Th\Ac_{\nat}\ ,
$$
so that the Chern character $\ch^p_1(\gh)$ is exactly the cycle $\frac{1}{\sqrt{2\pi i}}\, \Tr\nat[2q]\,  \uh^{-1}\dd\uh$. This simplifies drastically the computation of $\gamma\ch^p_1(\gh)$. The Goodwillie equivalence $\gamma$ is explicitly constructed in section \ref{scy}; it is based on the linear map $\phi: \Th\Ac\to \Om^2 \Th\Ac$ verifying the properties $\phi(xy)=\phi(x) y + x\phi(y) +\dd x\dd y$ for all $x,y\in \Th\Ac$, and $\phi(a)=0$ whenever $a\in \Ac$. We extend $\phi$ to a linear map 
$$
\phi: (\Th\cc\hotimes\Ic\Th\Ac)^+\to \Th\cc\hotimes \Ic\Om^2 \Th\Ac
$$
acting by the identity on the factor $\Th\cc\hotimes\Ic$ and setting $\phi(1)=0$. This implies $\phi(\uh\uh^{-1})=0=\phi(\uh)\uh^{-1}+\uh \phi(\uh^{-1}) + \dd\uh\dd\uh^{-1}$. Moreover $\uh$ lies in $(\cc\hotimes\Ic\Ac)^+$ so that $\phi(\uh)=0$, and one gets
$$
\phi(\uh^{-1})=-\uh^{-1}\dd\uh\dd\uh^{-1}\ .
$$
Now, extending $\phi$ in all degrees as in section \ref{scy} one gets a linear map $\phi:\Th\cc\hotimes \Ic\Om^i \Th\Ac\to \Th\cc\hotimes \Ic\Om^{i+2} \Th\Ac$ for any $i\in\nn$. The following computation is then straightforward (remark that $\Tr[2q]$ is a trace hence cyclic permutations are allowed; moreover the fundamental class $[2q]$ selects $(p+1)$ powers of $e$, hence of $\Ic$, so that $\Tr$ is well-defined):
$$
\gamma (\Tr\nat[2q]\, \uh^{-1}\dd\uh) = \sum_{i\geq 0}\Tr[2q]\, \phi^i(\uh^{-1}\dd\uh) = \sum_{i\geq 0}(-)^i i! \, \Tr[2q]\, \uh^{-1}\dd\uh (\dd\uh^{-1}\dd\uh)^i\ .
$$
Hence $\gamma\ch^p_1(\gh)$ is equal to this $(b+B)$-cycle over $\Th\Ac$, divided by a factor $\sqrt{2\pi i}$. It remains to apply the chain map $\chi^p\rho_*:\Omh \Th\Ac\to X(\Rch)$ associated to the quasihomomorphism $\rho:\Ac\to \Ec^s\triangleright \Ic^s\Bc$. In $2\times 2$ matrix notation, the image of any $x\in \Th\Ac$ under the lifted quasihomomorphism $\rho_*:\Th\Ac\to \Mch^s\triangleright \Ic^s\Rch$ and the odd multiplier $F$ read 
$$
\rho_*x=\left(\begin{matrix}
x_+ & 0 \\
0 & x_- \end{matrix} \right)\in \Mch^s_+\ ,\qquad F=\left(\begin{matrix}
0 & 1 \\
1 & 0 \end{matrix} \right)\ ,
$$
and the difference $x_+-x_-$ is therefore an element of the pro-algebra $\Ic\Rch$. On the other hand, the odd component of the chain map $\chih^p\rho_*$ evaluated on a $(p+1)$-form $x_0\dd x_1\ldots\dd x_{p+1}$ is given by Eqs. (\ref{chi}):
$$
\frac{q!}{(p+1)!} \sum_{i=1}^{p+1}  \tau'\nat(\rho_*x_0[F,\rho_*x_1]\ldots\dd (\rho_*x_i) \ldots [F,\rho_*x_{p+1}])
$$
where $\tau'=\frac{1}{2}\tau(F[F,\ ])$ is the modified supertrace of even degree. Then we extend canonically $\rho_*$ to a unital homomorphism $(\Th\cc\hotimes\Ic\Th\Ac)^+\to (\Th\cc\hotimes\Ic\Mch^s_+)^+$ by taking the identity on the factor $\Th\cc\hotimes\Ic$. One thus has $\rho_*\uh=\bigl( \begin{smallmatrix} \uh_+ & 0 \\ 0 & \uh_- \end{smallmatrix} \bigr)$ with $\uh_+-\uh_-\in \Th\cc\hotimes\Ic\Ic\Rch$. A direct computation gives
\beq
\ch^p(\rho)\,(\Tr\nat[2q]\, \uh^{-1}\dd\uh) &=& (-)^q q!\, \chih^p_1\rho_* \,\Tr [2q](\uh^{-1}\dd\uh (\dd\uh^{-1}\dd\uh)^q)\non\\
&=& \frac{(q!)^2}{p!}\, \Tr\nat[2q]\, \ut^{-1}[(\ut-1)(\ut^{-1}-1)]^q\dd\ut\ ,\non
\eeq
where $\ut=\uh_+\uh_-^{-1}\in (\Th\cc\hotimes \Ic\Ic\Rch)^+$ may be considered as an invertible element of the pro-algebra $(\Th\cc\hotimes \Ic\Rch)^+$ after applying the homomorphism $\boxtimes: \Ic\hotimes\Ic\to \Ic$. Dividing by a factor $\sqrt{2\pi i}$, the right-hand-side should be defined as the Chern character $\ch^p_1(\ut)$, cf. (\ref{hch}). One thus gets the identity
$$
\ch^p(\rho)\cdot \ch^p_1(\gh) = \ch^p_1(\ut)
$$
at the level of cycles in $X(\Rch)$. Now, observe that the projection of $\ut$ onto the quotient algebra $(\cc\hotimes\Ic\Bc)^+$ is $u_+u_-^{-1}=1+e\hotimes(\rho_+(g)\rho_-(g)^{-1}-1)$. It corresponds to the direct image $\rho_+(g)\rho_-(g)^{-1}=\rho_!(g)$ by virtue of the isomorphism $(\cc\hotimes\Ic\Bc)^+\cong (\Ic\Bc)^+$. Hence we expect that $\ch^p_1(\ut)$ is homologous to the Chern character of any invertible lift $\widehat{\rho_!(g)}\in (\Ic\Rch)^+$. To see this, consider an invertible path $\vh\in (\Th\cc\hotimes\Ic\Rch[0,1])^+$ connecting homotopically $\vh_0=\ut$ to $\vh_1=1+\eh\otimes(\widehat{\rho_!(g)}-1)$, and such that its projection onto $(\cc\hotimes\Ic\Bc[0,1])^+$ is the constant invertible function $1+e\otimes(\rho_!(g)-1)$ over $[0,1]$. Such a path always exists, for example the linear interpolation
\be
\vh_t=t\big(1+\eh\otimes(\widehat{\rho_!(g)}-1)\big)+ (1-t)\ut\ ,\qquad t\in[0,1] \label{vh}
\ee
works. Since the evaluation of the fundamental class $[2q]$ on the canonical idempotent lift $\eh$ is the unit, a little computation shows the equality
$$
\ch^p_1(\vh_1)=\ch^p_1\big(\widehat{\rho_!(g)}\big)\ \in \Om^1\Rch_{\nat}
$$
at the level of cycles. Moreover, the Chern-Simons form associated to $\vh$, defined in analogy with formulas (\ref{hcs})
$$
\cs^p_0(\vh) = \frac{1}{\sqrt{2\pi i}}\, \frac{(q!)^2}{p!} \int_0^1dt\, \Tr [2q]\, \vh^{-1}[(\vh-1)(\vh^{-1}-1)]^q\frac{\d\vh}{\d t} 
$$
fulfills the transgression relation in the complex $X(\Rch)$
$$
\nat\dd \cs^p_0(\vh) = \ch^p_1(\vh_1)-\ch^p_1(\vh_0)=\ch^p_1(\widehat{\rho_!(g)})-\ch^p_1(\ut)
$$
as wanted. We are now in a position to define the map $\rho_!$ on multiplicative $K$-theory. Let a pair $(\gh,\te)$ represent a class in $MK^{\Ic}_n(\Ac)$ of odd degree $n=2k+1$. From Remark \ref{rcan}, we know that $\gh\in (\Ic\Th\Ac)^+$ can be taken as the \emph{canonical lift} of some invertible element $g\in (\Ic\Ac)^+$. Then, the transgression $\te$ is a chain of even degree in the quotient complex $X_{n-1}(T\Ac,J\Ac)$, and the relation $\ch^p_1(\gh)=\nat\dd\te$ holds in $X_{n-1}(T\Ac,J\Ac)$. We set
\be
\boxed{ \rho_!(\gh,\te) = \big(\widehat{\rho_!(g)}\ ,\ \ch^p(\rho)\cdot\te +\cs^p_0(\vh)\big)\ \in MK^{\Ic}_{n-p}(\Bc) }
\ee
where $\widehat{\rho_!(g)}\in (\Ic\Rch)^+$ is any invertible lift of $\rho_!(g)$ and $\vh$ is an invertible path constructed as above. Let us explain why this defines a multiplicative $K$-theory class. First, the bivariant Chern character $\ch^p(\rho)\in HC^p(\Ac,\Bc)$ induces a morphism of quotient complexes $X_{n-1}(T\Ac,J\Ac)\to X_{n-p-1}(\Rc,\Jc)$, hence $\ch^p(\rho)\cdot\te$ is a well-defined chain of even degree in $X_{n-p-1}(\Rc,\Jc)$. Regarding also $\cs^p_0(\vh)$ as an element of $X_{n-p-1}(\Rc,\Jc)$, we see that the relation
$$
\nat\dd(\ch^p(\rho)\cdot\te +\cs^p_0(\vh))=\ch^p(\rho)\cdot \ch^p_1(\gh)+\ch^p_1\big(\widehat{\rho_!(g)}\big) - \ch^p_1(\ut)=\ch^p_1\big(\widehat{\rho_!(g)}\big)
$$
holds in this quotient complex, hence $\rho_!(\gh,\te)$ represents a class in $MK^{\Ic}_{n-p}(\Bc)$. In fact, the latter does not depend on the choice of lifting $\widehat{\rho_!(g)}$, nor on the invertible path $\vh$. This can be proved simultaneously with the fact that the equivalence class of $\rho_!(\gh,\te)$ depends only on the equivalence class of $(\gh,\te)$. To show this, consider two equivalent pairs $(\gh_0,\te_0)$ and $(\gh_1,\te_1)$ representing the same element of $MK^{\Ic}_n(\Ac)$.  It means there exists a homotopy $\gh\in (\Ic\Th\Ac[0,1]_x)^+$ between $\gh_0$ and $\gh_1$ (we denote by $x$ the variable of this homotopy, which should not be confused with the variable $t$ used in the definition of the interpolation (\ref{vh})), and a chain $\la \in X_{n-1}(T\Ac,J\Ac)$ such that $\te_1-\te_0=\cs^p_0(\gh)+\bb\la$. From Remark \ref{rcan}, we may suppose that $\gh_0$, $\gh_1$ and $\gh$ are respectively the canonical lifts of invertibles $g_0,g_1 \in (\Ic\Ac)^+$ and $g\in (\Ic\Ac[0,1]_x)^+$. By definition one has
$$
\rho_!(\gh_i,\te_i)=\big(\widehat{\rho_!(g_i)}\ ,\ \ch^p(\rho)\cdot\te_i +\cs^p_0(\vh(g_i))\big)\ ,\quad i=0,1\ ,
$$
where $\vh(g_i)\in (\Th\cc\hotimes\Ic\Rch[0,1]_t)^+$ is a choice of invertible path associated to $g_i$, for example by Eq. (\ref{vh}). Choose an invertible path $\widehat{\rho_!(g)}\in (\Ic\Rch[0,1]_x)^+$ interpolating $\widehat{\rho_!(g_0)}$ and $\widehat{\rho_!(g_1)}$: it can be chosen as a lift of the path $\rho_!(g)=\rho_+(g)\rho_-(g)^{-1}$. Our goal is to show that the relation
$$
\ch^p(\rho)\cdot(\te_1-\te_0) +\cs^p_0(\vh(g_1))-\cs^p_0(\vh(g_0))\equiv \cs^p_0(\widehat{\rho_!(g)}) \mod \bb
$$
holds in $X_{n-p-1}(\Rc,\Jc)$. As before we identify the canonical lift $\gh$ of $g$ with the invertible element
$$
\uh = 1+ e\otimes(\gh-1)  \ \in (\cc\hotimes\Ic\Ac[0,1]_x)^+\hookrightarrow (\Th\cc\hotimes\Ic\Th\Ac[0,1]_x)^+\ ,
$$
and we remark that the higher Chern-Simons form $\cs^p_0(\gh)$ given by Lemma \ref{lchcs}, Eqs. (\ref{hcs}), can be written as
$$
\cs^p_0(\gh) = \frac{1}{\sqrt{2\pi i}} \int_0^1 \Tr[2q]\, \uh^{-1}s\uh 
$$
where $s=dx\frac{\d}{\d x}$ is the de Rham coboundary acting on the space of differential forms $\Om[0,1]_x$. It follows that the computation of $\ch^p(\rho)\cdot \cs^p_0(\gh)$ requires first to evaluate the Goodwillie equivalence $\gamma$ on the one-form $\omh=\uh^{-1}s\uh$. To this end, we extend $\phi$ to a linear map 
$$
\phi:\Th\cc\hotimes\Ic\Om^i\Th\Ac\hotimes \Om[0,1]_x\to \Th\cc\hotimes\Ic\Om^{i+2} \Th\Ac\hotimes \Om[0,1]_x
$$
acting by the identity on the factors $\Th\cc\hotimes\Ic$ and $\Om[0,1]_x$. The algebraic property of $\phi$ implies $\phi(\uh^{-1}s\uh)=\phi(\uh^{-1})s\uh + \uh^{-1}s\phi(\uh) + \dd\uh^{-1}\dd(s\uh)$. From $\phi(\uh)=0$, $\phi(\uh^{-1})=-\uh^{-1}\dd\uh\dd\uh^{-1}$ and $\dd \uh^{-1}=-\uh^{-1}\dd\uh \uh^{-1}$ we deduce
$$
\phi(\omh)=-\uh^{-1}\dd\uh\dd\omh\ .
$$
Then the image of $[2q]\,\omh$ under the Goodwillie equivalence is a straightforward computation, taking into account the tracial property of the fundamental class $\Tr[2q]$ and the fact that $\Om[0,1]_x$ is a commutative algebra:
\beq
\lefteqn{\gamma (\Tr[2q]\,\omh) = \sum_{i\geq 0}\Tr[2q]\phi^i(\omh) =}\non\\
&&\Tr[2q]\,\omh+ \sum_{i\geq 1}(-)^i(i-1)!\sum_{j=0}^{i-1}\Tr[2q]\, \uh^{-1}\dd\uh(\dd\uh^{-1}\dd\uh)^j\dd\omh (\dd\uh^{-1}\dd\uh)^{i-j-1}\ .\non
\eeq
This is a chain in the bicomplex $\Omh\Th\Ac\hotimes \Om[0,1]_x$ endowed with the boundary maps $(b+B)$ and $s$. Its is related to the $(b+B)$-cocycle $\gamma(\Tr\nat[2q]\uh^{-1}\dd\uh)$ via the descent equation
$$
(b+B)\gamma(\Tr[2q]\omh)+s\gamma(\Tr\nat[2q]\uh^{-1}\dd\uh)=0\ ,
$$
which can be shown either by direct computation, or simply by observing that $\nat\dd(\Tr[2q]\om)+s(\Tr\nat[2q]\uh^{-1}\dd\uh)=0$ and $\gamma\nat\dd= (b+B)\gamma$, $\gamma s=s\gamma$. Next we have to evaluate the image of $\gamma (\Tr[2q]\,\omh)$ by the chain map $\chih^p\rho_*:\Omh\Th\Ac\to X(\Rch)$, whose even component evaluated on a $p$-form $x_0\dd x_1\ldots\dd x_p$ over $\Th\Ac$ reads
$$
\frac{q!}{(p+1)!} \sum_{\la\in S_{p+1}}\eps(\la)  \tau'(\rho_*x_{\la(0)}[F,\rho_*x_{\la(1)}] \ldots [F,\rho_*x_{\la(p)}])\ .
$$
Denote as before $\ut=\uh_+\uh_-^{-1} \in (\Th\cc\hotimes\Ic\Rch[0,1]_x)^+$, and $\omt=\ut^{-1}s\ut$. One finds:
\beq
\lefteqn{\ch^p(\rho)\,(\Tr[2q]\,\omh)= \frac{(q!)^2}{(p+1)!}\, \Tr[2q]\, \Big( \uh_-^{-1} \omt [(\ut^{-1}-1)(\ut-1)]^q\uh_-+ }\non\\
&&(\ut-1)\omt [(\ut^{-1}-1)(\ut-1)]^{q-1}(\ut^{-1}-1)+\ldots + \uh_-^{-1}[(\ut^{-1}-1)(\ut-1)]^q\omt \uh_- \Big)\non
\eeq
After evaluation on the current $\frac{1}{\sqrt{2\pi i}} \int_{x=0}^1$, the right-hand-side may be identified, modulo commutators, with the Chern-Simons form $\cs^p_0(\ut)$ defined in analogy with (\ref{hcs}). One thus gets
$$
\ch^p(\rho)\cdot \cs^p_0(\gh) \equiv \cs^p_0(\ut)\mod\bb\ .
$$
Now, introduce a parameter $t\in [0,1]$ and choose an invertible interpolation $\vh \in (\Th\cc\hotimes\Ic\Rch[0,1]_x[0,1]_t)^+$ between $\vh_{t=0}=\ut$ and $\vh_{t=1}=1+\eh\otimes(\widehat{\rho_!(g)}-1)$, with the property that it restricts to $\vh(g_0)$ for $x=0$ and to $\vh(g_1)$ for $x=1$. The projection of $\vh$ on the algebra $(\cc\hotimes\Ic\Bc[0,1]_x[0,1]_t)^+$ may be chosen constant with respect to $t$. In the proof of Lemma \ref{lbott2} we established at any point $(x,t)\in [0,1]^2$ the identity
$$
\frac{\d}{\d t}\big( \Tr[2q]\, \vh^{-1}[(\vh-1)(\vh^{-1}-1)]^q s\vh \big) \equiv s \big(\Tr[2q]\, \vh^{-1}[(\vh-1)(\vh^{-1}-1)]^q\frac{\d\vh}{\d t}\big)\mod\bb\ ,
$$
and integrating over $[0,1]^2$ this implies
$$
\cs^p_0(\vh_{t=1})-\cs^p_0(\vh_{t=0}) \equiv \cs^p_0(\vh_{x=1})-\cs^p_0(\vh_{x=0}) \mod\bb\ .
$$
Taking into account that $\vh_{x=0}=\vh(g_0)$ and $\vh_{x=1}=\vh(g_1)$, we calculate mod $\bb$ in the complex $X_{n-p-1}(\Rc,\Jc)$
\beq
\lefteqn{\ch^p(\rho)\cdot(\te_1-\te_0) +\cs^p_0(\vh(g_1))-\cs^p_0(\vh(g_0))}\non\\
&&\qquad\qquad \equiv\ch^p(\rho)\cdot \cs^p_0(\gh)+ \cs^p_0(\vh_{t=1})-\cs^p_0(\vh_{t=0}) \mod\bb\non\\
&&\qquad\qquad \equiv\cs^p_0(\ut) + \cs^p_0(1+\eh\otimes(\widehat{\rho_!(g)}-1))-\cs^p_0(\ut) \mod\bb\non\\
&&\qquad\qquad \equiv\cs^p_0(\widehat{\rho_!(g)})\mod\bb \ .\non
\eeq
Hence the direct images $\rho_!(\gh_0,\te_0)$ and $\rho_!(\gh_1,\te_1)$ are equivalent and the map $\rho_!: MK^{\Ic}_n(\Ac)\to MK^{\Ic}_{n-p}(\Bc)$ for $n=2k+1$ and $p=2q$ is well-defined. It is obviously compatible with the push-forward map in topological $K$-theory $\rho_!:\Kt_1(\Ic\Ac)\to\Kt_1(\Ic\Bc)$. The compatibility with the push-forward map in cyclic homology $\ch^p(\rho): HC_{n-1}(\Ac)\to  HC_{n-p-1}(\Bc)$ is clear once we remark that the Chern-Simons form $\cs^{2q}_0(\vh)$ vanish whenever $\vh=1$. Hence the diagram (\ref{rr}) is commutative. \\ 
We have to check the invariance of $\rho_!$ with respect to conjugation of quasihomomorphisms. Let $\rho_0$ and $\rho_1$ be two conjugate quasihomomorphisms $\Ac\to \Ec^s\triangleright \Ic^s\Bc$. Hence there exists an invertible element $U\in (\Ec^s_+)^+$ such that $\rho_1=U^{-1}\rho_0 U$. We follow the proof of Proposition \ref{pinv} and remark that the lifting homomorphisms $\rho_{0*},\rho_{1*}: \Th\Ac\to \Mc^s_+$ factor through homomorphisms $\varphi_0,\varphi_1:\Th\Ac\to\Th\Ec^s_+$. The maps $\rho_{0!},\rho_{1!}:MK^{\Ic}_n(\Ac)\to MK^{\Ic}_{n-p}(\Bc)$ are obtained by composition of the pushforward maps $\varphi_{0!},\varphi_{1!}: MK^{\Ic}_n(\Ac)\to MK^{\Ic}_n(\Ec^s_+)$ induced by the homomorphisms
$$
\varphi_{i!}(\gh, \te) = (\varphi_i(\gh), \varphi_i(\te))
$$
with the map $MK^{\Ic}_n(\Ec^s_+)\to MK^{\Ic}_{n-p}(\Bc)$ associated with the natural $(p+1)$-summable quasihomomorphism of even degree $\Ec^s_+\to \Ec^s\triangleright \Ic^s\Bc$. Hence it is sufficient to check that the maps $\varphi_{0!}$ and $\varphi_{1!}: MK^{\Ic}_n(\Ac)\to MK^{\Ic}_n(\Ec^s_+)$ coincide. From the proof of \ref{pinv} we know that $\varphi_1$ is smoothly homotopic to $\Uh^{-1}\varphi_0 \Uh$, where $\Uh\in (\Th\Ec^s_+)^+$ is a lifting of $U$, and the interpolating homomorphism $\varphi:\Th\Ac\to \Th\Ec^s_+[0,1]$ is constant modulo the ideal $\Jh\Ec^s_+$. Consequently the morphisms $X(\varphi_1)$ and $X(\Uh^{-1}\varphi_0 \Uh): X(\Th\Ac)\to X(\Th\Ec^s_+)$ are homotopic,
$$
X(\varphi_1) - X(\Uh^{-1}\varphi_0 \Uh) = [\d,H]
$$
with $H\in \hom^0(X(T\Ac,J\Ac),X(T\Ec^s_+,J\Ec^s_+))$ a cochain \emph{of order zero}. Let us now compare the images of $(\gh,\te)\in MK^{\Ic}_n(\Ac)$ under the pushworwards $\varphi_{1!}$ and $(\Uh^{-1}\varphi_0 \Uh)_!$. The images $\varphi_1(\gh)$ and $\Uh^{-1}\varphi_0(\gh)\Uh$ are smoothly homotopic, with interpolation $\varphi(\gh)\in (\Ic\Th\Ec^s_+[0,1])^+$. If moreover $\gh$ is the canonical lift of an invertible element $g\in (\Ic\Ac)^+$, Chern-Simons form associated to $\varphi(\gh)$ can be written as
\beq
\cs^{p}_0(\varphi(\gh)) &=& \frac{1}{\sqrt{2\pi i}}\frac{(q!)^2}{p!}\int_0^1 \Tr\, \varphi(\gh^{-1})\varphi[(\gh-1)(\gh^{-1}-1)]^q s\varphi(\gh)  \non\\
&=& \frac{1}{\sqrt{2\pi i}}\int_0^1 \Tr[2q]\, \varphi(\uh^{-1}) s\varphi(\uh) \ , \non
\eeq
where as usual $\uh=1+e\otimes(\gh-1)$ is the invertible associated to $\gh$, with $\varphi(\uh)=1+e\otimes(\varphi(\gh)-1) \in (\Th\cc\hotimes\Ic\Th\Ec^s_+[0,1])^+$. By construction (Proposition \ref{pinv} i)), the r.h.s. coincides with the evaluation of $H$ on the Chern character $\ch^{p}_1(\gh)=\frac{1}{\sqrt{2\pi i}} \Tr\nat[2q]\, \uh^{-1}\dd \uh$, and because $H$ is of order zero one has
$$
\cs^{p}_0(\varphi(\gh)) = H\ch^{p}_1(\gh) = H (\nat\dd \te) \equiv \varphi_1(\te) - (\Uh^{-1}\varphi_0 \Uh)(\te) \mod \bb 
$$
in the complex $X_{n-1}(T\Ec^s_+,J\Ec^s_+)$. This proves that $\varphi_{1!}(\gh,\te)$ is equivalent to $(\Uh^{-1}\varphi_0 \Uh)_!(\gh,\te)$. Now it remains to show that the image $(\Uh^{-1}\varphi_0 \Uh)_!(\gh,\te)=(\Uh^{-1}\varphi_0(\gh)\Uh, \Uh^{-1}\varphi_0(\te)\Uh)$ is equivalent to $\varphi_{0!}(\gh,\te)=(\varphi_0(\gh),\varphi_0(\te))$. Here we mimic the proof of Lemma \ref{linv} and construct a homotopy between the invertible matrices $\bigl( \begin{smallmatrix} \varphi_0(\gh) & 0 \\ 0 & 1 \end{smallmatrix} \bigr)$ and $\bigl( \begin{smallmatrix} \Uh^{-1}\varphi_0(\gh)\Uh & 0 \\ 0 & 1 \end{smallmatrix} \bigr)$, whose associated Chern-Simons form is a $\bb$-boundary. Since $\Uh^{-1}\varphi_0(\te)\Uh\equiv \varphi_0(\te) \mod \bb$, we conclude that $\rho_{0!}$ and $\rho_{1!}$ agree on topological and multiplicative $K$-theories. \\ 
Finally we have to check the compatibility with the negative Chern character $MK^{\Ic}_*\to HN_*$. For any pair $(\gh,\te)\in MK^{\Ic}_n(\Ac)$, one has
$$
\ch_n(\gh,\te)=\ch^p_1(\gh)-\nat\dd \tilde\te \ \in F^{n-1}\Xh(T\Ac,J\Ac) ,
$$
where $\tilde\te$ is any lift of $\te$ in $\Xh(T\Ac,J\Ac)$. On the other hand, if $\gh$ is the canonical lift of some $g\in (\Ic\Ac)^+$, its image $\rho_!(\gh,\te) \in MK^{\Ic}_{n-p}(\Bc)$ is represented by the pair $(\widehat{\rho_!(g)}, \ch^p(\rho)\cdot\te +\cs^p_0(\vh))$ constructed above, so that
$$
\ch_{n-p}\big( \rho_!(\gh,\te)\big)=\ch^p_1\big(\widehat{\rho_!(g)}\big) - \nat\dd(\ch^p(\rho)\cdot\tilde\te +\cs^p_0(\vh))\ \in F^{n-p-1}\Xh(\Rc,\Jc)\ .
$$
But we know that the relation $\ch^p_1\big(\widehat{\rho_!(g)}\big)- \nat\dd \cs^p_0(\vh) = \ch^p(\rho)\cdot \ch^p_1(\gh) $ actually holds in the complex $\Xh(\Rc,\Jc)=X(\Rch)$. Therefore
$$
\ch_{n-p}\big( \rho_!(\gh,\te)\big) = \ch^p(\rho)\cdot(\ch^p_1(\gh)-\nat\dd\tilde\te)= \ch^p(\rho)\cdot \ch_n(\gh,\te)\ ,
$$
and (\ref{rr}) is compatible with the diagram of Theorem \ref{tbiv}.\\

\noindent {\bf ii) $n=2k$ is even and $p=2q$ is even.} As in the case of topological $K$-theory we pass to the suspensions of $\Ac$ and $\Bc$. We shall only sketch the procedure. The multiplicative $K$-theory group of even degree $MK^{\Ic}_n(\Ac)$ has an alternative description in terms of the set $MK'^{\Ic}_n(\Ac)$ of equivalence classes of pairs $(\gh,\te)$, where $\gh\in (\Kc\hotimes\Ic S\Th\Ac)^+$ is an invertible and $\te\in X_{n-1}(T\Ac,J\Ac)$ is a chain of odd degree such that $\cs^p_0(\gh)=\bb \te$. The equivalence relation is based on a higher transgression of the Chern-Simons form: $(\gh_0,\te_0)$ is equivalent to $(\gh_1,\te_1)$ iff there exists an invertible interpolation $\gh\in (\Kc\hotimes\Ic S\Th\Ac[0,1])^+$ and a chain of even degree $\la$ such that 
$$
\te_1-\te_0=\cs'^p_1(\gh) + \nat\dd\la\  \in X_{n-1}(T\Ac,J\Ac)\ ,
$$
where the odd chain $\cs'^p_1(\gh)\in X(\Th\Ac)$ is defined modulo $\nat\dd$ by the higher transgression formula (see the proof of Lemma \ref{lbott2})
$$
\bb \cs'^p_1(\gh)=\cs^p_0(\gh_1)-\cs^p_0(\gh_0)\ .
$$
Like $MK$, one can show that $MK'^{\Ic}_n(\Ac)$ is an abelian group inserted between $HC_*(\Ac)$ and $\Kt_*(\Ic S\Ac)$ in an exact sequence. More precisely there is a commutative diagram with exact rows:
$$
\xymatrix@C=1.2pc{
\Kt_1(\Ic\Ac) \ar[r]^{\ch_1} \ar[d]^{\al} & HC_{n-1}(\Ac) \ar[r]^{\delta} \ar[d]^{\times \sqrt{2\pi i}} & MK^{\Ic}_n(\Ac)  \ar[r] \ar[d]  & \Kt_0(\Ic\Ac) \ar[r]^{\ch_0} \ar[d]_{\beta} & HC_{n-2}(\Ac) \ar[d]_{\times \sqrt{2\pi i}}  \\
\Kt_0(\Ic S\Ac) \ar[r]^{\cs_1}  & HC_{n-1}(\Ac) \ar[r]^{\delta}  & MK'^{\Ic}_{n}(\Ac)  \ar[r]  & \Kt_1(\Ic S\Ac) \ar[r]^{\cs_0} & HC_{n-2}(\Ac) }
$$
Because for even $n$ the group $MK'_n$ is constructed from invertibles, it has \emph{odd} parity by convention. The (odd) map $MK^{\Ic}_n(\Ac)\to MK'^{\Ic}_n(\Ac)$ sends a pair $(\eh,\te)$ to $(\beta(\eh),\sqrt{2\pi i}\te + l^p_1(\eh))$, where $\beta(\eh)=(1+(z-1)\eh)(1+(z-1)p_0)^{-1}$ is the invertible image of $\eh$ under the Bott map, and $l^p_1(\eh)$ is the transgressed cochain defined modulo $\nat\dd$ by $\bb (l^p_1(\eh))=\cs^p_0(\beta(\eh))- \sqrt{2\pi i}\, \ch^p_0(\eh)$. The map $MK'^{\Ic}_n(\Ac)\to \Kt_1(\Ic S\Ac)$ is the forgetful map, and $HC_{n-1}(\Ac)\to MK'^{\Ic}_n(\Ac)$ sends a cycle $\te$ to $(1, \sqrt{2\pi i}\, \te)$. By the five lemma, $MK'^{\Ic}_n(\Ac)$ is thus isomorphic to $MK^{\Ic}_n(\Ac)$. One easily checks that the negative Chern character $\ch_n:MK'^{\Ic}_n(\Ac)\to HN_n(\Ac)$ given by $\ch_n(\gh,\te)=\cs^p_0(\gh)-\bb\tilde\te$ coincides with the negative Chern character on $MK^{\Ic}_n(\Ac)$ up to a factor $\sqrt{2\pi i}$.\\ Hence it suffices to construct the pushforward morphism $\rho_!$ for the groups $MK'_n$, whose elements are represented by \emph{invertibles} of the suspended algebras:
$$
\xymatrix{
\Kt_0(\Ic S\Ac) \ar[r] \ar[d]^{\rho_!} & HC_{n-1}(\Ac) \ar[r] \ar[d]^{\ch^p(\rho)} & MK'^{\Ic}_n(\Ac)  \ar[r] \ar[d]^{\rho_!}  & \Kt_1(\Ic S\Ac)  \ar[d]^{\rho_!}  \\
\Kt_{-p}(\Ic S\Bc) \ar[r]  & HC_{n-1-p}(\Bc) \ar[r]  & MK'^{\Ic}_{n-p}(\Bc)  \ar[r]  & \Kt_{1-p}(\Ic S\Bc)  }
$$
This can be done explicitly as in case i), with the only difference that the Chern character $\ch^p_1$ and Chern-Simons transgression $\cs^p_0$ are now replaced respectively by $\cs^p_0$ and the higher transgression $\cs'^p_1$. The needed formulas were already established in i): let $g\in (\Ic S\Ac)^+$ be any invertible with canonical lift $\gh\in (\Ic S\Th\Ac)^+$. One can write
$$
\ch^p(\rho)\cdot \cs^p_0(\gh) = \cs^p_0(\ut) - \bb k^p_1(\uh)
$$
with the invertibles $\uh=1+e\otimes(\gh-1) \in (\Th\cc\hotimes \Ic S\Th\Ac)^+$ and $\ut=\uh_+\uh_-^{-1} \in (\Th\cc\hotimes \Ic S\Rch)^+$, and $k^p_1(\uh)$ is a chain defined mod $\nat\dd$. Let $\widehat{\rho_!(g)} \in (\Ic S\Rch)^+$ be any invertible lift of $\rho_!(g)=\rho_+(g)\rho_-(g)^{-1}$, and $\vh \in (\Th\cc\hotimes \Ic S\Rch [0,1])^+$ be an invertible interpolation between $\vh_0=\ut$ and $\vh_1=1+\eh\otimes(\widehat{\rho_!(g)}-1)$. Then one has
$$
\bb \cs'^p_1(\vh)= \cs^p_0(\widehat{\rho_!(g)}) - \cs^p_0(\ut)\ .
$$
Therefore if $(\gh,\te)$ represents a class in $MK'^{\Ic}_n(\Ac)$ we define its pushforward as the multiplicative $K$-theory class over $\Bc$
\be
\boxed{
\rho_!(\gh,\te) = \big(\widehat{\rho_!(g)}\ ,\ \ch^p(\rho)\cdot \te + k^p_1(\uh) + \cs'^p_1(\vh) \big) \ \in MK'^{\Ic}_{n-p}(\Bc)\  }
\ee 
with the odd chain $\ch^p(\rho)\cdot \te + k^p_1(\uh) + \cs'^p_1(\vh)$ sitting in $X_{n-p-1}(\Rc,\Jc)$. One shows the consistency of $\rho_!$ with the various equivalence relations using the properties of higher Chern-Simons transgressions. Details are left to the reader. \\

\noindent {\bf iii) $n=2k+1$ is odd and $p=2q+1$ is odd.} We first establish an explicit formula for the composition of the topological Chern character $\ch^{2q}_1:\Kt_1(\Ic\Ac)\to HP_1(\Ac)$ with the bivariant Chern character $\ch^p(\rho)\in HC^p(\Ac,\Bc)$. Remark that $\Ic$ is $(2q+1)$-summable by hypothesis hence $\ch^{2q}_1$ is well-defined. As in case i) let $g\in (\Ic\Ac)^+$ be an invertible, $\gh\in (\Ic\Th\Ac)^+$ its canonical lift and $\uh=1+e\otimes(\gh-1)\in (\Th\cc\hotimes\Ic\Th\Ac)^+$ the associated invertible. Recall that
$$
\ch^{2q}_1(\gh)=\frac{1}{\sqrt{2\pi i}}\, \Tr\nat[2q]\, \uh^{-1}\dd\uh\ \in \Om^1\Th\Ac_{\nat}\ ,
$$
and the image of $\Tr\nat[2q]\, \uh^{-1}\dd\uh$ under the Goodwillie equivalence is the $(b+B)$-cycle over $\Th\Ac$
$$
\gamma (\Tr\nat[2q]\, \uh^{-1}\dd\uh) = \sum_{i\geq 0}(-)^i i! \, \Tr[2q]\, \uh^{-1}\dd\uh (\dd\uh^{-1}\dd\uh)^i\ .
$$
Now the quasihomomorphism $\rho:\Ac\to \Ec^s\triangleright \Ic^s\Bc$ is of odd degree. Hence, the image of an element $x\in \Th\Ac$ under the lifted quasihomomorphism $\rho_*:\Th\Ac\to \Mch^s\triangleright \Ic^s\Rch$ is a $2\times 2$ matrix over $\Mch$ whose off-diagonal entries lie in $\Ic\Rch$. Moreover the multiplier $F$ is given by the matrix 
$$
F=\eps\left(\begin{matrix}
1 & 0 \\
0 & -1 \end{matrix} \right)=\eps(2p_0-1)
$$
where $\eps$ is the odd generator of the Clifford algebra $C_1$. Thus the commutator $[p_0,\rho_*x]$ lies in the matrix algebra $M_2(\Ic\Rch)$ for any $x\in\Th\Ac$. On the other hand, the component of the chain map $\chih^p\rho_*:\Omh \Th\Ac\to X(\Rch)$ evaluated on a $p$-form $x_0\dd x_1\ldots\dd x_p$ reads
$$
 -\frac{\Gamma(q+\frac{3}{2})}{(p+1)!} \sum_{\la\in S_{p+1}}\eps(\la)  \tau(\rho_*x_{\la(0)}[F,\rho_*x_{\la(1)}] \ldots [F,\rho_*x_{\la(p)}])\ ,
$$
where $\tau(\eps \cdot)=-\sqrt{2 i}\,\Tr(\cdot)$ is the odd supertrace (see section \ref{sbiv}). As in case i), let us extend $\rho_*$ to a unital homomorphism $(\Th\cc\hotimes\Ic\Th\Ac)^+\to (\Th\cc\hotimes\Ic\Mch^s_+)^+$. Using $\Gamma(q+\frac{3}{2})=\sqrt{\pi}\, p!/(2^pq!)$ with $p=2q+1$, one gets by direct computation
\beq
\ch^p(\rho)\cdot\ch^{2q}_1(\gh) &=& \frac{1}{2}\Tr[2q]\, (\ut^{-1}[p_0,\ut])^p + \frac{1}{2}\Tr[2q]\, ([p_0,\ut]\ut^{-1})^p\non\\
&=& \Tr[2q]\, (\ut^{-1}[p_0,\ut])^p -\bb \frac{1}{2}\Tr\nat[2q]\, (\ut^{-1}[p_0,\ut])^p\ut^{-1}\dd\ut \non
\eeq
where $\ut=\rho_*\uh$ is an invertible element of the algebra $(\Th\cc\hotimes\Ic\Mch^s_+)^+$, and the commutator $[p_0,\ut]\in M_2(\Th\cc\hotimes\Ic\Ic\Rch)^+$ may be considered as an element of $M_2(\Th\cc\hotimes\Ic\Rch)^+$ after applying the homomorphism $\boxtimes:\Ic\hotimes\Ic\to\Ic$. The first term of the r.h.s. is recognized as the higher Chern character $\ch^{2q}_0(\ft)=\Tr[2q](\ft-p_0)^p$ given by (\ref{hch}) for the idempotent $\ft=\ut^{-1}p_0\ut \in M_2(\Th\cc\hotimes\Ic\Ic\Rch)^+$ (or $M_2(\Th\cc\hotimes\Ic\Rch)^+$), whence the equality
$$
\ch^p(\rho)\cdot\ch^{2q}_1(\gh) = \ch^{2q}_0(\ft) - \bb\, \frac{1}{2}\Tr\nat[2q]\, (\ft-p_0)^p\ut^{-1}\dd\ut
$$
of cycles in $X(\Rch)$. Then, observe that the projection of $\ft$ to the algebra $M_2(\cc\hotimes\Ic\Bc)^+$ is the idempotent $p_0+e\otimes(\rho(g)^{-1}p_0\rho(g)-p_0)$. Using the isomorphism $(\cc\hotimes\Ic\Bc)^+\cong (\Ic\Bc)^+$, this idempotent may be identified with the direct image $\rho_!(g)=\rho(g)^{-1}p_0\rho(g)$. Hence, it is possible to relate $\ch^{2q}_0(\ft)$ with the Chern character of a given idempotent lift $\widehat{\rho_!(g)}\in M_2(\Ic\Rch)^+$, via a homotopy with parameter $t\in[0,1]$. Let $\fh\in M_2(\Th\cc\hotimes\Ic\Rch[0,1])^+$ be an idempotent path lifting the constant family $p_0+e\otimes(\rho_!(g)-p_0)$ and connecting the two endpoints 
$$
\fh_0=\ft\ ,\qquad \fh_1=p_0+\eh\otimes(\widehat{\rho_!(g)}-p_0)\ .
$$
$\eh\in\Th\cc$ is the canonical idempotent lift of the unit $e\in\cc$ as in case i). The lifting $\fh$ is thus defined up to homotopy (at least after stabilization by the matrix algebra $\Kc$). The property $[2q]\,\eh=1$ implies the equality
$$
\ch^{2q}_0(\fh_1)=\ch^{2q}_0\big(\widehat{\rho_!(g)}\big) \in \Rch
$$
at the level of cycles. Furthermore, in analogy with Eqs. (\ref{hcs}) the Chern-Simons form associated to the idempotent $\fh$ is defined by
$$
\cs^{2q}_1(\fh)=\int_0^1dt\, \Tr\nat[2q]\, (-2\fh+1)\sum_{i=0}^q(\fh-p_0)^{2i}\frac{\d\fh}{\d t} (\fh-p_0)^{2(q-i)}\dd\fh\ ,
$$
and fulfills the transgression relation in $X(\Rch)$
$$
\bb\cs^{2q}_1(\fh)=\ch^{2q}_0(\fh_1)-\ch^{2q}_0(\fh_0)=\ch^{2q}_0(\widehat{\rho_!(g)})-\ch^{2q}_0(\ft)\ .
$$
This leads to the definition of the map $\rho_!$ on multiplicative $K$-theory. Let $(\gh,\te)$ represent a class in $MK^{\Ic}_n(\Ac)$ of odd degree $n=2k+1$. By Remark \ref{rcan} we may suppose that $\gh$ is the canonical lift of some invertible $g\in(\Ic\Ac)^+$, and $\te\in X_{n-1}(T\Ac,J\Ac)$ is a transgression of the Chern character $\ch^{2q}_1(\gh)=\nat\dd\te$. We set
\be
\boxed{ \rho_!(\gh,\te)=\big(\widehat{\rho_!(g)}\ ,\ -\ch^p(\rho)\cdot \te + h^{2q}_1(\ut) + \cs^{2q}_1(\fh)\big)\ \in MK^{\Ic}_{n-p}(\Bc) }
\ee
where $\widehat{\rho_!(g)}\in M_2(\Ic\Rch)^+$ is any idempotent lift of $\rho_!(g)$, $h^{2q}_1(\ut)$ is the chain $\frac{1}{2}\Tr\nat[2q]\, (\ft-p_0)^p\ut^{-1}\dd\ut$, and $\fh$ is an idempotent path constructed as above. The minus sign in front of $\ch^p(\rho)\cdot\te$ is necessary because the bivariant Chern character $\ch^p(\rho)$ is of odd degree $p=2q+1$. This ensures the correct transgression relation
\beq
\lefteqn{\bb(-\ch^p(\rho)\cdot \te + h^{2q}_1(\ut)+ \cs^{2q}_1(\fh))}\non\\
&=& \ch^p(\rho)\cdot\ch^{2q}_1(\gh) +\bb\,\frac{1}{2}\Tr\nat[2q]\, (\ft-p_0)^p\ut^{-1}\dd\ut + \ch^{2q}_0(\widehat{\rho_!(g)})-\ch^{2q}_0(\ft)\non\\
&=& \ch^{2q}_0(\widehat{\rho_!(g)}) \non
\eeq
in the quotient complex $X_{n-p-1}(\Rc,\Jc)$, which shows that $\rho_!(\gh,\te)$ indeed defines an element of $MK^{\Ic}_{n-p}(\Bc)$. Its class does not dependent on the chosen idempotent lift $\widehat{\rho_!(g)}$ nor on the path $\fh$, and moreover $\rho_!$ is compatible with the equivalence relation on multiplicative $K$-theory. We proceed as in case i) and let $(\gh_0,\te_0)$ and $(\gh_1,\te_1)$ be two equivalent representatives of a class in $MK^{\Ic}_n(\Ac)$, provided with an interpolation $\gh\in (\Ic\Th\Ac[0,1]_x)^+$ and a chain $\la\in X_{n-1}(T\Ac,J\Ac)$ such that $\te_1-\te_0=\cs^{2q}_0(\gh)+\bb\la$. From Remark \ref{rcan} the elements $\gh_0$, $\gh_1$ and $\gh$ can be taken as the canonical lifts of $g_0,g_1\in (\Ic\Ac)^+$ and $g\in (\Ic\Ac[0,1]_x)^+$. Denoting by $\rho_*\uh(g_i)=\ut(g_i)\in (\Th\cc\hotimes\Ic\Mch^s_+)^+$ the invertible and $\fh(g_i)\in M_2(\Th\cc\hotimes\Ic\Rch[0,1]_t)^+$ the idempotent path associated to $g_i$, we have to establish the relation
\beq
&&-\ch^p(\rho)\cdot(\te_1-\te_0) + h^{2q}_1(\ut(g_1)) - h^{2q}_1(\ut(g_0)) + \cs^{2q}_1(\fh(g_1)) - \cs^{2q}_1(\fh(g_0))\non\\
&&\qquad\qquad\qquad \equiv \cs^{2q}_1\big(\widehat{\rho_!(g)}\big) \mod \nat\dd \non
\eeq
in the complex $X_{n-p-1}(\Rc,\Jc)$, where $\widehat{\rho_!(g)}\in M_2(\Ic\Rch[0,1]_x)^+$ is a choice of idempotent interpolation between the liftings $\widehat{\rho_!(g_i)}$'s. As usual, let $\uh=1+e\otimes(\gh-1)\in (\cc\hotimes\Ic\Ac[0,1]_x)^+\hookrightarrow  (\Th\cc\hotimes\Ic\Th\Ac[0,1]_x)^+$ be the invertible identification with $\gh$. We know the equality
$$
\cs^{2q}_0(\gh) = \frac{1}{\sqrt{2\pi i}} \int_0^1 \Tr[2q]\, \uh^{-1}s\uh \ ,
$$
where $s$ is the de Rham differential on $\Om[0,1]_x$. Set $\omh=\uh^{-1}s\uh$. The computation of $\ch^p(\rho)\cdot\cs^{2q}_0(\gh)$ involves the formula 
\beq
\lefteqn{\gamma (\Tr[2q]\,\omh) =} \non\\
&& \Tr[2q]\,\omh+ \sum_{i\geq 1}(-)^i(i-1)!\sum_{j=0}^{i-1}\Tr[2q]\, \uh^{-1}\dd\uh(\dd\uh^{-1}\dd\uh)^j\dd\omh (\dd\uh^{-1}\dd\uh)^{i-j-1}\ ,\non
\eeq
as well as the component of the chain map $\chih^p\rho_*$ evaluated on a $(p+1)$-form $x_0\dd x_1\ldots\dd x_{p+1}$ over $\Th\Ac$:
$$
 -\frac{\Gamma(q+\frac{3}{2})}{(p+1)!} \sum_{i=1}^{p+1}  \tau\nat(\rho_*x_0[F,\rho_*x_1]\ldots\dd (\rho_*x_i) \ldots [F,\rho_*x_{p+1}])
$$
Denote as before $\ut=\rho_*\uh$ the invertible image in $(\Th\cc\hotimes\Ic\Mch^s_+[0,1]_x)^+$, the associated idempotent $\ft=\ut^{-1}p_0\ut \in M_2(\Th\cc\hotimes\Ic\Rch[0,1]_x)^+$, and the Maurer-Cartan form $\omt=\ut^{-1}s\ut$. One gets
\beq
\lefteqn{\ch^p(\rho)\,(\Tr[2q]\,\omh)=}\non\\
 && -\frac{\sqrt{2\pi i}}{2}\Tr\nat[2q]\, \Big( \sum_{i=0}^{2q}(\ft-p_0)^i[p_0,\omt](\ft-p_0)^{2q-i}\ut^{-1}\dd\ut +(\ft-p_0)^p\dd\omt\Big)\ .  \non
\eeq
Now observe that $\ut_{x=0}=\ut(g_0)$ and $\ut_{x=1}=\ut(g_1)$, so that after integration over the current $\frac{1}{\sqrt{2\pi i}} \int_{x=0}^1$ we get the identity (recall $\ch^p(\rho)$ is odd)
\beq
\lefteqn{-\ch^p(\rho)\cdot \cs^{2q}_0(\gh)+ h^{2q}_1(\ut(g_1))-h^{2q}_1(\ut(g_0))}\non\\
&=& \frac{1}{\sqrt{2\pi i}}\int_0^1 \ch^p(\rho)\,(\Tr[2q]\,\omh) +  \frac{1}{2}\int_0^1 s\, \Tr\nat[2q]\, (\ft-p_0)^p\ut^{-1}\dd\ut \non\\
&=& \int_0^1 \Tr\nat[2q]\, \Big( \sum_{i=1}^{q}(\ft-p_0)^{2i-1}[p_0,\omt](\ft-p_0)^{2(q-i)+1}\ut^{-1}\dd\ut - (\ft-p_0)^p\dd\omt \Big) \non
\eeq
On the other hand, let us calculate the Chern-Simons form associated to the idempotent $\ft$,
$$
\cs^{2q}_1(\ft)=\int_0^1 \Tr\nat[2q]\, (-2\ft+1)\sum_{i=0}^q(\ft-p_0)^{2i} s\ft (\ft-p_0)^{2(q-i)}\dd\ft\ ,
$$
in terms of $\omt$. Since by definition $\ft=\ut^{-1}p_0\ut$, the structure equation $s\ft=[\ft,\omt]$ follows and one finds
\beq
\lefteqn{\cs^{2q}_1(\ft) = -\nat\dd\int_0^1\Tr[2q]\, (\ft-p_0)^p\omt}\non\\
&+& \int_0^1 \Tr\nat[2q]\, \Big( \sum_{i=1}^{q}(\ft-p_0)^{2i-1}[p_0,\omt](\ft-p_0)^{2(q-i)+1}\ut^{-1}\dd\ut - (\ft-p_0)^p\dd\omt \Big)\non
\eeq
Thus holds the fundamental relation
$$
-\ch^p(\rho)\cdot \cs^{2q}_0(\gh)+ h^{2q}_1(\ut(g_1))-h^{2q}_1(\ut(g_0)) \equiv \cs^{2q}_1(\ft)\mod\nat\dd\ .
$$
Now let $\fh\in M_2(\Th\cc\hotimes\Ic\Rch[0,1]_x[0,1]_t)^+$ be an idempotent interpolation between $\fh_{t=0}=\ft$ and $\fh_{t=1}=p_0+\eh\otimes(\widehat{\rho_!(g)}-p_0)$, with the property that it restricts to $\fh(g_0)$ for $x=0$ and to $\fh(g_1)$ for $x=1$. The projection of $\fh$ to the algebra $M_2(\cc\hotimes\Ic\Bc[0,1]_x[0,1]_t)^+$ may be chosen constant with respect to $t$. In the proof of Lemma \ref{lbott2} we established the following identity at any point $(x,t)\in[0,1]^2$:
\beq
\lefteqn{\frac{\d}{\d t}\big( \Tr\nat[2q]\, (-2\fh+1)\sum_{i=0}^q(\fh-p_0)^{2i} s\fh (\fh-p_0)^{2(q-i)}\dd\fh \big) }\non\\
&\equiv& s \big( \Tr\nat[2q]\, (-2\fh+1)\sum_{i=0}^q(\fh-p_0)^{2i}\frac{\d\fh}{\d t} (\fh-p_0)^{2(q-i)}\dd\fh\big) \mod \nat\dd \ , \non
\eeq
and integration over the square $[0,1]^2$ yields 
$$
\cs^{2q}_1(\fh_{t=1})-\cs^{2q}_1(\fh_{t=0}) \equiv \cs^{2q}_1(\fh_{x=1})-\cs^{2q}_1(\fh_{x=0}) \mod\nat\dd\ .
$$
Since $\fh_{x=0}=\fh(g_0)$ and $\fh_{x=1}=\fh(g_1)$ we calculate, modulo $\nat\dd$ in the complex $X_{n-p-1}(\Rc,\Jc)$
\beq
\lefteqn{-\ch^p(\rho)\cdot(\te_1-\te_0) + h^{2q}_1(\ut(g_1)) - h^{2q}_1(\ut(g_0)) + \cs^{2q}_1(\fh(g_1)) - \cs^{2q}_1(\fh(g_0))}\non\\
&\equiv& -\ch^p(\rho)\cdot \cs^{2q}_0(\gh)+ h^{2q}_1(\ut(g_1))-h^{2q}_1(\ut(g_0)) + \cs^{2q}_1(\fh_{t=1})-\cs^{2q}_1(\fh_{t=0}) \non\\
&\equiv& \cs^{2q}_1(\ft) + \cs^{2q}_1(p_0+\eh\otimes(\widehat{\rho_!(g)}-p_0)) - \cs^{2q}_1(\ft)  \non\\
&\equiv& \cs^{2q}_1(\widehat{\rho_!(g)}) \mod\nat\dd\non
\eeq
as wanted. Hence $\rho_!(\gh_0,\te_0)$ and $\rho_!(\gh_1,\te_1)$ are equivalent and the map $\rho_!: MK^{\Ic}_n(\Ac)\to MK^{\Ic}_{n-p}(\Bc)$ for $n=2k+1$ and $p=2q+1$ is well-defined. Its compatibility with the map $\rho_!$ on topological $K$-theory is obvious. Concerning its compatibility with the map $\ch^p(\rho):HC_{n-1}(\Ac)\to HC_{n-p-1}(\Bc)$, we should take care of a minus sign which shows that the middle square of (\ref{rr}) is actually anticommutative; this has to be so because all the maps involved in this square are of odd degree. Hence the diagram (\ref{rr}) is graded commutative.\\
The invariance of $\rho_!$ with respect to conjugation of quasihomomorphisms is proved exactly as in case i), by decomposing $\rho_!$ as the pushforward map $\varphi_!:MK^{\Ic}_n(\Ac)\to MK^{\Ic}_n(\Ec^s_+)$ induced by the homomorphism $\varphi: \Th\Ac\to \Th\Ec^s_+$, followed by the map $MK^{\Ic}_n(\Ec^s_+)\to MK^{\Ic}_{n-p}(\Bc)$ associated with the natural $p$-summable quasihomomorphism of odd degree $\Ec^s_+\to \Ec^s\triangleright \Ic^s\Bc$. Also the compatibility with the negative Chern character is easily established.\\ 

\noindent {\bf iv) $n=2k$ is even and $p=2q+1$ is odd.} As in case ii) we pass to the suspensions of $\Ac$ and $\Bc$ and work with the group $MK'^{\Ic}_n(\Ac)$. Hence a multiplicative $K$-theory class of degree $n$ over $\Ac$ is represented by a pair $(\gh,\te)$ of an invertible $\gh\in (\Kc\hotimes\Ic S\Th\Ac)^+$ and an odd chain $\te\in X_{n-1}(T\Ac,J\Ac)$ such that $\cs^{2q}_0(\gh)=\bb \te$. We are thus led to build a morphism
$$
\rho_!: MK'^{\Ic}_n(\Ac) \to MK'^{\Ic}_{n-p}(\Bc)\ ,
$$
where the group $MK'^{\Ic}_{n-p}(\Bc)$, for $n-p$ odd, is represented by pairs $(\eh,\te)$ of idempotent $\eh\in M_2(\Kc\hotimes\Ic S\Rch)^+$ and chain of even degree $\te\in X_{n-p-1}(\Rc,\Jc)$ such that $ \cs^{2q}_1(\eh)=\nat\dd \te$. Note that the parity of $MK'^{\Ic}_{n-p}(\Bc)$ is \emph{even}. We already established the needed formulas in case iii): let $g\in (\Ic S\Ac)^+$ be any invertible with canonical lift $\gh\in (\Ic S\Th\Ac)^+$. One can write
$$
-\ch^p(\rho)\cdot \cs^{2q}_0(\gh) = \cs^{2q}_1(\ft) - \nat\dd k^{2q}_0(\uh)
$$
with the invertible $\uh=1+e\otimes(\gh-1) \in (\Th\cc\hotimes \Ic S\Th\Ac)^+$, the idempotent $\ft=\ut^{-1}p_0\ut \in M_2(\Th\cc\hotimes \Ic S\Rch)^+$ where $\ut=\rho_*\uh$, and the chain $k^{2q}_0(\uh)=-\int_0^1 \Tr[2q]\, (\ft-p_0)^p\omt$ where $\omt=\ut^{-1}s\ut$. Let $\widehat{\rho_!(g)} \in M_2(\Ic S\Rch)^+$ be any idempotent lift of $\rho_!(g)=\rho(g)^{-1}p_0\rho(g)$, and $\fh \in M_2(\Th\cc\hotimes \Ic S\Rch [0,1])^+$ be an idempotent interpolation between $\fh_0=\ft$ and $\fh_1=p_0+\eh\otimes(\widehat{\rho_!(g)}-p_0)$. Then one has by means of the higher transgressions (see the proof of Lemma \ref{lbott2})
$$
\nat\dd \cs'^{2q}_0(\fh)= \cs^{2q}_1(\widehat{\rho_!(g)}) - \cs^{2q}_1(\ft)\ ,
$$
with $\cs'^{2q}_0(\fh)$ defined modulo $\bb$. Therefore if $(\gh,\te)$ represents a class in $MK'^{\Ic}_n(\Ac)$ we define its pushforward as the multiplicative $K$-theory class over $\Bc$
\be
\boxed{
\rho_!(\gh,\te) = \big(\widehat{\rho_!(g)}\ ,\ \ch^p(\rho)\cdot \te + k^{2q}_0(\uh) + \cs'^{2q}_0(\fh) \big) \ \in MK'^{\Ic}_{n-p}(\Bc)\  }
\ee 
with the chain $\ch^p(\rho)\cdot \te + k^{2q}_0(\uh) + \cs'^{2q}_0(\fh)$ of even degree sitting in the quotient complex $X_{n-p-1}(\Rc,\Jc)$.  \cqfd\\

\begin{example}\label{ereg}\textup{For $\Bc=\cc$ a quasihomomorphism $\Ac\to \Ec^s\triangleright \Ic^s$ induces a map $MK^{\Ic}_n(\Ac)\to MK^{\Ic}_{n-p}(\cc)$. Thus if $\Ic$ is a Schatten ideal on a Hilbert space, Example \ref{ec} yields index maps or regulators, depending on the degrees:
\beq
MK^{\Ic}_n(\Ac)\ \to &\zz &\quad \mbox{if}\ n\leq p\ ,\ n\equiv p \mod 2\ , \non\\
MK^{\Ic}_n(\Ac)\ \to &\cc^{\times} &\quad \mbox{if}\  n> p\ ,\ n\equiv p+1 \mod 2\ . \non
\eeq }
\end{example}

\section{Assembly maps and crossed products}\label{sassem}

In this section we illustrate the general theory of secondary characteristic classes with the specific example of crossed product algebras, and build an "assembly map" for multiplicative $K$-theory modelled on the Baum-Connes construction \cite{BC}.\\

Let $\Ac$ be a unital Fr\'echet $m$-algebra and $\Gamma$ a countable discrete group acting on $\Ac$ from the right by automorphisms. The action of an element $\gamma\in \Gamma$ on $a\in \Ac$ reads $a^{\gamma}$. We impose the action to be almost isometric in the following sense: for each submultiplicative seminorm $\|\cdot \|_{\al}$ on $\Ac$ there exists a constant $C_{\al}$ such that
\be
\|a^{\gamma}\|_{\al} \leq C_{\al} \|a \|_{\al} \qquad \forall\ a\in \Ac\ ,\ \gamma\in \Gamma\ . \label{cond}
\ee
The algebraic tensor product $\Ac\otimes \cc \Gamma$ is identified with the space of $\Ac$-valued functions with finite support over $\Gamma$. Thus any element of $\Ac\otimes \cc \Gamma$ is a finite linear combination of symbols $a \gamma^*$ with $a_{\gamma}\in \Ac$ and $\gamma \in  \Gamma$. The star refers to a contravariant notation. We endow $\Ac \otimes \cc \Gamma$ with the crossed product defined in terms of symbols by
$$
(a_1\gamma_1^*)(a_2\gamma_2^*)= a_1 (a_2)^{\gamma_1} (\gamma_2 \gamma_1)^* \qquad \forall\ a_i\in\Ac\ ,\ \gamma_i\in \Gamma\ .
$$
The crossed product algebra $\Ac\cp \Gamma$ is an adequate completion of the above space consisting of $\Ac$-valued functions with ``rapid decay" over $\Gamma$. This requires to fix once and for all a right-invariant distance function $d:\Gamma\times \Gamma\to \rr_+$. Endow the space $\Ac\otimes \cc \Gamma$ with the seminorms
$$
\|b\|_{\al,\beta}= C_{\al} \sum_{\gamma\in \Gamma} \si_{\beta}(\gamma) \|b(\gamma)\|_{\al}\qquad \forall\ b\in \Ac\otimes \cc \Gamma\ ,
$$
where the $\rr_+$-valued function $\si_{\beta}(\gamma):= (1+d(\gamma,1))^{\beta}$, for $\beta \geq 0$, fulfills the property $\si_{\beta}(\gamma_1\gamma_2)\leq \si_{\beta}(\gamma_1) \si_{\beta}(\gamma_2)$. One easily checks that $\|\cdot \|_{\al,\beta}$ is submultiplicative with respect to the crossed product, hence the completion of $\Ac\otimes \cc \Gamma$ for the family of seminorms $(\|\cdot \|_{\al,\beta})$ yields a Fr\'echet $m$-algebra $\Ac\cp \Gamma$.\\
Multiplicative $K$-theory classes of $\Ac\cp \Gamma$ may be obtained by adapting the assembly map construction of \cite{BC}. The idea is to replace the noncommutative space $\Ac\cp \Gamma$ by a more classical space, for which the secondary invariants are presumably easier to describe. Consider a compact Riemannian manifold $M$ without boundary, and let $P\stackrel{\Gamma}{\longrightarrow} M$ be a $\Gamma$-covering. $\Gamma$ acts on $P$ from the left by deck transformations. Denote by
$$
\Ac_P := \cinf(P;\Ac)^\Gamma
$$
the algebra of $\Gamma$-invariant smooth $\Ac$-valued functions over $P$: any function $a\in \Ac_P$ verifies $a(\gamma^{-1}\cdot x)=(a(x))^{\gamma}$, $\forall x\in P,\gamma\in\Gamma$. Thus $\Ac_P$ is the algebra of smooth sections of a non-trivial bundle with fibre $\Ac$ over $M$. It can be represented as a subalgebra of matrices over $\cinf(M)\hotimes (\Ac\cp \Gamma)=\cinf(M;\Ac\cp \Gamma)$ as follows. Let $(U_i)$, $i=1,\ldots, m$ be a finite open covering of $M$ trivializing the bundle $P$, via a set of sections $s_i:U_i\to P$ and locally constant transition functions $\gamma_{ij}:U_i\cap U_j\to \Gamma$:
$$
\gamma_{ij}\gamma_{jk}=\gamma_{ik}\ \mbox{over}\ U_i\cap U_j\cap U_k\ ,\quad s_i(x)=\gamma_{ij}\cdot s_j(x)\quad \forall x\in U_i\cap U_j\ .
$$
Choose a partition of unity $c_i\in \cinf(M)$ relative to this covering: $\supp c_i\subset U_i$ and $\sum_{i=1}^m c_i(x)^2=1$. From these data consider the homomorphism $\rho:\Ac_P\to M_m(\cinf(M)\hotimes(\Ac\cp \Gamma))$ sending an element $a\in \Ac_P$ to the $m\times m$ matrix $\rho(a)$ whose components, as $(\Ac\cp \Gamma)$-valued functions over $M$, read
$$
\rho(a)_{ij}(x):= c_i(x)c_j(x) a(s_i(x)) \gamma_{ji}^*\quad \forall\ i,j=1,\ldots,m\ ,\ \forall x\in M\ .
$$
Of course $\rho$ depends on the choice of trivialization $(U_i,s_i)$ and partition of unity $(c_i)$, but different choices are related by conjugation in a suitably large matrix algebra. Indeed, if $(U'_{\al},s'_{\al})$, $\al=1,\ldots,\mu$ denotes another trivialization with transition functions $\gamma'_{\al\beta}$ and partition of unity $(c'_{\al})$, we get a corresponding homomorphism $\rho': \Ac_P\to M_{\mu}(\cinf(M)\hotimes(\Ac\cp \Gamma))$. Introduce the rectangular matrices $u$, $v$ over $\cinf(M)\hotimes(\Ac\cp \Gamma)$ with components
$$
u_{i\al}(x)=c_i(x)c'_{\al}(x) \gamma_{\al i}^*\ ,\qquad v_{\al i}(x)=c'_{\al}(x)c_i(x) \gamma_{i\al}^*\ ,
$$
(recall $\Ac$ is unital by hypothesis hence $\cc\Gamma \subset \Ac\cp\Gamma$), where the mixed transition functions $\gamma_{i\al}$, $\gamma_{\al i}$ are defined by $s_i(x)=\gamma_{i\al}\cdot s'_{\al}(x)$ and $s'_{\al}(x)=\gamma_{\al i}\cdot s_{i}(x)$ for any $x\in U_i\cap U'_{\al}$. Then $uv$ and $vu$ are idempotent square matrices, and for any element $a\in \Ac_P$ one has $\rho(a)=u\rho'(a) v$ and $\rho'(a)=v\rho(a)u$. Moreover, the invertible square matrix of size $m+\mu$
$$
W=\left( \begin{matrix}
1-uv & -u \\
v & 1-vu \end{matrix} \right)\ ,
\qquad
W^{-1}=\left( \begin{matrix}
1-uv & u \\
-v & 1-vu \end{matrix} \right)
$$
verifies $W^{-1}\bigl( \begin{smallmatrix} \rho(a) & 0 \\
                                                    0 & 0 \end{smallmatrix} \bigr) W = \bigl( \begin{smallmatrix} 0 & 0 \\
                           0 & \rho'(a) \end{smallmatrix} \bigr)$, which shows that the homomorphisms $\rho$ and $\rho'$ are stably conjugate.\\
In order to get a quasihomomorphism from $\Ac_P$ to $\Ac\cp \Gamma$, we need a $K$-cycle for the Fr\'echet algebra $\cinf(M)$ (see Example \ref{ek}). By a standard procedure \cite{C0, C1}, such a $K$-cycle $D$ may be constructed from an elliptic pseudodifferential operator or a Toeplitz operator over $M$. We shall suppose that $D$ is of parity $p$ mod 2, and of summability degree $p+1$ (even case) or $p$ (odd case). Hence (see Example \ref{ek}) in the even case one has an infinite-dimensional separable Hilbert space $H$ with two continuous representations $\cinf(M)\rightrightarrows \Lc=\Lc(H)$ which agree modulo the Schatten ideal $\Ic=\Lc^{p+1}(H)$, whereas in the odd case the algebra $\cinf(M)$ is represented in the matrix algebra $\bigl( \begin{smallmatrix} 
\Lc & \Ic \\
\Ic & \Lc \end{smallmatrix} \bigr) $ with $\Ic=\Lc^p(H)$. Therefore, upon choosing an isomorphism $H\cong H\hotimes \cc^m$ the composition of $\rho:\Ac_P\to M_m(\cinf(M)\hotimes (\Ac\cp \Gamma))$ with the Hilbert space representation induced by the K-cycle $D$ leads to a quasihomomorphism of parity $p$ mod 2
$$
\rho_D : \Ac_P\to \Ec^s\triangleright \Ic^s\hotimes (\Ac\cp \Gamma)\ ,
$$
with intermediate algebra $\Ec=\Lc\hotimes(\Ac\cp \Gamma)$ (or $(\Lc\ltimes \Ic)\hotimes(\Ac\cp \Gamma)$, see Example \ref{ebiv}). Note that $\Lc$ and $\Ic$ may be replaced by other suitable operator algebras, if needed. From the discussion above we see that $\rho_D$ depends only on $D$ up to conjugation by an invertible element $W\in \Ec^s_+$. Taking into account the Chern characters in negative and periodic cyclic homology, the Riemann-Roch-Grothendieck Theorem \ref{trr} thus yields cube diagrams of the following kind:

\begin{corollary}
A $K$-cycle $D$ over $\cinf(M)$ as above yields for any integer $n\in\zz$ a commutative diagram
$$
\vcenter{\xymatrix@!0@=3.7pc{ & \Kt_n(\Ic\Ac_P) \ar[rr] \ar'[d][dd] & & \Kt_{n-p}(\Ic(\Ac\cp \Gamma)) \ar[dd] \\
MK^{\Ic}_n(\Ac_P) \ar[ur] \ar[rr] \ar[dd] & & MK^{\Ic}_{n-p}(\Ac\cp \Gamma) \ar[ur] \ar[dd] \\
 & HP_n(\Ac_P) \ar'[r][rr] & & HP_{n-p}(\Ac\cp \Gamma) \\
HN_n(\Ac_P) \ar[ur] \ar[rr] & & HN_{n-p}(\Ac\cp \Gamma) \ar[ur] & }}
$$
where the horizontal arrows are induced by the quasihomomorphism $\rho_D : \Ac_P\to \Ec^s\triangleright \Ic^s\hotimes (\Ac\cp \Gamma)$. 
\end{corollary}
The background square describes the topological side of the Riemann-Roch-Grothendieck theorem, namely the compatibility between the push-forward in topological $K$-theory and the bivariant Chern character in periodic cyclic homology. One may choose $D$ as a representative of the fundamental class in the $K$-homology of $M$. If moreover $M$ is a model for the classifying space $B\Gamma$, one may choose $P$ as the universal bundle $E\Gamma$. For torsion-free groups $\Gamma$ the morphism $\Kt_n(\Ic\Ac_P)\to \Kt_{n-p}(\Ic(\Ac\cp \Gamma))$ thus obtained is related to the Baum-Connes assembly map \cite{BC} and exhausts many (in some cases, all the) interesting topological $K$-theory classes of $\Ac\cp \Gamma$. \\
The foreground square provides a lifting of the topological situation at the level of multiplicative $K$-theory and negative cyclic homology, i.e. secondary characteristic classes. Hence a part of $MK_*^{\Ic}(\Ac\cp \Gamma)$ may be obtained by direct images of multiplicative $K$-theory over $\Ac_P$. Note that in contrast to the topological situation, the push-forward map in multiplicative $K$-theory does \emph{not} exhaust all the interesting classes over $\Ac\cp \Gamma$.\\

Let us now deal with the case $\Ac=\cinf(N)$, for a compact smooth Riemannian manifold $N$, endowed with a left action of $\Gamma$ by diffeomorphisms. We provide $\Ac$ with its usual Fr\'echet topology, and condition (\ref{cond}) forces the $\Gamma$-action be "almost isometric" on $N$. The crossed product $\Ac\cp \Gamma$ is then isomorphic to a certain convolution algebra of functions over the smooth \'etale groupoid $\Gamma\ltimes N$, describing a highly noncommutative space when the action of $\Gamma$ is not proper. The commutative algebra $\Ac_P$ is the subalgebra of smooth $\Gamma$-invariant functions $a\in \cinf(P\times N)$, $a(\gamma\cdot x, \gamma\cdot y)=a(x,y)$ for any $(x,y)\in P\times N$, and is thus isomorphic to the algebra of smooth functions over the (compact) quotient manifold $Q=\Gamma\backslash(P\times N)$. \\
The problem is therefore reduced to the computation of secondary invariants for the classical space $Q$. The cyclic homology of $\Ac_P=\cinf(Q)$ has been determined by Connes \cite{C0} and is computable from the de Rham complex of differential forms over $Q$. We will see that the multiplicative $K$-theory $MK^{\Ic}_*(\Ac_P)$ is closely related (though not isomorphic) to Deligne cohomology. We first recall some definitions. Let $\Om^n(Q)$ denote the space of complex, smooth differential $n$-forms over $Q$, $d$ the de Rham coboundary, $\Zdr^n(Q)=\ker(d:\Om^n\to\Om^{n+1})$ the space of closed $n$-forms and $\Bdr^n(Q)=\im(d:\Om^{n-1}\to\Om^n)$ the space of exact $n$-forms. By de Rham's theorem, the de Rham cohomology $\Hdr^n(Q)=\Zdr^n(Q)/\Bdr^n(Q)$ is isomorphic to the \v{C}ech cohomology of $Q$ with complex coefficients $\Hch^n(Q;\cc)$. For any half-integer $q$ we define the additive group $\zz(q):=(2\pi i)^q\zz\subset \cc$ (the square root of $2\pi i$ must be chosen consistently with the Chern character on $\Kt_1$). Let $\Omb^k$ denote the sheaf of differential $k$-forms over $Q$ and consider for any $n\in\nn$ the complex of sheaves
\be
0\longrightarrow \zzb(n/2) \longrightarrow \Omb^0 \stackrel{d}{\longrightarrow} \Omb^1 \stackrel{d}{\longrightarrow} \ldots \stackrel{d}{\longrightarrow} \Omb^{n-1} \longrightarrow 0 \label{sheaf}
\ee
where the constant sheaf $\zzb(n/2)$ sits in degree $0$ and $\Omb^k$ in degree $k+1$. The map $\zzb(n/2)\to \Omb^0$ is induced by the natural inclusion of constant functions into complex-valued functions. By definition the (smooth) Deligne cohomology $H^n_{\Dc}(Q;\zz(n/2))$ is the hyperhomology of (\ref{sheaf}) in degree $n$. The natural projection onto the constant sheaf $\zzb(n/2)$ yields a well-defined map from $H^n_{\Dc}(Q;\zz(n/2))$ to the \v{C}ech cohomology with integral coefficients $\check{H}^n(Q;\zz(n/2))$. On the other extreme, the de Rham coboundary $d:\Omb^{n-1}\to \Omb^n$ sends a Deligne $n$-cocycle to a globally defined closed $n$-form over $Q$, called its \emph{curvature}, which only depends on the Deligne cohomology class. It follows from the definitions that the image of the curvature in de Rham cohomology coincides with the complexification of the \v{C}ech cohomology class of the Deligne cocycle. One thus gets a commutative diagram in any degree $n$
\be
\vcenter{\xymatrix{
H^n_{\Dc}(Q;\zz(n/2)) \ar[r] \ar[d]_d  & \Hch^n(Q;\zz(n/2)) \ar[d]^{\otimes\cc}  \\
\Zdr^n(Q) \ar[r]  & \Hdr^n(Q) }}  \label{del}
\ee
This has to be compared with the commutative square involving the multiplicative and topological $K$-theories of the algebra $\Ac_P=\cinf(Q)$, with their Chern characters:
\be
\vcenter{\xymatrix{
MK^{\Ic}_n(\Ac_P) \ar[r] \ar[d] & \Kt_n(\Ic\Ac_P) \ar[d]  \\
HN_n(\Ac_P) \ar[r]  & HP_n(\Ac_P) }}  \label{square}
\ee
In fact one can construct, at least in low degrees $n$, an explicit transformation from Deligne cohomology to multiplicative $K$-theory, and the curvature morphism captures the lowest degree part of the negative Chern character. Let us explain this with some details. Firstly, it is well-known that the bottom line of (\ref{del}) is included as a direct summand in the bottom line of (\ref{square}). Since we deal essentially with the $X$-complex description of cyclic homology (section \ref{scy}), we recall how the latter is related to the de Rham cohomology of $Q$. Choose the universal free extension $0\to J\Ac_P \to T\Ac_P \to \Ac_P\to 0$. The cyclic homology of $\Ac_P$ is computed by the $X$-complex $X(\Th\Ac_P)$ of the pro-algebra $\Th\Ac_P=\varprojlim_n T\Ac_P/(J\Ac_P)^n$, together with its filtration by the subcomplexes $F^n\Xh(T\Ac_P,J\Ac_P)$. As a pro-vector space, $X(\Th\Ac_P)$ is isomorphic to the completed space of noncommutative differential forms $\Omh\Ac_P$, and the $J\Ac_P$-adic filtration coincides with the Hodge filtration $F^n\Omh\Ac_P$. A canonical chain map $X(\Th\Ac_P)\to \Om^*(Q)$ is given by projecting noncommutative differential forms to commutative ones, up to a rescaling: 
$$
\Om^n\Ac_P \ni a_0da_1\ldots da_n \to \la_n a_0da_1\ldots da_n \in \Om^n(Q)
$$
with $\la_n=(-)^k\frac{k!}{(2k)!}$ if $n=2k$ and $\la_n=(-)^k\frac{k!}{(2k+1)!}$ if $n=2k+1$. These factors are fixed in order to get exactly a chain map. Clearly it sends the Hodge filtration of $\Omh\Ac_P$ onto the natural filtration by degree on $\Om^*(Q)$. The following proposition is a reformulation of Connes' version of the Hochschild-Kostant-Rosenberg theorem \cite{C0}.

\begin{proposition}
The chain map $X(\Th\Ac_P)\to \Om^*(Q)$ is a homotopy equivalence compatible with the filtrations. Hence follow the isomorphisms
\beq
HP_n(\Ac_P) &=& \bigoplus_{k\in\zz} \Hdr^{n+2k}(Q)\ ,\non\\
HC_n(\Ac_P) &=& \frac{\Om^n(Q)}{\Bdr^n(Q)}\oplus \bigoplus_{k< 0} \Hdr^{n+2k}(Q)\ ,\\
HN_n(\Ac_P) &=& \Zdr^n(Q) \oplus \bigoplus_{k>0} \Hdr^{n+2k}(Q)\ .\non
\eeq
\end{proposition}

\noindent Of course the injections $\Zdr^n(Q)\to HN_n(\Ac_P)$ and $\Hdr^n(Q)\to HP_n(\Ac_P)$ are compatible with the forgetful maps $\Zdr^n(Q)\to \Hdr^n(Q)$ and $HN_n(\Ac_P)\to HP_n(\Ac_P)$. It is useful to calculate the image of the Chern character of idempotents and invertibles under the chain map $X(\Th\Ac_P)\to\Om^*(Q)$. Let $e\in M_2(\Kc\hotimes\Ac_P)^+$ be an idempotent such that $e-p_0\in M_2(\Kc\hotimes\Ac_P)$, and let $\eh\in  M_2(\Kc\hotimes\Th\Ac_P)^+$ be its canonical lift. The image of the Chern character $\ch_0(\eh)=\Tr(\eh-p_0)\in X(\Th\Ac_P)$ is the differential form of even degree
$$
\chdr(e)=\Tr(e-p_0) +\sum_{k\geq 1}\frac{(-)^k}{k!}\,\Tr((e-\frac{1}{2})(dede)^k)\ \in \Om^+(Q)\ .
$$
Let $g\in (\Kc\hotimes\Ac_P)^+$ be an invertible such that $g-1\in\Kc\hotimes\Ac_P$, and let $\gh\in (\Kc\hotimes\Th\Ac_P)^+$ be its canonical lift. The image of the Chern character $\ch_1(\gh)=\frac{1}{\sqrt{2\pi i}}\nat \gh^{-1}\dd\gh\in X(\Th\Ac_P)$ is the differential form of odd degree
$$
\chdr(g)=\frac{1}{\sqrt{2\pi i}} \sum_{k\geq 0}(-)^k \frac{k!}{(2k+1)!}\, \Tr(g^{-1}dg(dg^{-1}dg)^k)\ \in \Om^-(Q)\ .
$$
We now construct explicit morphisms $H_{\Dc}^n(Q;\zz(n/2))\to MK_n^{\Ic}(\Ac_P)$ in degrees $n=0,1,2$. In fact the ideal $\Ic$ is irrelevant and the previous morphisms factor through the multiplicative group $MK_n(\Ac_P):= MK^{\cc}_n(\Ac_P)$ associated to the $1$-summable algebra $\cc$. Then choosing any rank one injection $\cc\to\Ic$ induces a unique map $MK_n(\Ac_P)\to MK^{\Ic}_n(\Ac_P)$ by virtue of Lemma \ref{linv}. \\

\noindent $n=0$: Then $\zz(0)=\zz$ and the complex $0\to \zzb(0) \to 0$ calculates the \v{C}ech cohomology of $Q$ with coefficients in $\zz$. Hence $H_{\Dc}^0(Q;\zz(0))=\Hch^0(Q;\zz)$ is the additive group of $\zz$-valued locally constant functions over $Q$. The map
\be
H_{\Dc}^0(Q;\zz(0)) \to MK_0(\Ac_P)\cong \Kt_0(\Ac_P)
\ee
associates to such a function $f$ the $K$-theory class of the trivial complex vector bundle of rank $f$ over $Q$. \\

\noindent $n=1$: Then $\zz(1/2)=\sqrt{2\pi i}\, \zz$ and $H_{\Dc}^1(Q;\zz(1/2))$ is the hyperhomology in degree 1 of the complex $0\to \zzb(1/2)\to \Omb^0\to 0$. Choose a good covering $(U_i)$ of $Q$. A Deligne $1$-cocycle is given by a collection $(f_i, n_{ij})$ of smooth functions $f_i:U_i\to \cc$ and locally constant functions $n_{ij}:U_i\cap U_j\to \zz$ related by the descent equations
$$
f_i-f_j=\sqrt{2\pi i}\, n_{ij}\ \mbox{over}\ U_i\cap U_j\ ,\quad n_{jk}-n_{ik}+n_{ij}=0 \ \mbox{over}\ U_i\cap U_j\cap U_k\ .
$$
The cocycle is trivial if the $f_i$'s are $\zz(1/2)$-valued. Taking the exponentials $g_i=\exp(\sqrt{2\pi i} f_i)$ one gets invertible smooth functions which agree on the overlaps $U_i\cap U_j$, hence define a global invertible function $g$ over $Q$. The latter is equal to 1 exactly when the cocycle $(f_i,n_{ij})$ is trivial. Hence $H_{\Dc}^1(Q;\zz(1/2))$ is the multiplicative group $\cinf(Q)^{\times}$ of complex-valued invertible functions over $Q$. On the other hand, the elements of $MK_1(\Ac_P)$ are represented by pairs $(\gh,\te)$ of an invertible $\gh\in (\Kc\hotimes\Th\Ac_P)^+$ and a cochain $\te\in X_0(T\Ac_P,J\Ac_P)\cong\cinf(Q)$. We get a map
\be
H_{\Dc}^1(Q;\zz(1/2))\cong\cinf(Q)^{\times}\to MK_1(\Ac_P)
\ee
by sending an invertible $g\in \cinf(Q)^{\times}$ to the multiplicative $K$-theory class of $(\gh,0)$, with $\gh$ the canonical lift of $g$ (to be precise one should replace $g$ by $1+(g-1_{\Ac_P}) \in (\Ac_P)^+$). This map identifies the curvature morphism $H_{\Dc}^1(Q;\zz(1/2))\to \Zdr^1(Q)$ with the lowest degree part of the negative Chern character $\ch_1:MK_1(\Ac_P)\to HN_1(\Ac_P)$. Indeed the curvature of an element $g\in \cinf(Q)^{\times}$ is by definition the closed one-form 
$$
df_i=\frac{1}{\sqrt{2\pi i}} \, g_i^{-1}dg_i= \frac{1}{\sqrt{2\pi i}} \, g^{-1}dg \qquad \forall i\ ,
$$
globally defined over $Q$. But this coincides with the $\Zdr^1(Q)$-component of the negative Chern character $\ch_1(\gh,0)$.\\

\noindent $n=2$: Then $\zz(1)=2\pi i\, \zz$ and $H_{\Dc}^2(Q;\zz(1))$ is the hyperhomology in degree 2 of the complex $0\to \zzb(1)\to \Omb^0\to \Omb^1 \to 0$. A Deligne cocycle relative to the finite good covering $(U_i)$, $i=1,\ldots,m$ is a collection $(A_i,f_{ij},n_{ijk})$ of one-forms $A_i$ over $U_i$, smooth functions $f_{ij}: U_i\cap U_j\to \cc$ and locally constant functions $n_{ijk}:U_i\cap U_j\cap U_k\to \zz$, subject to the descent equations
$$
A_i-A_j=df_{ij}\ ,\quad -f_{jk}+f_{ik}-f_{ij}=2\pi i\, n_{ijk}\ ,\quad n_{jkl}-n_{ikl}+n_{ijl}-n_{ijk}=0\ .
$$
Equivalently, passing to the exponentials $g_{ij}=\exp f_{ij}$ a cocycle is a collection $(A_i,g_{ij})$ such that $A_i-A_j=g_{ij}^{-1}dg_{ij}$ and $g_{ij}g_{jk}=g_{ik}$. Two cocycles $(A_i,g_{ij})$ and $(A'_i,g'_{ij})$ are cohomologous iff there exists a collection of smooth invertible functions (gauge transformations) $\al_i:U_i\to \cc^{\times}$ such that 
$$
A'_i=A_i+\al_i^{-1}d\al_i\ ,\qquad g'_{ij}=\al_ig_{ij}\al_j^{-1}\ .
$$
One sees that a Deligne cohomology class is nothing else but a complex line bundle over $Q$, described by the smooth transition functions $g_{ij}:U_i\cap U_j\to \cc^{\times}$, together with a connection given locally by the one-forms $A_i$, up to gauge transformation. Hence
$$
H_{\Dc}^2(Q;\zz(1))=\{\mbox{isomorphism classes of complex line bundles with connection}\}
$$
The group law is the tensor product of line bundles with connections. The curvature morphism $H_{\Dc}^2(Q;\zz(1))\to \Zdr^2(Q)$ maps a cocycle $(A_i,g_{ij})$ to the globally defined closed two-form $dA_i$ $\forall i$, i.e. the curvature of the connection of the corresponding line bundle. The construction of the morphism from Deligne cohomology to multiplicative $K$-theory requires to fix a partition of unity $(c_i)$ relative to the finite covering: $\supp c_i\subset U_i$ and $\sum_i c_i(x)^2=1$ $\forall x\in Q$. Given a Deligne cocycle $(A_i,g_{ij})$, we construct the idempotent $e_+ \in M_m(\cinf(Q))$ of rank 1 whose matrix elements are the functions
$$
(e_+)_{ij}=c_ig_{ij}c_j\ \in \cinfc(U_i\cap U_j)\ ,
$$
and let $e_-=1_{\Ac_P}$ be the unit of $\cinf(Q)$ (the constant function $1$ over $Q$). Define $e\in M_2(\Ac_P)^+$ as the idempotent matrix $\bigl( \begin{smallmatrix} 1-e_- & 0 \\ 0 & e_+  \end{smallmatrix} \bigr)$. The 0th degree of the Chern character $\chdr(e)$ is $\Tr(e-p_0)=\Tr(e_+)-1_{\Ac_P}=0$, so that the class of $e$ in $\Kt_0(\Ac_P)$ is a virtual bundle of rank 0. To get a multiplicative $K$-theory class $(\eh,\te)\in MK_2(\Ac_P)$, we must adjoin to the canonical idempotent lift $\eh\in M_2(\Th\Ac_P)^+$ an odd cochain $\te\in X_1(T\Ac_P,J\Ac_P)\cong X(\Ac_P)$. Since $\Ac_P$ is the commutative algebra of smooth functions over a compact manifold, its $X$-complex reduces to the de Rham complex of $Q$ truncated in degrees $\geq 2$, $X(\Ac_P):\cinf(Q)\stackrel{d}{\to} \Om^1(Q)$. The fact that $e$ is of virtual rank zero insures that any choice of one-form $\te\in \Om^1(Q)$ satisfies the correct transgression relation $\ch_0(\eh)=\bb\te=0$ in the complex $X(\Ac_P)$. We set
$$
\te= -\sum_i c_i^2 A_i \ \in \Om^1(Q)\ .
$$

\begin{lemma}
The assignement $(A_i,g_{ij})\mapsto (\eh,\te)$ yields a well-defined morphism
\be
H_{\Dc}^2(Q;\zz(1))\to MK_2(\Ac_P)\ ,
\ee
and the curvature of $(A_i,g_{ij})$ corresponds to the $\Zdr^2(Q)$ component of the negative Chern character $\ch_2(\eh,\te)\in HN_2(\Ac_P)$.
\end{lemma}
{\it Proof:} We have to show that the multiplicative $K$-theory class of $(\eh,\te)$ only depends on the Deligne cohomology class of $(A_i,g_{ij})$. Thus let $(A'_i,g'_{ij})$ be another representative, with $A'_i=A_i+\al_i^{-1}d\al_i$ and $g'_{ij}=\al_i g_{ij}\al_j^{-1}$. This yields a new pair $(\eh',\te')$ with $(e'_+)_{ij}=c_ig'_{ij}c_j$ and $\te'=-\sum_i c_i^2A'_i$. We show that $(\eh,\te)$ and $(\eh',\te')$ represent the same class in $MK_2(\Ac_P)$ by using the following general fact: if $e$ and $e'=W^{-1}eW$ are conjugate by an invertible $W$, then $(\eh,\te)$ is equivalent to $(\eh', \te+\cs_1(f))$, where $f$ is the idempotent interpolation between the matrices $\bigl( \begin{smallmatrix} e & 0 \\ 0 & p_0  \end{smallmatrix} \bigr)$ and $\bigl( \begin{smallmatrix} W^{-1}eW & 0 \\ 0 & p_0  \end{smallmatrix} \bigr)$ constructed as in Lemma \ref{linv}. One has
$$
\cs_1(f) \equiv \Tr(W^{-1} (e-p_0)d W) \mod d\quad \mbox{in}\ \Om^1(Q)\ ,
$$
so that finally $(\eh, \te)$ is equivalent to $(\eh',\te +\Tr(W^{-1} (e-p_0)d W) )$. In the present situation $g'_{ij}=\al_i g_{ij}\al_j^{-1}$, hence the idempotent $e'_+$ is conjugate to $e_+$ via the diagonal matrix $W_+=$ diag $(\al_1^{-1}\,\ldots,\al_m^{-1})$, and of course $e'_-=e_-=1_{\Ac_P}$ so one can choose $W_-=1$. One calculates
\beq
\Tr(W^{-1} (e-p_0)d W) &=& \Tr(W_+^{-1} e_+d W_+)\ =\ \sum_i\al_ic_i^2d(\al_i^{-1})\non\\
&=& -\sum_ic_i^2(A'_i-A_i)\ =\ \te'-\te\ ,\non
\eeq
hence $(\eh,\te)$ and $(\eh',\te')$ represent the same multiplicative $K$-theory class. \\
Now we leave the cocycle $(A_i,g_{ij})$ fixed and change the partition of unity $(c_i)$ to $(c'_i)$, $\sum_i(c'_i)^2=1$, whence a new idempotent $(e'_+)_{ij}=c'_ig_{ij}c'_j$ and a new cochain $\te'=-\sum_i(c'_i)^2A_i$. Introduce the matrices $u_{ij}=c_ig_{ij}c'_j$ and $v_{ij}=c'_ig_{ij}c_j$. Then one has $e_+=uv$, $e'_+=vu$, and the invertible matrix $W_+= \bigl( \begin{smallmatrix} 1-uv & -u \\ v & 1-vu \end{smallmatrix} \bigr) $ stably conjugates $e_+$ and $e'_+$ in the sense that $W_+^{-1} \bigl( \begin{smallmatrix} e_+ & 0 \\ 0 & 0 \end{smallmatrix} \bigr) W_+=\bigl( \begin{smallmatrix} 0 & 0 \\ 0 & e'_+ \end{smallmatrix} \bigr) $. A direct computation yields
$$
\Tr(W_+^{-1}\bigl( \begin{smallmatrix} e_+ & 0 \\ 0 & 0 \end{smallmatrix} \bigr) dW_+) = \sum_{i,j}(c'_i)^2c_j^2 g_{ij}dg_{ji}=\sum_{i,j}(c'_i)^2c_j^2(A_j-A_i)=\te'-\te\ ,
$$
hence the multiplicative $K$-theory class of $(\eh,\te)$ does not depend on the choice of partition of unity. A similar argument shows that it does not depend on the good covering $(U_i)$.\\
It remains to calulate the lowest component of the negative Chern character. By definition $\ch_2(\eh,\te)$ is the cycle of even degree in $X(\Th\Ac_P)$ given by $\Tr(\eh-p_0)-\bb\tilde{\te}$, where $\tilde{\te}$ is an arbitrary lifting of $\te\in X(\Ac_P)$. Thus, the image of $\ch_2(\eh,\te)$ under the chain map $X(\Th\Ac_P)\to \Om^*(Q)$ has a component of degree two given by 
$$
-\Tr((e-\frac{1}{2})dede)-d\te = -\sum_i d(c_i^2)A_i +d\sum_i c_i^2A_i=\sum_ic_i^2dA_i=dA_i\ ,
$$
and coincides with the curvature of the line bundle $(A_i,g_{ij})$. \cqfd\\

\noindent If one forgets the connection $A_i$, the morphism $H_{\Dc}^2(Q;\zz(1))\to MK_2(\Ac_P)$ just reduces to the elementary map $\Hch^2(Q;\zz(1))\to \Kt_2(\Ac_P)$ which associates to an isomorphism class of line bundles over $Q$ its topological $K$-theory class. \\

\begin{example}\textup{The simplest non-trivial example is provided by the celebrated non-commutative torus \cite{C}. Here $\Ac$ is the algebra of smooth functions over the circle $N=S^1=\zz\backslash\rr$. Conventionnally we parametrize the points of $N$ by the variable $y$. The group $\Gamma=\zz$ acts on $N$ by rotations of angle $\al\in \rr$: 
$$
y\mapsto y+n\al\qquad \forall\ y\in N\ ,\ n\in\zz\ .
$$
When $\zz$ is provided with its natural distance, the crossed product $\Ac\cp \zz$ is isomorphic to the algebra $\Ac_{\al}$ of the non-commutative torus, presented for example in $\cite{C}$ by generators and relations: let $V_1\in \Ac$ be the function $V_1(y)=e^{2\pi i y}$ over $N$ and $V_2=1^*\in\cc\zz$ be the element corresponding to the generator $1\in\zz$. Then $V_1$ and $V_2$ are invertible elements of $\Ac_{\al}$ and fulfill the noncommutativity relation
\be
V_2V_1=e^{2\pi i \al}V_1V_2\ .
\ee
Moreover any element of $\Ac_{\al}$ is a power series $\sum_{(n_1,n_2)\in\zz^2}a_{n_1n_2}V_1^{n_1}V_2^{n_2}$ with coefficients $a_{n_1n_2}\in\cc$ of rapid decay. For $\al\in \qq$ this algebra is Morita equivalent to the smooth functions over an ordinary (commutative) $2$-torus, and its multiplicative $K$-theory in any degree turns out to be completely determined by Deligne cohomology in this case. The situation is more interesting for $\al\notin\qq$. Following the discussion above we introduce the universal principal $\zz$-bundle $P=E\zz=\rr$ over the classifying space $M=B\zz=\zz\backslash\rr$. Conventionnally we parametrize the points of $P$ by the variable $x$. Thus, $\Ac_P=\cinf(P;\Ac)^{\zz}$ is the mapping torus algebra
$$
\Ac_P=\{ a\in \cinf(P\times N)\ |\ a(x+1,y+\al)=a(x,y)\ ,\ \forall\ x\in P\ ,\ y\in N\}\ .
$$
Equivalently it is the algebra of smooth functions over the commutative 2-torus $Q=\zz\backslash(P\times N)$, quotient of $\rr^2$ by the lattice generated by the vectors $(1,\al)$ and $(0,1)$. Now to get a quasihomomorphism from $\Ac_P$ to $\Ac_{\al}$ we need a $K$-cycle $D$ for the circle manifold $M$. Let $H=L^2(M)$ be the Hilbert space of square-integrable complex-valued functions. The algebra $\cinf(M)$ is represented in the algebra of bounded operators $\Lc(H)$ by pointwise multiplication. $D$ will be represented by the Toeplitz operator in $\Lc(H)$ which projects $H$ onto the Hardy space $H_+\subset H$:
$$
D(e^{2\pi i\, nx})=\left\{ \begin{array}{lcc} 
e^{2\pi i\, nx} & \mbox{if} & n\geq 0 \\
0 & \mbox{if} & n< 0 \end{array} \right.
$$ 
One thus obtains a polarization of the Hilbert space $H=H_+\oplus H_-$, with $H_-$ the kernel of $D$. In the Fourier basis $e^{2\pi i\, nx}$ of $H$, the representation $\cinf(M)\to \Lc(H)$ is easily seen to factor through the matrix subalgebra
$$
\left( \begin{matrix} 
\Tc & \Kc \\
\Kc & \Tc \end{matrix} \right)\subset \left( \begin{matrix} 
\Lc(H_+) & \Lc(H_-,H_+) \\
\Lc(H_+,H_-) & \Lc(H_-) \end{matrix} \right) =\Lc(H)\ ,
$$
where $\Tc$ is the \emph{smooth Toeplitz algebra} (the elementary non-trivial extension of $\cinf(S^1)$ by th ealgebra $\Kc$ of smooth compact operators, see \cite{Cu1} ). The induced quasihomomorphism $\rho_D:\Ac_P\to \Ec^s\triangleright \Kc^s\hotimes\Ac_{\al}$, with $\Ec=\Tc\hotimes\Ac_{\al}$, is therefore $1$-summable and of odd parity. Theorem \ref{trr} yields a graded-commutative diagram (remark that $MK^{\Kc}_n\cong MK_n$)
\be
\vcenter{\xymatrix{
\Kt_{n+1}(\Ac_P) \ar[r] \ar[d] & HC_{n-1}(\Ac_P) \ar[r] \ar[d]^{\ch^1(\rho_D)} & MK_n(\Ac_P)  \ar[r] \ar[d]  & \Kt_n(\Ac_P)  \ar[d]  \\
\Kt_{n}(\Ac_{\al}) \ar[r]  & HC_{n-2}(\Ac_{\al}) \ar[r]  & MK_{n-1}(\Ac_{\al})  \ar[r]  & \Kt_{n-1}(\Ac_{\al}) }} \label{toto}
\ee
in any degree $n\in\zz$. The group $\Kt_n(\Ac_P)$ is isomorphic to the topological $K$-theory of the 2-torus $Q$. Hence in even degree, $\Kt_0(\Ac_P)=\zz\oplus\zz$ is generated by the trivial line bundle over $Q$ together with the Bott class, whereas in odd degree $\Kt_1(\Ac_P)=\zz\oplus\zz$ is generated by the invertible functions $g_1(x,y)=e^{2\pi i\, x}$ and $g_2(x,y)=e^{2\pi i(\al x-y)}$. The pushforward map $\Kt_n(\Ac_P) \to \Kt_{n-1}(\Ac_{\al})$ is known to be an isomorphism (Baum-Connes). In particular the Bott class and the trivial line bundle over $Q$ are mapped respectively to the classes of the invertible elements $V_1$ and $V_2$ in $\Kt_1(\Ac_{\al})$. For multiplicative $K$-theory the situation is more involved. In degree $n=1$ the map 
$$
MK_1(\Ac_P)\cong H^1_{\Dc}(Q;\zz(1/2))\cong \cinf(Q)^{\times}\to MK_0(\Ac_{\al})\cong \Kt_0(\Ac_{\al})
$$
simply factors through the topological $K$-theory group $\Kt_1(\Ac_P)$.  In degree $n=2$ one still has an isomorphism $MK_2(\Ac_P)\cong H^2_{\Dc}(Q;\zz(1))$, and (\ref{toto}) amounts to
\beq
\vcenter{\xymatrix{
\Kt_1(\Ac_P) \ar[r]^{\ch_1} \ar[d] & \frac{\Om^1(Q)}{d\Om^0(Q)} \ar[r] \ar[d]^{\ch^1(\rho_D)} & H^2_{\Dc}(Q;\zz(1))  \ar[r] \ar[d]  & \Kt_0(\Ac_P)  \ar[d] & \\
\Kt_0(\Ac_{\al}) \ar[r]  & HC_{0}(\Ac_{\al}) \ar[r]  & MK_{1}(\Ac_{\al})  \ar[r]  & \Kt_{1}(\Ac_{\al}) \ar[r] & 0 }} \label{toto2}
\eeq
The map $\Om^1/d\Om^0 \to H^2_{\Dc}(Q;\zz(1))$ associates to a one-form $\frac{1}{\sqrt{2\pi i}}\,A$ over $Q$ the isomorphism class of the trivial line bundle with connection $-A$, while the range of $H^2_{\Dc}(Q;\zz(1))\to \Kt_0(\Ac_P)$ is generated by the Bott class. For generic values of $\al\notin \qq$ the commutator subspace $[\Ac_{\al},\Ac_{\al}]$ may not be closed in $\Ac_{\al}$, therefore the quotient $HC_0(\Ac_{\al})=HH_0(\Ac_{\al})= \Ac_{\al}/[\Ac_{\al},\Ac_{\al}]$ may not be separated. However the quotient $\overline{HC}_0(\Ac_{\al})$ by the \emph{closure} of the commutator subspace turns out to be isomorphic to $\cc$, via the canonical trace 
$$
\Ac_{\al}\to \cc\ , \qquad V_1^{n_1}V_2^{n_2} \mapsto \left\{ \begin{array}{cl}
1 & \mbox{if} \ n_1=n_2=0\ , \\
0 & \mbox{otherwise}\ . \end{array} \right. 
$$
With these identifications the evaluation of the Chern character $\ch^1(\rho_D): \Om^1/d\Om^0 \to \overline{HC}_0(\Ac_{\al})\cong \cc$ on a one-form $A=A_xdx+A_ydy$ is easily performed and one finds
$$
\ch^1(\rho_D)\left( \frac{A}{\sqrt{2\pi i}} \right)= \frac{1}{2\pi i} \int_0^1dy\int_0^1 dx \, A_x(x,y) \ .
$$
In particular $\ch^1(\rho_D)\cdot\ch_1(g_1)=1$ and $\ch^1(\rho_D)\cdot\ch_1(g_2)=\al$, and one recovers the well-known fact (\cite{C0}) that the image of $\Kt_0(\Ac_{\al})$ in $\overline{HC}_0(\Ac_{\al})$ is the subgroup $\zz+\al\zz \subset \cc$. \\
We may analogously define a new multiplicative $K$-theory group $\overline{MK}_1(\Ac_{\al})$ whose elements are represented by pairs $(\gh,\te)$ with $\te\in \cc$ instead of $\te \in \Ac_{\al}/[\Ac_{\al},\Ac_{\al}]$. Because $\Kt_1(\Ac_{\al})$ is generated by the invertibles $V_1$ and $V_2$ any class in $\overline{MK}_1(\Ac_{\te})$ is represented by a pair $(V_1^{n_1}V_2^{n_2},\te)$ for some integers $n_1,n_2$ and a complex number $\te$. Using a homotopy one shows that this pair is equivalent to $(e^{-\sqrt{2\pi i}\te}V_1^{n_1}V_2^{n_2},0)$, and by exactness $\overline{MK}_1(\Ac_{\al})$ is the quotient of the multiplicative group $\cc^{\times}\langle V_1\rangle \langle V_2 \rangle\subset \mbox{GL}_1(\Ac_{\al})$ by its commutator subgroup $\langle e^{2\pi i \al}\rangle\subset \cc^{\times}$, or equivalently the abelianization
\be
\overline{MK}_1(\Ac_{\al}) = (\cc^{\times}\langle V_1\rangle \langle V_2 \rangle)_{\mathrm{ab}}\ .
\ee
Since the Bott class of $K_0(\Ac_P)$ is sent to $[V_1]\in \Kt_1(\Ac_{\al})$, one sees that the range of $H^2_{\Dc}(Q;\zz(1))\to \overline{MK}_1(\Ac_{\al})$ coincides with the subgroup $\cc^{\times}\langle V_1 \rangle$. }
\end{example}

\end{document}